\documentclass[reqno,a4paper,oneside]{amsart}

\usepackage{amsmath, amsthm, amscd, amsfonts, amssymb, graphicx, color,mathrsfs,mathtools,enumerate}
\usepackage{lineno}
\usepackage{dsfont}
\modulolinenumbers[5]

\usepackage{bm}
\usepackage{graphicx}%
\usepackage{multirow}%
\usepackage{amsmath,amssymb,amsfonts}%
\usepackage{amsthm}%
\usepackage{mathrsfs}%
\usepackage[title]{appendix}%
\usepackage{xcolor}%
\usepackage{textcomp}%
\usepackage{manyfoot}%
\usepackage{booktabs}%
\usepackage{algorithm}%
\usepackage{algorithmicx}%
\usepackage{algpseudocode}%
\usepackage{listings}%
\usepackage{refcheck}
\newtheorem{thm}{Theorem}[section]
\newtheorem{lemma}[thm]{Lemma}
\newtheorem{proposition}[thm]{Proposition}
\newtheorem{corollary}[thm]{Corollary}

\newtheorem{hyp1}[thm]{Hypotheses}

\newtheorem{rmk}[thm]{Remark}
\newtheorem{remark}[thm]{Remark}

\newtheorem{defn}[thm]{Definition}

\newcommand{\N}{\mathbb N}
\newcommand{\R}{\mathbb R}
\newcommand{\X}{\mathcal{X}}
\newcommand{\K}{\mathcal{K}}
\newcommand{\C}{\mathcal{C}}
\newcommand{\Id}{{\operatorname{I}}}
\newcommand{\norm}[1]{{\left\|#1\right\|}}
\newcommand{\scal}[2]{{\left\langle #1,#2\right\rangle}}

\allowdisplaybreaks

\makeatletter
\@namedef{subjclassname@2020}{%
  \textup{2020} Mathematics Subject Classification}
\makeatother

\begin{document}

\frenchspacing

\title[Pathwise uniqueness for stochastic heat and damped equations]{Pathwise uniqueness for stochastic heat and damped equations with H\"older continuous drift}

\author[D. Addona and D. A. Bignamini]{Davide Addona, Davide A. Bignamini}

\address{D. Addona: Dipartimento di Scienze Matematiche, Fisiche e Informatiche, Universit\`a degli studi di Parma, Parco Area delle Scienze 53/a (Campus), 43124 PARMA, Italy}
\email{\textcolor[rgb]{0.00,0.00,0.84}{davide.addona@unipr.it}}

\address{D. A. Bignamini: Dipartimento di Scienza e Alta Tecnologia (DISAT), Universit\`a degli Studi dell'In\-su\-bria, Via Valleggio 11, 22100 COMO, Italy}
\email{\textcolor[rgb]{0.00,0.00,0.84}{da.bignamini@uninsubria.it}}


\begin{abstract}
In this paper, we prove pathwise uniqueness for stochastic differential equations in infinite dimension. Under our assumptions, we are able to consider the stochastic heat equation up to dimension $3$, the stochastic damped wave equation in dimension $1$ and the stochastic Euler-Bernoulli damped beam equation up to dimension $3$. 
\end{abstract}


\keywords{pathwise uniqueness; regularization by noise; H\"older continuous drift; It\^o-Tanaka trick; stochastic damped wave equation; stochastic heat equation}


\subjclass[2020]{60H15,60H50}

\maketitle

\section{Introduction}
The aim of this paper is to prove pathwise uniqueness for mild solutions to a class of stochastic differential equations in a real separable Hilbert space $H$, given by
\begin{align}
\label{intro_SDE}
\left\{
\begin{array}{ll}
dX(t)=AX(t)dt+B(X(t))dt+GdW(t), & t\in[0,T], \\ [1mm]
X(0)=x\in H.   
\end{array}
\right.
\end{align}
Here, $A:D(A)\subseteq H\to H$ is the infinitesimal generator of a strongly continuous analytic semigroup $\{e^{tA}\}_{t\geq0}$ on $H$, $B:H\to H$ is a bounded and $\theta$-H\"older continuous function for some $\theta\in(0,1)$, $G:U\rightarrow H$ is a linear bounded operator and $W:=\{W(t)\}_{t\geq 0}$ is a $U$-cylindrical Wiener process, where $U$ is another real separable Hilbert space.

Pathwise uniqueness plays a crucial role in the investigation of existence of strong solutions to equation \eqref{intro_SDE} since the pioneering result due to Yamada and Watanabe \cite{Ya-Wa71}, who proved that if a stochastic differential equation in finite dimension admits existence of weak solutions and pathwise uniqueness, then existence of strong solutions follows at once. In \cite{Zv74}, the author introduces the so-called Zvonkin transformation, which allows to remove a drift term by means of a suitable change of coordinates using the It\^o formula, and then applies the result of \cite{Ya-Wa71} to construct strong solutions to a class of stochastic differential equations with rough drift coefficients. Generalizations of the results of \cite{Zv74} can be found in \cite{Ver} and in \cite{Kr-Ro05}, where the authors prove strong uniqueness under weaker assumptions on the drift term.

The first extension of Yamada-Watanabe result to infinite dimension appears in \cite{On04}, where it is shown that an analogous result holds for stochastic differential equations with values in $2$-smooth Banach spaces. 

In the following years the problem of pathwise uniqueness and of regularization by noise for stochastic evolution equations as \eqref{intro_SDE} has been widely studied, see for instance \cite{Cer-Dap-Fla2013, Dap-Fla2010,Dap-Fla2014,DPFPR13,DPFPR15,DPFRV16, GyPa93, MasPri17,MasPri23,Pri2021,Zam00}. One of the main tools to prove pathwise uniqueness in infinite dimension is the so-called It\^o-Tanaka trick, which consists in replacing the bad drift term with the solution to a suitable Kolmogorov equation. This trick has been introduced in \cite{FlGuPr10}, where the authors obtain well-posedness of the transport equation perturbed with a multiplicative noise.

In \cite{Dap-Fla2010} Da Prato and Flandoli prove pathwise uniqueness for a class of stochastic parabolic equations with a (bounded) H\"older continuous perturbation in the drift. In the quoted paper the authors perform the It\^o-Tanaka trick by means of a finite-dimensional approximation of \eqref{intro_SDE} and It\^o formula.

A different technique has been exploited in \cite{AddMasPri23, MasPri17,MasPri23}, where the pathwise uniqueness is gained for a class of semilinear stochastic damped beam equations and stochastic heat equations in \cite{AddMasPri23} and wave equations in \cite{MasPri17, MasPri23}. In \cite{AddMasPri23,MasPri17, MasPri23}, the It\^o-Tanaka trick is provided by means of systems of forward-backward stochastic differential equations. We stress that, in this approach, it is necessary to assume the so-called {\it structure condition}, i.e., $B=GC$ in \eqref{eqFO} for some function $C:H\to U$. Under similar conditions, in \cite{AddBig2} a Lipschitz dependence on the initial datum of the solution to \eqref{intro_SDE} is proved, i.e., for every $T>0$ there exists a positive constant $C_T$ such that, if $X$ and $Y$ are two weak mild solutions defined on the same probability space with respect to the same $U$-valued cylindrical Wiener process $W$ with initial datum $x$ and $y$, respectively, then
\begin{align*}
\sup_{t\in[0,T]}\mathbb E\|X(t)-Y(t)\|_H^2\leq C_T\|x-y\|_H^2.    
\end{align*}

In this paper, we prove pathwise uniqueness for families of both stochastic heat and damped equations with bounded H\"older continuous perturbation in the drift term which are not covered by the previous papers by means of a unified approach that does not require the {\it structure condition}. In particular, we apply our results to stochastic damped wave equation in dimension $1$ (see Corollary \ref{coro:damped_wave_dim_1}), to stochastic Euler-Bernoulli damped beam equation up to dimension $3$ (see Corollary \ref{coro:beam_eq}) and to stochastic heat equation up to dimension $3$ (see Theorem \ref{Heatequation}). To the best of our knowledge, this is the first time such a result has been achieved for the stochastic damped wave equation, as well as for the stochastic Euler-Bernoulli damped beam equation up to dimension $3$. Further, we are able to consider the stochastic heat equation up to dimension $3$, which is not reached in \cite{Dap-Fla2010} and in \cite{AddMasPri23} it is obtained assuming the structure condition. These results can be classified into {\it regularization by noise theory} in view of counterexamples presented in Subsection \ref{subsec:counter} and \cite[Subsection 6.2]{Dap-Fla2010}.
The approach which we introduce is partially inspired by \cite{Dap-Fla2010}, but there are some substantial differences, which we briefly list below and will be discussed in detail throughout the paper.
\begin{itemize}
\item In \cite{Dap-Fla2010} the authors only consider the case $U=H$. However, the possibility to allow $U\neq H$ is crucial in view of the applications, since in the abstract formulation of damped wave and beam equation it is necessary to take $U\neq H$, see Subsection \ref{Damped}.
\item We do not require that $A$ is self-adjoint and that there exists an orthonormal basis of $H$ consisting of eigenvectors of $A$.
\item We develop a finite-dimensional modified It\^o-Tanaka trick, in which we do not replace the H\"older nonlinearity by means of the solution to a Kolmogorov equation, but we provide an alternative formulation with some correction terms which eventually compensate the bad behavior of the drift part. This modification arises from the choice of the approximating sequence, which is different from that considered in \cite{Dap-Fla2010}. We refer to Remarks \ref{rmk:modified_I-T} and \ref{rmk:appr_seq} for a detailed discussion.
\item Under our assumptions the solution of the integral Kolmogorov equation associated to \eqref{intro_SDE} is not necessarily twice Frechét differentiable, see Remark \ref{2-derivate}.
\end{itemize}
In a sense, this paper can be seen as an extension of \cite{Dap-Fla2010} to the case where $A$ is a generic operator associated with a parabolic problem.

Let us spend few words on the assumption that the semigroup $\{e^{tA}\}_{t\geq0}$ is analytic. This requirement appears since we look for a positive constant $c$ such that
\begin{align}
\label{funz_per_dis_max}
\left\|A\int_0^\cdot e^{(\cdot-s)A}f(s)ds\right\|_{L^2([0,T];H)}\leq c\|f\|_{L^2([0,T];H)}   
\end{align}
holds true for every $\|f\|_{L^2([0,T];H)}$. An estimate of the form of \eqref{funz_per_dis_max} is known as $L^2$-{\it maximal regularity}, and it is verified if and only if $\{e^{tA}\}_{t\geq0}$ is analytic (see \cite{desimon64,weis01}). The $L^2$-{\it maximal regularity} is essential both in this paper and in \cite{Dap-Fla2010}.

The paper is organized as follows. In Section \ref{sec:notation} we fix the notation and we recall the main definitions which will be used in the paper.
In Section \ref{sec:ass_main_res} we state the assumptions and the main theorems. 
In Section \ref{sec:main_results} we develop our finite dimensional approximations and we prove the results stated in Section \ref{sec:ass_main_res}.
In Section \ref{sec:appl} we exhibit two classes of SPDEs to which our results apply: a family of stochastic damped equations which describe elastic systems, such as the stochastic damped wave equation and the stochastic Euler-Bernoulli beam equation, and a family of stochastic heat equations.
Appendix \ref{sec:Equazioni-Kolmogorov} is devoted to the study of Kolmogorov equations in Hilbert spaces and $L^2$-{\it maximal regularity}.

\section{Notation}
\label{sec:notation}

Let $\K$ be a Banach space endowed with the norm $\norm{\cdot}_\K$. We denote by $\mathcal{B}(\K)$ the Borel $\sigma$-algebra associated to the norm topology in $\K$. 

Let $(\Omega,\mathcal{F},\mathbb{P})$ be a probability space and let $\xi:(\Omega,\mathcal{F},\mathbb{P}) \rightarrow (\K,\mathcal{B}(\K))$ be a random variable. We denote by 
\[
\mathbb{E}[\xi]:=\int_\Omega \xi(w)\mathbb{P}(d\omega)=\int_\K x[\mathbb{P}\circ\xi^{-1}](dx)
\]
the expectation of $\xi$ with respect to $\mathbb{P}$. Let $\{Y(t)\}_{t\geq 0}$ be a $\K$-valued stochastic process defined on a normal filtered probability space $(\Omega,\mathcal{F},\{\mathcal{F}_t\}_{t\geq 0},\mathbb{P})$. We say that $\{Y(t)\}_{t\geq 0}$ is pathwise continuous $a.s.$ (almost surely) if there exists $\Omega_0\subseteq \Omega$ such that $\mathbb{P}(\Omega_0)=1$ and for every $\omega\in\Omega_0$ the function $t\rightarrow Y(t)(\omega)$ is continuous. 

Let $\X$ be a separable Hilbert space and let $\{g_k\}_{k\in\N}$ be an orthonormal basis of $\X$. We call $\X$-cylindrical Wiener process a stochastic process $\{W(t)\}_{t\geq 0}$ defined by 
\[
W(t):=\sum_{k=1}^{\infty} \beta_k(t)g_k \qquad \forall t\geq0,
\]
where  $\{\beta_1(t)\}_{t\geq 0}, \{\beta_2(t)\}_{t\geq 0},...,\{\beta_k(t)\}_{t\geq 0},...$ are real independent Brownian motions on a probability space $(\Omega,\mathcal{F},\mathbb{P})$.

Let $\K_1$ and $\K_2$ be two real separable Banach spaces equipped with norms $\norm{\cdot}_{\K_1}$ and $\norm{\cdot}_{\K_2}$, respectively. We denote by $B_b(\K_1;\K_2)$ the set of bounded and Borel measurable functions from $\K_1$ into $\K_2$. If $\K_2=\R$, then we simply write $B_b(\K_1)$. We denote by $C_b(\K_1;\K_2)$ ($UC_b(\K_1;\K_2)$, respectively) the space of bounded and continuous (uniformly continuous, respectively) functions from $\K_1$ into $\K_2$. We endowed $C_b(\K_1;\K_2)$ and $UC_b(\K_1;\K_2)$ with the norm
\[
\norm{f}_{\infty}=\sup_{x\in \K_1}\|f(x)\|_{\K_2}.
\]
If $\K_2=\R$ we simply write $C_b(\K_1)$ and $UC_b(\K_1)$, respectively.

Let $\theta\in (0,1)$. We denote by $C_b^\theta(\K_1;\K_2)$ the subspace of $C_b(\K_1;\K_2)$ of the $\theta$-H\"older continuous functions. The space $C_b^\theta(\K_1;\K_2)$ is a Banach space if it is endowed with the norm
\begin{align*}
\norm{f}_{C_b^\theta(\K_1;\K_2)}:=\norm{f}_\infty+[f]_{C_b^\theta(\mathcal{K}_1;\mathcal{K}_2)},
\end{align*}
where $[\cdot]_{C_b^\theta(\mathcal{K}_1;\mathcal{K}_2)}$ denote the standard seminorms on $C_b^\theta(\mathcal{K}_1;\mathcal{K}_2)$.
If $\K_2=\R$ then we simply write $C_b^\theta(\K_1)$.

We denote by $\Id_{\K_1}$ the identity operator on $\K_1$.
For $k\in\N$ we set $\mathcal{L}^{(k)}(\K_1;\K_2)$ the space of continuous multilinear mappings from  $\K_1^k:=\K_1\times\ldots\times\K_1$ into $\K_2$. If $k=1$ we simply write $\mathcal{L}(\K_1;\K_2)$, while if $\K_1=\K_2$ then we write $\mathcal{L}^{(k)}(\K_1)$.

Let $k\in\N$ and let $f:\K_1\rightarrow \K_2$ be a $k$-times Fr\'echet differentiable function. We denote by $D^i f(x)$, $i=1,\ldots,k$, its Fr\'echet derivative of order $i$ at $x\in\X$. In the case $\K_2=\mathbb{R}$ we denote by $\nabla f (x)$ and $\nabla^2f(x)$ the Fr\'echet gradient and Hessian at $x\in\X$, respectively. 
For $k\in\N$, we denote by $C_b^{k}(\K_1;\K_2)$ ($UC^k_b(\K_1;\K_2)$, respectively) the space of bounded, uniformly continuous and $k$ times Fr\'echet differentiable functions $f:\K_1\rightarrow\K_2$ such that $D^if\in C_b(\K_1;\mathcal{L}^{(i)}(\K_1;\K_2))$ ($D^if\in UC_b(\K_1;\mathcal{L}^{(i)}(\K_1;\K_2))$, respectively), for $i=1,\ldots,k$. We endow $C_b^{k}(\K_1;\K_2)$ and $UC_b^{k}(\K_1;\K_2)$  with the norm
\[
\norm{f}_{C_b^{k}(\K_1;\K_2)}:=\norm{f}_{\infty}+\sum^k_{i=1}\sup_{x\in\K_1}\|D^if(x)\|_{\mathcal{L}^{(i)}(\K_1;\K_2)}.
\]
We set $C^\infty_b(\K_1,\K_2)=\bigcap_{k\geq 1}C^k_b(\K_1;\K_2)$. If $\K_2=\R$ then we simply write $C_b^k(\K_1)$ and $UC_b^k(\K_1)$, respectively. Let $T>0$. We denote by $C_b^{0,1}([0,T]\times\K_1;\K_2)$ the space of functions $f:[0,T]\times\K_1\rightarrow\K_2$ such that $f(t,\cdot)\in C^1_b(\K_1;\K_2)$ for any $t\in [0,T]$ and $f(\cdot,x)\in C([0,T];\K_2)$ for any $x\in \K_1$. If $\K_2=\R$ then we simply write $C_b^{0,1}([0,T]\times\K_1)$.

We still denote by $\K_1$ the complexification of $\K_1$. Let $T:{\rm Dom}(T)\subseteq\K_1\rightarrow\K_1$ be a linear operator. We define the resolvent set of $T$ as 
\[
\rho(T):=\{\lambda\in\mathbb{C}\; : \; (T-\lambda\Id):{\rm Dom}(T)\rightarrow\K_1 \mbox{ is bijective and its inverse is bounded}  \}.
\]
The set $\sigma(T):=\mathbb{C}\backslash\rho(T)$ is called spectrum of $A$. Moreover, for every $\lambda\in \rho(T)$ we define the resolvent $R(\lambda,T)$ of $T$ as
\[
R(\lambda,T)=(T-\lambda\Id)^{-1}.
\]

Let $\X$ be a separable Hilbert space equipped with the inner product $\scal{\cdot}{\cdot}_{\X}$. We say that $Q\in\mathcal{L}(\X)$ is \emph{non-negative} (\emph{positive}) if for every $x\in \X\backslash\{0\}$
\[
\langle Qx,x\rangle_\X\geq 0\ (>0).
\]
On the other hand, $Q \in \mathcal{L}(\X)$ is a \emph{non-positive} (respectively, \emph{negative}) operator if $-Q$ is non-negative (respectively, positive). Let $Q\in\mathcal{L}(\X)$ be a non-negative and self-adjoint operator. We say that $Q$ is a trace-class operator if
\begin{align}\label{trace_defn}
{\rm Trace}_{\mathcal X}[Q]:=\sum_{n=1}^{\infty}\langle Qe_n,e_n\rangle_\X<\infty
\end{align}
for some (and hence for all) orthonormal basis $\{e_n:n\in\N\}$ of $\X$. We recall that the definition of trace operator given by \eqref{trace_defn} is independent of the choice of the orthonormal basis. 
Let $Y$ be another separable Hilbert space and let $R\in\mathcal{L}(\X;Y)$. We say that $R$ is a Hilbert--Schmidt operator if 
\[
\norm{R}^2_{\mathcal L_2(\K_1;K_2)}:=\sum^{\infty}_{k=1}\norm{Rg_k}^2_Y<\infty
\]
for some (and hence for all) orthonormal basis $\{g_k:k\in\N\}$ of $\X$. It follows that if $R$ is a Hilbert--Schmidt operator, then $RR^*$ and $R^*R$ are trace-class operators and
\begin{align*}
{\rm Trace}_{Y}[RR^*]={\rm Trace}_{\mathcal X}[R^*R]=\|R\|^2_{\mathcal L_2(\mathcal K_1;\mathcal K_2)}.    
\end{align*}

\section{Assumptions and main results}
\label{sec:ass_main_res}

Let $H$ and $U$ be real separable Hilbert spaces, let $B:H\to H$ be a bounded and $\theta$-H\"older continuous function for some $\theta\in (0,1)$, let $G\in\mathcal{L}(U;H)$, let $A:D(A)\subseteq H\to H$ be the infinitesimal generator of a strongly continuous and analytic semigroup $\{e^{tA}\}_{t\geq 0}$ on $H$ and let $W$ be a $U$-cylindrical Wiener process. For every fixed $T>0$, we aim to study well-posedness for equations of the form
\begin{align}\label{eqFO}
\left\{
\begin{array}{ll}
\displaystyle  dX(t)=AX(t)dt+B(X(t))dt+GdW(t), \quad t\in[0,T],  \vspace{1mm} \\
X(0)=x\in H.
\end{array}
\right.
\end{align}
We begin to define the notions of weak and strong solutions to \eqref{eqFO}.

\begin{defn}\label{weak-solution}
 Let $T>0$ and $x\in H$.
\begin{itemize}
\item[(Weak)] A weak (mild) solution to \eqref{eqFO} is a couple $(X,W)$ where $W=\{W(t)\}_{t\in [0,T]}$ is a $U$-cylindrical Wiener process defined on a filtered probability space $(\Omega,\mathcal{F},\{\mathcal{F}_t\}_{t\in [0,T]},\mathbb{P})$ and $X=\{X(t,x)\}_{t\in[0,T]}$ is a $H$-valued $\{\mathcal{F}_t\}_{t\in [0,T]}$-adapted stochastic process such that for every $t\in [0,T]$ 
\begin{align}\label{mild}
X(t,x)=e^{tA}x+\int_0^te^{(t-s)A}B(X(s,x))ds+W_A(t),\qquad \mathbb{P}{\rm -a.s.},    
\end{align}
where $\{W_A(t)\}_{t\geq 0}$ is the stochastic convolution process given by 
\begin{equation*}
W_A(t):=\int^t_0e^{(t-s)A}GdW(s),\qquad \mathbb P{\rm-a.s.}
\end{equation*}
for every $t\geq0$.
\item[(Strong)] We say that strong existence holds true for \eqref{eqFO} if for every $U$-cylindrical Wiener process $W=\{W(t)\}_{t\in [0,T]}$, defined on a complete filtered probability space $(\Omega,\mathcal{F},\{\mathcal{F}_t\}_{t\in [0,T]},\mathbb{P})$, there exists a $H$-valued process $X=\{X(t,x)\}_{t\in[0,T]}$, adapted with respect to the filtration $\{\mathcal{F}_t\}_{t\in [0,T]}$, such that for every $t\in [0,T]$ it satisfies \eqref{mild}.
\end{itemize}
\end{defn}

\begin{rmk}
We refer to \cite[Theorem 3.6 and Proposition 6.8]{Kun2013} for the equivalence between martingale, analytic weak, analytic weak mild and mild solutions.   
\end{rmk}

We introduce the concepts of weak and pathwise uniqueness for \eqref{eqFO}.

\begin{defn}\label{uniqueness}
Let $T>0$ and $x\in H$. 

\begin{itemize}

\item[(Weak)] We say that weak uniqueness holds for \eqref{eqFO} if whenever $(X_1,W)$ and $(X_2,W)$ are two weak mild solutions to \eqref{eqFO}, then $X_1$ and $X_2$ have the same law on $C([0,T];H)$, namely that for every continuous and bounded $\psi:C([0,T];H)\to \R$ we have
\[
\mathbb{E}\left[\psi(X_1)\right] = \mathbb{E}\left[\psi(X_2)\right].
\]

\item[(Pathwise)] We say that pathwise (or strong) uniqueness holds true for \eqref{eqFO} if whenever $(X_1,W)$ and $(X_2,W)$ are two weak mild solutions to \eqref{eqFO} defined on the same probability space $(\Omega,\mathcal F,\{\mathcal F_t\}_{t\geq0},\mathbb P)$ with same $U$-cylindrical process $W$, then 
\[
\mathbb{P}\left(\left\lbrace \omega\in \Omega\; :\; X_1(t,x)(\omega)=X_2(t,x)(\omega),\; \forall t\in [0,T]   \right\rbrace\right)=1.
\]
\end{itemize}
\end{defn}

We now state the assumptions exploited to perform the finite-dimensional approximations and the modified It\^o--Tanaka trick which will be crucial to prove the main result of this paper.
\begin{hyp1}\label{hyp:finito-dimensionale}
The following conditions hold true.
\begin{enumerate}[\rm(i)]
\item $A:{\rm Dom}(A)\subseteq H\rightarrow H$
is the infinitesimal generator of a strongly continuous analytic semigroup $\{e^{tA}\}_{t\geq0}$.

\item $G\in\mathcal L(U;H)$ satisfies
\[
G=\widetilde{G}\mathcal{V}
\]
where $\widetilde{G}\in\mathcal{L}(U;H)$ and $\mathcal{V}\in \mathcal{L}(U;U)$.

\item $B\in C_b^\theta(H;H)$ with $\theta\in (0,1)$ such that
\[
B=\widetilde{G}\widetilde{B},
\]
where $\widetilde{B}\in C_b^\theta(H;U)$.
\item There exists $\eta\in (0,1)$ such that for every $t>0$ we have
\[
\int^t_0\frac{1}{s^\eta}{\rm Trace}_H\left[e^{sA}GG^*e^{sA^*}\right]ds<\infty.
\]

\item \label{Accan}There exists a sequence of finite-dimensional subspaces $\{H_n\}_{n\in\N}\subseteq H$ such that $H=\overline{\cup_{n\in\N}H_n}$, $H_0:=\{0\}$ and for every $n\in\N$ we have
\begin{align*}
&H_{n-1}\subseteq H_{n},\qquad H_{n-1}\subseteq {\rm Dom}(A),\qquad A(H_{n}\cap H_{n-1}^\perp)\subseteq \left(H_{n}\cap H_{n-1}^\perp\right).
\end{align*}

\item\label{contrin} For every $t>0$ we have
\begin{align}
   &e^{tA}(H)\subseteq Q^{\frac12}_t(H),\qquad\qquad \qquad\qquad\quad Q_{t}:=\int^t_0e^{sA}GG^*e^{sA^*}ds;\label{contron}\\
&\int^t_0\norm{\Gamma_s}^{1-\theta}_{\mathcal{L}(H)}\|\Gamma_s\widetilde{G}\|_{\mathcal{L}(U;H)}ds<\infty, \quad\qquad\Gamma_t:=Q^{-\frac12}_te^{tA}\label{supercontron},
\end{align}
where $s\mapsto \norm{\Gamma_s}^{1-\theta}_{\mathcal{L}(H)}$ and $s\mapsto\|\Gamma_s\widetilde{G}\|_{\mathcal{L}(U;H)}$ are bounded from below functions in $(0,t)$ for every $t>0$. Further, we assume that there exists $\theta'<\theta$ such that
\begin{align*}
\int^t_0\norm{\Gamma_s}^{1-\theta'}_{\mathcal{L}(H)}ds<\infty.  
\end{align*}

\end{enumerate}
\end{hyp1}

\begin{remark}\label{remarkhyp1}
$ $
\begin{enumerate}[\rm(i)]
\item  The operator $\mathcal{V}$ in Hypotheses \ref{hyp:finito-dimensionale}$(i)$ represents the color of the noise driving the SPDE \eqref{eqFO}. On the other hand, the operator $\widetilde{G}$ is an auxiliary operator that allows us to cover a wide class of SPDEs in the abstract form \eqref{eqFO}, such as the class of damped equations discussed in Subsection \ref{Damped}.

\item When we consider perturbed versions of the Heat equation we can set $U=H$ and $\widetilde{G}=\Id_H$, see  Subsection \ref{Heat-case}.

\item The aim of Hypotheses \ref{hyp:finito-dimensionale}\eqref{Accan} is to generalize the case when there exists an orthonormal basis $\{e_k:k\in\N\}$ of $H$ consisting of eigenvectors of $A$. This is the case when $A$ is a realization of the Laplace operator in $H=L^2(\mathcal{O})$ for some smooth and bounded subset $\mathcal{O}$ of $\mathbb{R}^n$. In this case, Hypotheses \ref{hyp:finito-dimensionale}\eqref{Accan} is verified with $H_n:={\rm span}\{e_1,...,e_n\}$ for every $n\in\N$. However, there are some significant situations where Hypotheses \ref{hyp:finito-dimensionale}\eqref{Accan} holds true but there is no orthonormal basis of $H$ consisting of eigenvectors of $A$; for instance, this happens when $A$ is a realization of the differential operator driven a damped wave equation or a damped Euler-Bernoulli beam equation in $H=L^2(\mathcal O)\times L^2(\mathcal O)$ $($see Subsection \ref{Damped}$)$.
\item We note that there exists an orthonormal basis $\{g_k:k\in\N\}$ of $H$ such that 
\[
H_n={\rm span}\{g_1,...,g_{s_n}\}\quad \forall n\in\N,
\]
where $s_n:={\rm Dim}(H_n)$. However, in general the basis $\{g_k:k\in\N\}$ does not consist of eigenvectors of $A$ $($see point {\rm(iii)}$)$.
\item From Hypotheses \ref{hyp:finito-dimensionale}(vi), it follows that
\[
\int^T_0\|\Gamma_t\|^{1-\theta}_{\mathcal{L}(H)}dt<\infty.
\]
In the quoted assumption, the condition on $\theta'$ is just technical and automatically verified in our examples, where $\|\Gamma_t\|_{\mathcal L(H)}$ behaves near $0$ like $t^{-\sigma}$ for some positive $\sigma<1$.
\end{enumerate}
\end{remark}

If $\{e^{tA}\}_{t\geq0}$ is compact, then under the previous assumptions, the SPDE \eqref{eqFO} is well-posed in the weak sense.
\begin{proposition}[Theorem 2.6 of \cite{AddBig3} and Proposition 3 of \cite{Cho-Gol1995}]
Assume that Hypotheses \ref{hyp:finito-dimensionale} hold true and that $\{e^{tA}\}_{t\geq 0}$ is a compact semigroup. For every $x\in H$ and $T> 0$,
the SPDE \eqref{eqFO} admits a weak mild solution and weak uniqueness holds true.
\label{prop:weak_tutto}
\end{proposition}

\begin{rmk}
We refer to \cite{Kun2013,Kun2013-2} for other results about the weak well-posedness. For a discussion on weak uniqueness we refer to \cite{BerOrrSca2024, Pri2021} for the stochastic heat equation with singular drift and to \cite{Han2024} for the stochastic wave equation with multiplicative noise.
\end{rmk}

In view of Proposition \ref{prop:weak_tutto}, in the present paper we are interested in an abstract result which guarantees that \eqref{eqFO} is well-posed in strong sense.

Now we present the modified It\^o-Tanaka trick which we will perform to reach the desired result. Let $(X,W)$ be a weak solution to \eqref{eqFO}. 
For every $n\in\N$ and $x\in H$ we consider the $H_n$-valued stochastic process $\{X_{n}(t,x)\}_{t\in [0,T]}$ which for every $t\in[0,T]$ fulfills
\begin{equation}\label{app-intro}
X_n(t,x):=e^{tA_n} P_nx+\int^t_0e^{(t-s)A_n}B_n(X(s,x))ds+\int^t_0 e^{(t-s)A_n}G_ndW(s), \qquad \mathbb{P}{\rm -a.s.},
\end{equation}
where, for every $n\in\N$,

\begin{align*}
 B_n(\cdot)=P_nB(P_n(\cdot)),\quad A_n=AP_n=P_nA,\quad G_n=P_nG,\quad n\in\N
\end{align*}
and $P_n$ is the orthogonal projection on  $H_n$.
We notice that $\{X_{n}(t,x)\}_{t\in [0,T]}$ solves in the classical It\^o sense the SDE
\begin{align}\label{SDEn-intro}
\left\{
\begin{array}{ll}
\displaystyle  dX_n(t)=A_nX_n(t)dt+B_n( X(t,x))dt+G_ndW(t), \quad t\in[0,T],  \vspace{1mm} \\
X_n(0)=P_nx\in H_n,
\end{array}
\right.
\end{align}
namely, for every $t\in [0,T]$, it holds
\[
X_n(t,x)=P_nx+\int_0^t(A_nX_n(s,x)+B_n(X(s,x)))ds +G_nW(t), \qquad \mathbb{P}{\rm -a.s.}
\]
Problem \eqref{SDEn-intro} is a linear non homogeneous equation in $X_n$ of the form $dX_n(t)=A_nX_ndt+f(t)dt+G_ndW(t)$, where $f$ is a stochastic perturbation. We emphasize that, in general, the process $\{X_n(t,x)\}_{t \in [0,T]}$ is not a strong solution to equation \eqref{SDEn-intro}. This is because the weak mild solution $(X, W)$ is fixed a-priori; consequently, if we consider another cylindrical Wiener process and another filtration, it is not guaranteed that the process $\{B_n(X(t,x))\}_{t \in [0,T]}$ is adapted to the new filtration. However, in the method presented in this paper we stress that it is not necessary that $\{X_n(t,x)\}_{t \in [0,T]}$ is a strong solution to \eqref{eqFO}.

We consider the following backward Kolmogorov equation
\begin{equation}\label{Back-Kolmo-intro}
    U_n(t,x)=\int^T_t \mathcal{R}_n(r-t)\left( DU_n(r,\cdot)B_n(\cdot)+B_n(\cdot)\right)(x)dr ,\quad t\in [0,T], \ x\in H_n,
\end{equation}
where $\{\mathcal{R}_n(t)\}_{t\geq 0}$ is the vector-valued Ornstein-Uhlenbeck semigroup associated to \eqref{SDEn-intro} with $B_n\equiv 0$. We will show that \eqref{Back-Kolmo-intro} has a unique solution $U_n\in C^{0,1}_b([0,T]\times H_n;H_n)$ (see Appendix \ref{sec:Equazioni-Kolmogorov}).
Under our assumptions, it is not guaranteed that $U_n$ admits neither first order time derivative nor spatial derivatives of order $2$. However, by applying the It\^o formula to a smooth approximation of the processes $\{U_{n}(t,X_n(t,x))\}_{t\in[0,T]}$ we get the following representation formula for $X_n$. 
\begin{thm}\label{Ito-Tanaka}
Assume that Hypotheses \ref{hyp:finito-dimensionale} hold true. For every $n\in\N$ and $x\in H$ the solution $\{X_n(t,x)\}_{t\in [0,T]}$ to \eqref{app-intro} satisfies the following equality: for every $t\in[0,T]$,
\begin{align}
X_n(t,x)&=e^{tA_n}( P_nx+U_n(0,P_nx))-U_n(t,X_n(t,x))-A_n\int^t_0e^{(t-s)A_n}U_{n}(s,X_n(s,x))ds\notag\\
&+\int^t_0e^{(t-s)A_n} \left[B_n(X(s,x))-B_n(X_n(s,x))\right]ds\notag\\
&+\int^t_0e^{(t-s)A_n}DU_{n}(s,X_n(s,x))\left[B_n(X(s,x))-B_n(X_n(s,x)))\right]ds\notag\\
&+\int^t_0e^{(t-s)A_n}DU_n(s,X_n(s,x))G_ndW(s)+ \int^t_0e^{(t-s)A_n}G_ndW(s),\qquad \mathbb{P}-{\rm a.s.}\label{I-T}
\end{align}
where $U_n\in C^{0,1}_b([0,T]\times H;H)$ is the unique solution to \eqref{Back-Kolmo-intro}.
\end{thm}


\begin{remark}
\label{rmk:modified_I-T}
At this stage, one expects to let $n$ go to infinity obtaining that, for every $t\in[0,T]$,
\begin{align}
X(t,x)
= & e^{tA}(x+U(0,x))-U(t,X(t,x))-A\int^t_0e^{(t-s)A}U(s,X(s,x))ds\notag\\
&+\int^t_0e^{(t-s)A}DU(s,X(s,x))GdW(s)+ \int^t_0e^{(t-s)A}GdW(s),\qquad  \mathbb{P}-{\rm a.s.},\label{I-T_X}    
\end{align}
where $U$ is the unique solution to the backward integral equation 
\begin{equation*}
U(t,x)=\int^T_t \mathcal{R}(r-t)\left[ DU(r,\cdot)B(\cdot)+B(\cdot)\right](x)dr ,\quad t\in [0,T], \ x\in H
\end{equation*}
and $\{\mathcal{R}(t)\}_{t\geq 0}$ is the vector-valued Ornstein-Uhlenbeck semigroup associated to \eqref{eqFO} with $B\equiv 0$, and then exploit direct computations on \eqref{I-T_X} $($see Appendix \ref{sec:Equazioni-Kolmogorov}$)$. Unfortunately, even if in Subsection \ref{sezione-finito-dimensionali} we prove that $\{X_{n}(t,x)\}_{t\in [0,T]}$ converges to $\{X(t,x)\}_{t\in [0,T]}$ in $L^2([0,T]\times \Omega, \mathcal{B}([0,T])\times \mathcal{F},\lambda\times\mathbb{P})$, where $\lambda$ is the Lebesgue measure on $[0,T]$, it is not clear if $U_n$ converges to $U$ with respect to some suitable norm which allows us to deduce the convergence of \eqref{I-T} to \eqref{I-T_X}. Hence, in Section \ref{sec:main_results} we perform our computations on \eqref{I-T} $($which is more involved with respect to \eqref{I-T_X} since it includes additional addends$)$ taking advantage of estimate  \eqref{supercontron}, which is independent of $n$. 

Further, the difference
$B_n(X(s,x))-B_n(X_n(s,x))$ appears in the second and in the third integral of \eqref{I-T} since in \eqref{app-intro} the argument of $B_n$ is $X$ and not $X_n$. This is not a big deal, since we will prove that
\begin{align*}
& \int^T_0 \mathbb E\left\|\int_0^te^{(t-s)A_n} \left(B_n(X(s,x))-B_n(X_n(s,x))\right)ds\right\|_H^2dt\to 0, \\
&\int_0^T\mathbb E\left\|\int^t_0e^{(t-s)A_n}DU_{n}(s,X_n(s,x))\left(B_n(X(s,x))-B_n(X_n(s,x)))\right)ds\right\|^2_Hdt \to 0    
\end{align*}
as $n$ tends to $\infty$.
\end{remark}


\begin{remark}\label{rmk:appr_seq}
We stress that, given a weak solution $(X,W)$, the processes $\{X_{n}(t,x)\}_{t\in [0,T]}$ used in this paper $($see \eqref{app-intro}$)$ explicitly depend on $\{X(t,x)\}_{t\geq 0}$.
Instead of $\{X_{n}(t,x)\}_{t\in [0,T]}$, one should be tempted to consider the sequence given by the solutions $\widehat{X}_n:=\{ \widehat X_n(t,x)\}_{t\geq 0}$ to the following finite-dimensional version of \eqref{eqFO} 
\begin{align*}
\left\{
\begin{array}{ll}
\displaystyle  d\widehat{X}_n(t)=A_n\widehat{X}_n(t)dt+B_n(\widehat{X}_n(t))dt+{G_n}dW(t), \quad t\in[0,T],  \vspace{1mm} \\
\widehat{X}_n(0)=P_nx\in H_n,
\end{array}
\right.
\end{align*}
$($see, for instance, \cite[Lemma 6]{Dap-Fla2010}$)$. In this case, in \eqref{I-T} both $X$ and $X_n$ are replaced by $\widehat X_n$, the second and the third integral in \eqref{I-T} vanish and so \eqref{I-T} simplify. However, processes $\{\widehat X_n(t,x)\}_{t\geq 0}$ do not depend on $\{X(t,x)\}_{t\geq 0}$ and nothing ensures that $\widehat X_n$ converges to $X$ as $n$ goes to $\infty$ for every $x\in H$.

We stress that the convergence of $(\widehat X_n)_{n\in\N}$ to $X$ as $n$ goes to $\infty$ for every $x\in H$ is true a-posteriori, as a consequence of the pathwise uniqueness and our computations. Indeed, in Proposition \ref{prop:conv_appr_buone} we show that, if the assumptions of the main theorem of this paper are verified $($Theorem \ref{pathwiseuniqueness}$)$, then $(\widehat X_n)_{n\in\N}$ converges to $X$ in $L^2([0,T]\times \Omega, \mathcal{B}([0,T])\times \mathcal{F},\lambda\times\mathbb{P})$.
\end{remark}

To prove pathwise uniqueness, we consider two weak mild solutions $(X_1,W)$ and $(X_2,W)$ to \eqref{eqFO}  which are defined on the same probability space $(\Omega,\mathcal{F},\{\mathcal{F}_t\}_{t\in [0,T]},\mathbb{P})$ and same $U$-cylindrical Wiener process $W$, and the approximating sequences $\{X_{1,n}(t,x)\}_{t\in [0,T]}$ and $\{X_{2,n}(t,x)\}_{t\in [0,T]}$ of $\{X_{1}(t,x)\}_{t\in [0,T]}$ and $\{X_{2}(t,x)\}_{t\in [0,T]}$, respectively, given by \eqref{app-intro}. \\ 
Exploiting Theorem \ref{Ito-Tanaka}, we estimate the $L^2$-norm of the difference between $\{X_{1,n}(t,x)\}_{t\in [0,T]}$ and $\{X_{2,n}(t,x)\}_{t\in [0,T]}$ and we show that this difference vanishes as $n$ goes to $\infty$. This yields the desired result. To this end, it is necessary that either Hypotheses \ref{hyp:traccia-finita} or Hypotheses \ref{hyp:goal-addo} (which generalizes \cite[Assumption 3]{Dap-Fla2010}) stated below are satisfied

\begin{hyp1}\label{hyp:traccia-finita}
The operator $\mathcal{V}\in\mathcal L_2(U)$ is a Hilbert--Schmidt operator.
\end{hyp1}

If necessary, in the following hypothesis we consider (without changing the notation) the complexification of $H$ and we take our assumptions on the complexified space.  
\begin{hyp1}\label{hyp:goal-addo}

 There exists a family of normalized $($but not necessarily orthogonal$)$ vectors $\{f_n:n\in\N\}$ of $H$ consisting of eigenvectors of $A^*$ such that $H=\overline{{\rm span}\{f_n:n\in\N\}}$ and there exists a sequence $(d_n)_{n\in\N}\subseteq \N$ such that 
\begin{enumerate}[\rm(a)]
\item for every $n\in\N$ we have
\begin{align*}    
\{f_1,\ldots,f_{s_n}\}=\bigcup_{i=1}^n \{e_1^i,\ldots,e_{d_i}^i\},\quad s_n=d_1+\ldots + d_n,
\end{align*}
where for every $i,j\in\N$ with $i\neq j$, we have
\[
\langle e^i_k,e^j_h\rangle_H=0,\qquad k=1,\ldots,d_i, \ h=1,\ldots,d_j.
\]
\item There exists $d\in\N$ such that $d_n\leq d$ for every $n\in\N$.

\item For every $n\in\N$ and $j=1,\ldots,d_n$, the eigenvalue $\rho^n_j$ associated to the eigenvector $e^n_j$ has negative real part. Moreover,
\begin{align}
\label{conv_serie_holder}
-\sum_{n\in\N}\sum_{j=1}^{d_n}\frac{\|B^n_j\|_{C_b^\theta(H)}^2}{{\rm Re}(\rho^n_j)}<\infty,   
\end{align}
where $B^n_j(\cdot)=\langle B(\cdot),e^n_j\rangle_H$ for every $n\in\N$ and $j\in\{1,\ldots,d_n\}$.
\end{enumerate}
\end{hyp1}
\begin{remark}$ $
    
\begin{enumerate}[\rm(i)]

\item Hypotheses \ref{hyp:goal-addo} imply Hypotheses \ref{hyp:finito-dimensionale}\eqref{Accan} with
\begin{align*}
& H_n:=\bigcup^n_{k=1}{\rm span}\left\{e^k_1,\ldots,e^k_{d_k}\right\},\quad H_{n}\cap H_{n-1}^\perp={\rm span}\left\{e^{n}_1,\ldots e^{n}_{d_{n}}\right\},\quad  n\in\N,\\    
& s_n:={\rm Dim}(H_n),\quad d_n={\rm Dim}(H_{n}\cap H_{n-1}^\perp),\quad n\in\N.       
\end{align*}
Indeed for every $n\in\N$, $i\in \{1,\ldots,d_n\}$ and $x=\sum_{k\in\N}\sum_{j=1}^{d_k}x_j^ke_j^k\in D(A^*)$ we have
\begin{align*}
\scal{A^*x}{e^n_i}&=\scal{\sum_{k\in\N}\sum_{j=1}^{d_k}x_j^k A^{*}e^k_j}{e^n_i}=\scal{\sum_{k\in\N}\sum_{j=1}^{d_k}\rho^k_jx_j^ke^k_j}{e^n_i}=\scal{\sum_{j=1}^{d_n}\rho^n_jx_j^ke^n_j}{e^n_i},
\end{align*}
which gives $|\langle A^*x,e_i^n\rangle| \leq c_n\|x\|$ for every $x\in D(A^*)$ with $c_n=\|A_{|H_n}\|_{\mathcal L(H_n)}$. Since $D(A^*)$ is dense in $H$ and $A^{**}=A$ $(A$ is closed$)$, by definition it follows that $e^n_i\in{\rm Dom}(A)$ for every $n\in\N$ and $i\in \{1,..,d_n\}$. Moreover, fixed $n\in \N$, for every $i\in \{1,\ldots,d_n\}$ we have
\[
\scal{Ae^n_i}{e_j^k}=\scal{e^n_i}{A^*e_j^k}=0,\quad \forall k\neq n,\; j\in \{1,\ldots,d_k\}.
\]
Hence, $A\left({\rm span}\{e^n_1,\ldots,e^n_{d_n}\}\right)\subseteq {\rm span}\{e^n_1,\ldots,e^n_{d_n}\}$ and so Hypotheses \ref{hyp:finito-dimensionale}\eqref{Accan} holds true.

\item In our examples, Hypotheses \ref{hyp:goal-addo}$(a)$-$(b)$ are always verified. Indeed, if $A$ is a realization of the Laplace operator in $H=L^2(\mathcal{O})$ for some smooth and bounded subset of $\mathbb{R}^n$  $\mathcal{O}$, then Hypotheses \ref{hyp:goal-addo} is verified with $d_n=1$ for every $n\in\N$, while if $A$ is a suitable realization of the differential operator driven a damped wave equation or a Euler-Bernoulli damped beam equation in $H=L^2(\mathcal O)\times L^2(\mathcal O)$, then Hypotheses \ref{hyp:goal-addo} is verified with $d_n=2$ for every $n\in\N$ $($see Subsection \ref{Damped}$)$.

\item We underline that, in many significant cases, the so-called structure condition (namely $\widetilde{B}=\mathcal{V}F$ with $F\in C_b^\theta(H;U)$, see Hypotheses \ref{hyp:finito-dimensionale}) is strictly stronger than \eqref{conv_serie_holder}, see Proposition \ref{structure-heat}.
\end{enumerate}
\end{remark}

Finally, we can state the main result of this paper.

\begin{thm}\label{pathwiseuniqueness}
Assume that Hypotheses \ref{hyp:finito-dimensionale} and either Hypotheses \ref{hyp:traccia-finita} or \ref{hyp:goal-addo} hold true. Then, for every $T>0$ and $x\in H$ pathwise uniqueness holds true for equation \eqref{eqFO}.
\end{thm}
\begin{remark}
We point out that, using a localization argument as in  \cite{DPFPR15}, the boundedness of $B$ can be skipped. Further, we may also allow $B$ to depend on $t$, by assuming that the conditions on $B$ hold true uniformly with respect to $t\in[0,T]$. 
\end{remark}

In the next corollary we establish strong existence by exploiting the results in \cite{On04}. 

\begin{corollary}\label{Strong}
Assume that Hypotheses \ref{hyp:finito-dimensionale}, that either Hypotheses \ref{hyp:traccia-finita} or \ref{hyp:goal-addo} hold true and that $\{e^{tA}\}_{t\geq 0}$ is a compact semigroup. Then strong existence holds true for equation \eqref{eqFO}.    
\end{corollary}

\begin{rmk}
We underline that in all the examples presented in Section \ref{sec:appl}, $\{e^{tA}\}_{t\geq 0}$ is a compact semigroup.
\end{rmk}

\section{Proofs of the main results}\label{sec:main_results}

In this section we are going to prove the results stated in the previous section.
\subsection{Finite-dimensional approximation}\label{sezione-finito-dimensionali}
In this subsection, we will present the finite-dimensional procedure which we will perform. Fix $T>0$ and let $(X,W)$ be a weak solution to \eqref{eqFO}.

Assume that Hypotheses \ref{hyp:finito-dimensionale} hold true, let $n\in\N$ and let $P_n$ be the orthogonal projection on  $H_n$. For every $n\in\N$ we set
\begin{align}\label{coefficienti-n}
A_n:=AP_n,\quad \widetilde{G}_n:=P_n \widetilde{G},\quad G_n=\widetilde{G}_n\mathcal{V}, \quad \widetilde{B}_n=\widetilde{B}\circ P_n,\quad B_n=\widetilde{G}_n\widetilde{B}_n.
\end{align}
For every $n\in \N$, the operator $A_{|H_n}$ belongs to $\mathcal L(H_n)$ and $A_n$ is the infinitesimal generator of a uniformly continuous semigroup $\{e^{tA_n}\}_{t\geq 0}$ both in $H$ and $H_n$. In the following lemma we have collected some useful properties of $\{e^{tA_n}\}_{t\geq 0}$ that we will exploit in the rest of the paper.

\begin{lemma}\label{Lemma:semi-n}
Assume that Hypotheses \ref{hyp:finito-dimensionale} hold true. For every $n\in\N$ we have
\begin{align}
& AP_nx=P_nAx,\qquad \forall x\in {\rm Dom}(A),\label{P1}\\
& P_nR(\lambda,A)x=R(\lambda,A)P_nx,\qquad \forall x\in H,\ \forall \lambda\in\rho(A),\label{P2}\\
& \rho(A)\subseteq\rho(A_n)\  \textrm{ and } \  R(\lambda,A)_{|H_n}=R(\lambda,A_n), \quad \forall \lambda\in\rho(A),\label{P3}\\
& P_ne^{tA}x=e^{tA}P_nx,\qquad \forall x\in H, \ \forall t\geq0,\label{P4}\\
& e^{tA}_{|H_n}=e^{tA_n}, \qquad \forall t\geq0.\label{P5}
\end{align}
Moreover, if we introduce the operator $S_n:L^2(0,T;H)\to L^2(0,T;H)$, defined as
\begin{align*}
(S_nf)(t):=\int_0^te^{(t-s)A_n}f(s)ds \qquad \forall t\in [0,T],   
\end{align*}
for every $f\in L^2(0,T;H)$, then $S_nf\in L^2(0,T;D(A))$ and there exists a positive constant $C_{T,A}$, independent of $n\in\N$, such that
\begin{align}\label{stima_conv_fourier_part_n}
\left(\int_0^T\norm{A_nS_nf(t)}^2_{H_n}dt\right)^{1/2}\leq  \|S_nf\|_{L^2(0,T;D(A))}\leq C_{T,A}\|f\|_{L^2(0,T;H)}.
\end{align}
\end{lemma}

\begin{proof}
Let $n\in\N$. By Hypotheses \ref{hyp:finito-dimensionale}\eqref{Accan}, if $x\in H_n$ then
\begin{align}
\label{comm_A_P_n_H_n}
P_n Ax=P_nAP_nx=AP_nx.
\end{align}
On the other hand, if $x\in \left(\cup_{k\in\N}H_k\right)\cap H_n^\perp$ then $AP_nx=0=P_nAx$. This implies that $P_n$ and $A$ commutes on $\cup_{k\in\N}H_k$. Let  $\{g_k:k\in\N\}$ be the orthonormal basis of $H$ introduced in Remark \ref{remarkhyp1}\rm(iv). By taking \eqref{comm_A_P_n_H_n} into account, for every $x\in D(A)$ we get
\begin{align*}
P_nAx&=
\sum_{k\in\N}\langle x,g_k\rangle P_n(Ag_k)
= \sum_{k=1}^{s_n}\langle x,g_k\rangle P_n(Ag_k)\\
&= \sum_{k=1}^{s_n}\langle x,g_k\rangle A(P_ng_k)
=\sum_{k\in\N}\langle x,g_k\rangle A(P_ng_k)
=AP_nx,
\end{align*}
and so \eqref{P1} is proved.

We now show that \eqref{P2} and \eqref{P3} hold true. Let $n\in\N$ and let $\lambda\in \rho(A)$. For every $y\in H_n$ there exists a unique $x\in {\rm Dom}(A)$ such that $\lambda x-Ax=y$. Applying $P_n$ to both the sides of this equation, we get
\begin{align*}
y=P_n\lambda x-P_n A x=\lambda P_nx-AP_nx=(\lambda\Id-A)P_nx,    
\end{align*}
which means that $R(\lambda,A)P_nx=R(\lambda,A)x$. The injectivity of the resolvent implies that $P_nx=x$, i.e., $x\in H_n$. We have so proved that for every $n\in\N$ and $y\in H_n$, there exists a unique $x\in H_n$ such that $R(\lambda,A)y=x$, i.e.,
\[
y=\lambda x-Ax=\lambda x-AP_nx=\lambda x-A_nx,
\]
which means $\lambda\in\rho(A_n)$ and \eqref{P3} is verified. Moreover, for every $n\in\N$, if $y\in H_n$ then
\begin{align*}
P_nR(\lambda,A)y=P_nR(\lambda,A)P_ny=R(\lambda,A)P_ny,
\end{align*}
while if $y\in\cup_{k\in\N}H_k\cap H_n^\perp$ then $P_ny=0$ and from \eqref{P1} we deduce that
\begin{align*}
y=\lambda R(\lambda,A)y-AR(\lambda,A)y\Rightarrow 0=P_ny=\lambda P_nR(\lambda,A)y-AP_nR(\lambda,A)y.
\end{align*}
This gives $P_nR(\lambda,A)y=0=R(\lambda,A)P_ny$. Hence, for every $n\in\N$ we infer that $R(\lambda,A)$ commutes with $P_n$ on $\cup_{k\in\N}H_k$ and, from the density of $\cup_{k\in\N}H_k$ in $H$, we obtain \eqref{P2}. Formulae \eqref{P2} and \eqref{P3} yield \eqref{P4} and \eqref{P5} by means of the representation of $e^{tA}$ via resolvent.

Finally, we prove \eqref{stima_conv_fourier_part_n}. Since $A$ is the infinitesimal generator of a strongly continuous and analytic semigroup then it verifies \eqref{stima_risolvente} and \eqref{stima_ris_per_fourier_1} with some constant $c>0$ and $\omega\in\R$. Hence, fixed $\zeta>\omega$, by \eqref{P2} and \eqref{P3}, for every $\lambda\in \{\lambda\in\mathbb C:{\rm Re}\lambda\geq \zeta\}$ we have
\begin{align*}
\|R(\lambda,A_n)y\|_{H_n}
= & \|P_n R(\lambda, A)y\|_{H}
\leq c_{\zeta,1}\|y\|_H, \quad 
\|A_n R(\lambda,A_n)y\|_{H_n}
= \|P_nA R(\lambda, A)y\|_{H}
\leq c_{\zeta,2}\|y\|_H,
\end{align*}
where $c_{\zeta,1}>0$ and $c_{\zeta,2}>0$ are the constant given by \eqref{stima_ris_per_fourier_1} with $\mathcal{C}=A$. From Proposition \ref{lemm:fourier}, it follows that $S_nf\in  L^2(0,T;D(A))$ and 
\begin{align}\label{stima_conv_fourier_part_nnnn}
\left(\int_0^T\norm{A_nS_nf(t)}^2_{H_n}dt\right)^{1/2}\leq \|S_nf\|_{L^2(0,T;D(A))}\leq 2\pi(c_{\zeta,1}+c_{\zeta,2})e^{2|\zeta| T}\|f\|_{L^2(0,T;H)}.   
\end{align}

In particular, for every $\zeta>\omega$, estimate \eqref{stima_conv_fourier_part_nnnn} is independent of $n\in\N$ and estimate \eqref{stima_conv_fourier_part_n} follows at once. 
\end{proof}

Let $(X,W)$ be a weak solution to \eqref{eqFO} (see Definition \ref{weak-solution}). 
For every $n\in\N$ and $x\in H$ we consider the $H_n$-valued stochastic process $\{X_{n}(t,x)\}_{t\in[0,T]}$, which for every $t\in[0,T]$ satisfies
\begin{equation}
\label{approsimazione-mild}
X_n(t,x):=e^{tA_n} P_nx+\int^t_0e^{(t-s)A_n}B_n(X(s,x))ds+W_{A,n}(t), \qquad  \mathbb P{\rm-a.s.},
\end{equation}
where $\{W_{A,n}(t)\}_{t\geq 0}$ is the stochastic process defined as
\begin{align*}
W_{A,n}(t):=\int^t_0 e^{(t-s)A_n}G_ndW(s),
\qquad \mathbb P{\rm-a.s.}
\end{align*}
for every $t\geq0$. We recall that $\{X_{n}(t,x)\}_{t\in[0,T]}$ solves \eqref{SDEn-intro} in the classical It\^o sense.
\begin{lemma}
Assume that Hypotheses \ref{hyp:finito-dimensionale} hold true. Therefore, for every fixed $T>0$ and $x\in H$ we have
\begin{equation*}
\lim_{n \to \infty}\sup_{t\in [0,T]}\mathbb{E}\left[\|X_{n}(t,x)-X(t,x)\|_H^2\right]dt=0.
\end{equation*} 
\end{lemma}

\begin{proof}
Fix $T>0$ and $x\in H$. We begin by proving that 
\begin{align}\label{convsupWA}
\lim_{n\rightarrow \infty}\sup_{t\in [0,T]}\mathbb{E}\left[\|W_A(t)-W_{A,n}(t)\|_H^2\right]=0.
\end{align}
From the definition of $W_{A,n}$ and \eqref{P5}, for every $t\geq0$ we can write
\begin{align*}
W_A(t)-W_{A,n}(t)=\int_0^t e^{(t-s)A}\left(G-P_nG\right)dW(s), \qquad \mathbb P{\rm-a.s.}
\end{align*}
Let $\{u_n:n\in\N\}$ be an orthonormal basis of $U$. Therefore, for every $t\in[0,T]$ and every $n\in\N$ we get
\begin{align*}
\mathbb{E}\left[\|W_A(t)-W_{A,n}(t)\|^2_H\right]
&=\int_0^t \left\|e^{(t-s)A}\left(G-P_nG\right)\right\|_{\mathcal{L}_2(U;H)}^2ds \\
& =\int_0^t \left\|e^{(t-s)A}\left(I-P_n\right)G\right\|_{\mathcal{L}_2(U;H)}^2ds \\
& = \int_0^t\sum_{k=1}^\infty \|e^{(t-s)A}(I-P_n)Gu_k\|_H^2ds \\
& = \int_0^t\sum_{k=1}^\infty \|e^{sA}(I-P_n)Gu_k\|_H^2ds \\
& \leq \int_0^T\sum_{k=1}^\infty \|e^{sA}(I-P_n)Gu_k\|_H^2ds.
\end{align*}
Let us prove that, for a.s $s\in(0,T]$, the series under the integral sign vanishes as $n$ goes to $\infty$. Let us fix $k\in\N$. It follows that
\begin{align*}
\|e^{sA}(I-P_n)Gu_k\|_H^2
= & \langle e^{sA}(I-P_n)Gu_k,e^{sA}(I-P_n)Gu_k\rangle_H\to 0, \qquad n\to \infty,
\end{align*}
since $(I-P_n)h$ converges to $0$ in $H$ as $n$ tends to $\infty$ for every $h\in H$. Further, by \eqref{P4} we get
\begin{align}
\label{stima_pun_serie_conv_stoc}
\|e^{sA}(I-P_n)Gu_k\|_H^2\leq 2\|e^{sA}Gu_k\|_H^2+2\|e^{sA}P_n(Gu_k)\|_H^2
\leq 4\|e^{sA}Gu_k\|_H^2.
\end{align}
From Hypotheses \ref{hyp:finito-dimensionale}(iv), we deduce that for a.s. $s\in(0,T]$ we have
\begin{align*}
\|e^{sA}G\|_{\mathcal L_2(U;H)}^2
= {\rm Trace}_H\left[e^{sA}GG^*e^{sA^*}\right]<\infty.    
\end{align*}
If we apply the dominated convergence theorem with respect to the counting measure, we infer that
\begin{align*}
\sum_{k=1}^\infty \|e^{sA}(I-P_n)Gu_k\|_H^2\to 0, \qquad n\to\infty, \ \textrm{a.s.}\in(0,T].   
\end{align*}
To prove that the whole integral tends to $0$ as $n$ goes to $\infty$, we apply once again the dominated convergence theorem. We have already shown that the function under the integral sign pointwise a.s. converges to $0$ as $n$ diverges to $\infty$. Moreover, arguing as for \eqref{stima_pun_serie_conv_stoc} we deduce that
\begin{align*}
\sum_{k=1}^\infty \|e^{sA}(I-P_n)Gu_k\|_H^2
\leq \sum_{k=1}^\infty2\left(\|e^{sA}Gu_k\|_H^2+\|e^{sA}P_n(Gu_k)\|_H^2\right)
\leq 4\sum_{k=1}^\infty\|e^{sA}Gu_k\|_H^2
\end{align*}
for every $s\in[0,T]$. From Hypotheses \ref{hyp:finito-dimensionale}(iv) it follows that
\begin{align*}
s\mapsto 4\sum_{k=1}^\infty\|e^{sA}Gu_k\|_H^2\in L^1(0,T),
\end{align*}
which implies that
\begin{align}\label{con-ind}
\int_0^T\sum_{k=1}^\infty \|e^{sA}(I-P_n)Gu_k\|_H^2ds\to 0, \qquad n\to\infty.    \end{align}
Since the convergence \eqref{con-ind} is independent of $t$, we obtain \eqref{convsupWA}.

By Hypotheses \ref{hyp:finito-dimensionale} and formulae \eqref{coefficienti-n} and \eqref{P5}, for every $t\in [0,T]$ we get
\begin{align*}
&\norm{\int_0^te^{(t-s)A_n}B_n(X(s,x))ds-\int_0^te^{(t-s)A}B(X(s,x))ds}^2_H\\
&\phantom{aaaaaaaaaa}\qquad\leq\sup_{t\in[0,T]}\norm{e^{tA}}^2_{\mathcal{L}(H)}\int_0^T\norm{B_n(X(s,x))-B(X(s,x))}^2_Hds,\quad \mathbb{P}-{\rm a.s.}
\end{align*}
Since $\norm{B_n}_\infty\leq \norm{B}_\infty<\infty$, by applying the dominated convergence theorem we infer
\begin{equation}\label{convsupint}
    \lim_{n\rightarrow \infty}\sup_{t\in [0,T]}\mathbb{E}\left[\norm{\int_0^te^{(t-s)A_n}B_n(X(s,x))ds-\int_0^te^{(t-s)A}B(X(s,x))ds}^2_H\right]=0.
\end{equation}

Finally, by \eqref{P5} we have 
\[
\norm{e^{tA_n} P_nx-e^{tA}x}_H=\norm{e^{tA}(P_nx-x)}_H\leq \sup_{t\in [0,T]}\norm{e^{tA}}_{\mathcal{L}(H)}\norm{P_nx-x}_H.
\]
Hence, \eqref{convsupWA} and \eqref{convsupint} yield the statement.
\end{proof}
The previous lemma implies the following convergence result.
\begin{proposition}
Assume that Hypotheses \ref{hyp:finito-dimensionale} hold true. For every fixed $T>0$ and $x\in H$ we have
\begin{equation}\label{Convergenza-mild}
\lim_{n \to \infty}\int^T_0\mathbb{E}\left[\|X_{n}(t,x)-X(t,x)\|_H^2\right]dt=0.
\end{equation} 
\end{proposition}

Let $n\in\N$. We consider the linear version of the SPDE \eqref{SDEn-intro}, i.e., when $B_n=0$, given by 
\begin{align}\label{eqFOLn}
\left\{
\begin{array}{ll}
\displaystyle  Z_n(t)=A_nZ_n(t)dt+G_ndW(t), \quad t>0, \vspace{1mm} \\
Z_n(0)=z\in H_n.
\end{array}
\right.
\end{align}
The unique mild solution $\{Z_n(t,z)\}_{t\geq 0}$ to \eqref{eqFOLn} is the $H_n$-valued stochastic process which for every $t\geq0$ enjoys 
\begin{equation}
 Z_n(t,z)=e^{tA_n}z+\int^t_0e^{(t-s)A_n}G_ndW(t) \qquad \mathbb P{\rm-a.s.}
\end{equation}
We define the Markov transition semigroup $\{\mathcal{R}_n(t)\}_{t\geq 0}$ on $B_b(H_n;H_n)$ as
\begin{equation}\label{SOUn}
(\mathcal{R}_n(t)\Phi)(x)=\mathbb{E}\left[\Phi(Z_n(t,x))\right],\quad \Phi\in B_b(H_n;H_n),\; t>0,\; x\in H_n.
\end{equation}
We recall that for every $\Phi\in B_b(H_n;H_n)$, every $t>0$ and every $v,h,k,x\in H_n$ we have
\begin{align}
   &\scal{(\mathcal{R}_n(t)\Phi)(x)}{v}_H=(R_n(t)\phi_v)(x),\label{vet-real}
\end{align}
where $\{R_n(t)\}_{t\geq 0}$ is given by
\[
(R_n(t)\varphi)(x)=\int_{H_n}\varphi(e^{tA_n}x+y)\mu_{t,n}(dy),\quad t>0,\;\varphi\in B_b(H_n),\; x\in H_n,
\]
$\mu_{t,n}$ is the Gaussian measure on $\mathcal{B}(H_n)$ with mean $0$ and covariance operator 
\begin{align}\label{Qtn}
Q_{t,n}:=\int_0^te^{sA_n}G_nG^*_ne^{sA_n^*}ds, \qquad n\in\N, \ t\geq0
\end{align}
and $\phi_v(x)=\langle \Phi(x),v\rangle_H$ for every $x\in H_n$. Now we show that the infinite-dimensional control assumption \eqref{supercontron} implies an $n$-dimensional control assumption. To this aim, we set
\[
\Gamma_{t,n}=Q_{t,n}^{-\frac12}e^{tA_n},\quad t>0,\; n\in\N.
\]

\begin{proposition}
\label{prop:eq_ipotesi}
Conditions \eqref{contron} and \eqref{supercontron} are equivalent to
\begin{align}
&e^{tA_n}(H_n)\subseteq Q_{t,n}^{\frac12}(H_n),\qquad n\in\N,\label{Qtn-iniettivo}\\
&\int_0^t\sup_{n\in\N}\norm{\Gamma_{s,n}}_{\mathcal{L}(H_n)}^{1-\theta}\norm{\Gamma_{s,n}\widetilde{G}}_{\mathcal{L}(U;H_n)}ds<\infty.\label{Cesistenza}
\end{align} 
In particular for every $n\in\N$ and $t\geq 0$ we have
\begin{equation}\label{StimaUGamma}
    \norm{\Gamma_{t,n}}_{\mathcal{L}(H_n)}\leq \norm{\Gamma_{t}}_{\mathcal{L}(H)},\qquad\norm{\Gamma_{t,n}\widetilde{G}}_{\mathcal{L}(U;H_n)}\leq \norm{\Gamma_{t}\widetilde{G}}_{\mathcal{L}(U;H)}.
\end{equation}
\end{proposition}
\begin{proof}
From \eqref{coefficienti-n}, for every $t\geq0$ and $n\in\N$ we get
\begin{align*}
\|Q_{t,n}^{\frac12}x\|_H^2
= & \langle \int_0^t e^{sA_n}G_nG^*_n e^{sA^*_n}ds \ x,x\rangle_H
= \langle \int_0^t e^{sA}GG^* e^{sA^*}ds \ P_nx,P_nx\rangle_H 
=  \|Q_t^{\frac12}P_nx\|_H^2, \quad x\in H.
\end{align*}
By applying \cite[Chapter 3, Corollary 2.3]{Zab08} with $F=Q_{t,n}^{\frac12}$, $G=P_nQ_t^{\frac12}$ and $c=1$, it follows that $Q_{t,n}^{\frac12}(H)=P_nQ_{t}^{\frac12}(H)\subseteq Q_t^{\frac12}(H)$. Hence, for every $x\in Q_{t,n}^{\frac12}(H)$ we infer that $\|Q_{t,n}^{-\frac12}x\|_H=\|Q_{t}^{-\frac12}P_n^{-1}x\|_H$. 

Now, we prove that \eqref{contron} and \eqref{supercontron} imply \eqref{Qtn-iniettivo} and \eqref{Cesistenza}. Recalling that, for every $x\in P_n(H)$, $P_n^{-1}x=x$ as element of $H$, by \eqref{coefficienti-n}, for every $n\in\N$ we get $e^{tA_n}(H)=P_ne^{tA}(H)\subseteq P_nQ_t^{\frac12}(H)=Q_{t,n}^{\frac12}(H)\subseteq Q_t^{\frac12}(H)$. Since $e^{tA_n}(H)\subseteq H_n=P_n(H)$, it follows that
\begin{align}
\label{ug_norme_eq_cond_contr}
\|Q_t^{-\frac12}e^{tA_n}y\|_H
= & \|Q_t^{-\frac12}P_n^{-1}e^{tA_n}y\|_H
= \|Q_{t,n}^{-\frac12}e^{tA_n}y\|_H
\end{align}
for every $t>0$, every $n\in\N$ and every $y\in H$.
From \eqref{ug_norme_eq_cond_contr}, we infer that
\begin{align*}
\sup_{y\in H_n}\|Q_{t,n}^{-\frac12}e^{tA_n}y\|_{H_n}
=\sup_{y\in H_n}\|Q_{t}^{-\frac12}e^{tA}y\|_{H}\leq \norm{\Gamma_t}_{\mathcal{L}(H)}\|y\|_{H}, \qquad t\in(0,T), \ n\in\N.    
\end{align*}
In particular, $\Gamma_{t,n}$ is well-defined and satisfies $\norm{\Gamma_{t,n}}_{\mathcal{L}(H_n)}\leq \norm{\Gamma_t}_{\mathcal{L}(H)}$ for every $t\in(0,T)$ and every $n\in\N$. By similar arguments it follows that $\norm{\Gamma_{t,n}\widetilde{G}_n}_{\mathcal{L}(U;H_n)}\leq \norm{\Gamma_t\widetilde{G}}_{\mathcal{L}(U;H)}$, and so \eqref{contron} and \eqref{supercontron} imply \eqref{Qtn-iniettivo} and \eqref{Cesistenza}.

Assume that \eqref{Qtn-iniettivo} and \eqref{Cesistenza} are fulfilled. Therefore, $e^{tA_n}(H)\subseteq Q_{t,n}^{\frac12}(H)=P_nQ_t^{\frac12}(H)\subseteq Q_t^{\frac12}(H)$, the operator $Q_t^{-\frac12}e^{tA_n}$ belong to $\mathcal L(H)$ for every $n\in\N$ and every $t>0$ and \eqref{ug_norme_eq_cond_contr} holds true also in this case.\\
Fix $t>0$, $n\leq m$, $n,m\in\N$ and $y\in H_n$. From \eqref{ug_norme_eq_cond_contr} and the fact that $e^{tA_m}(H_n)=e^{tA_n}(H_n)$, we obtain
\begin{align*}
\|Q_{t,n}^{-\frac12}e^{tA_n}y\|_H
= \|Q_t^{-\frac12}e^{tA_n}y\|_H
= \|Q_t^{-\frac12}e^{tA_m}y\|_H
=\|Q_{t,m}^{-\frac12}e^{tA_m}y\|_H.
\end{align*}
Taking the supremum with respect to $y\in H_n\subseteq H_m$ with $\|y\|_{H_n}=\|y\|_H\leq 1$, we infer that $\norm{\Gamma_{t,n}}_{\mathcal{L}(H_n)}\leq \norm{\Gamma_{t,m}}_{\mathcal{L}(H_m)}$ for every $t>0$ and every $n\leq m$ with $m,n\in\N$. It can be analogously proved that $\norm{\Gamma_{t,n}\widetilde G_n}_{\mathcal{L}(U;H_n)}\leq \norm{\Gamma_{t,m}\widetilde G_m}_{\mathcal{L}(U;H_m)}$ for every $m,n\in\N$ with $n\leq m$. By the monotone convergence theorem and \eqref{Cesistenza}, we deduce that the function
\begin{align*}
L_t:=\lim_{n\to\infty}\norm{\Gamma_{t,n}}_{\mathcal{L}(H_n)}^{1-\theta}\norm{\Gamma_{t,n}\widetilde G_n}_{\mathcal{L}(U;H_n)}, \qquad t>0,    
\end{align*}
belongs to $L^{1}(0,T)$. We claim that $e^{tA}(H)\subseteq Q_t^{\frac12}(H)$ for every $t\in(0,T)$ such that $L_t<\infty$. Assume by contradiction that the claim is false. Therefore, the results \cite[Chapter 3, Section 2]{Zab08} and the density of $\cup_{n\in\N}H_n$ in $H$ imply that for every $k\in\N$ there exist $m_k\in\N$ and $x_k\in H_{m_k}$, with $\|x_k\|_H= 1$, such that $ \|Q_t^{-\frac12}e^{tA}x_k\|_H= \|Q_t^{-\frac12}e^{tA_{m_k}}x_k\|_H\geq k$. Without loss of generality we assume that $(m_k)$ is an increasing sequence. Since $\norm{\Gamma_{t,m_k}}_{\mathcal{L}(H_{{m_k}})}=\|Q_t^{-\frac12}e^{tA_{m_k}}\|_{\mathcal L(X)}\geq \|Q_t^{-\frac12}e^{tA_{m_k}}x_k\|_H$, letting $k$ go to infinity, we deduce that $L_t=\infty$, and this occurs only for $t$ which belongs to a subset of $(0,T)$ of null measure. 
To show that $e^{tA}(H)\subseteq Q_t^{\frac12}(H)$ for every $t\in(0,T)$, we notice that, if $\overline t\in(0,T)$ satisfies $\Gamma_{\overline t}=\infty$, then there exists $\hat t\in(0,\overline t)$ such that $e^{\hat tA}(H)\subseteq Q_{\hat t}^{\frac12}(H)$. Hence,
\begin{align*}
e^{\overline tA}(H)=e^{\hat tA}e^{(\overline t-\hat t)A}(H)\subseteq e^{\hat tA}(H)\subseteq Q_{\hat t}^{\frac12}(H)\subseteq Q_{\overline t}^{\frac12}(H).    
\end{align*}
It remains to prove that $t\mapsto \|\Gamma_t\|_{\mathcal L(H)}^{1-\theta}\|\Gamma_t\widetilde G\|_{\mathcal L(U;H)}\in L^{1}(0,T)$. Fix $t\in (0,T)$. The density of $\bigcup_{n\in\N}H_n$ in $H$ gives
\begin{align*}
\|Q_t^{-\frac12}e^{tA}\|_{\mathcal L(H)}
=\sup_{x\in\cup_{n\in\N}H_n, \|x\|_H\leq 1}\|Q_t^{-\frac12}e^{tA}x\|_H.
\end{align*}
We notice that, if $x\in H_n$, then
\begin{align*}
Q_t^{-\frac12}e^{tA}x=y \Longleftrightarrow e^{tA}x=Q^{\frac12}_t y \Longleftrightarrow 
e^{tA_n}x=Q^{\frac12}_t y \Longleftrightarrow
Q_{t,n}^{\frac12}e^{tA_n}x=Q_t^{-\frac12}e^{tA_n}x=y.
\end{align*}
Hence, for every $n\in\N$ and $x\in H_n$, it follows that $\|Q_t^{-\frac12}e^{tA}x\|_H= \|\Gamma_{t,n}x\|_{H_n}$. Further,
\begin{align*}
\Gamma_t \widetilde Gu
= & \lim_{m\to\infty} Q_t^{-\frac12}e^{tA}P_m\widetilde Gu =  \lim_{m\to\infty} Q_{t,m}^{-\frac12}e^{tA_m}\widetilde G_mu=\lim_{m\to\infty}\Gamma_{t,m}\widetilde G_mu, \qquad u\in U.
\end{align*}
This means that, for every $n\in\N$, every $x\in H_n$ and every $u\in U$, the fact that $\|\Gamma_{t,n}x\|_{H_n}=\|\Gamma_{t,m}x\|_{H_m}$ for every $m\geq n$ gives
\begin{align*}
\|\Gamma_tx\|_{H}^{1-\theta}\|\Gamma_t\widetilde Gu\|_{H}
= &\|\Gamma_{t,n}x\|_{H_n}^{1-\theta}\lim_{m\to\infty}\|\Gamma_{t,m}\widetilde G_mu\|_{H_m}
= \lim_{m\to\infty}\|\Gamma_{t,m}x\|_{H_m}^{1-\theta}\|\Gamma_{t,m}\widetilde G_mu\|_{H_m} \\
\leq & L_t\|x\|_{H}\|u\|_U.
\end{align*}
Therefore, $\|\Gamma_t\|_{\mathcal L(H)}\|\Gamma_t\widetilde G\|_{\mathcal L(U;H)}\leq L_t$ for every $t\in(0,T)$, which concludes the proof.
\end{proof}

Let $n\in\N$ and let $F\in C^\theta_b(H;H)$, where $\theta$ is defined in Hypotheses \ref{hyp:finito-dimensionale}. We set $F_n(\cdot)=P_nF(P_n(\cdot))$ and consider the backward integral equation
\begin{equation}\label{Back-Kolmon}
    U_n(t,x)=\int^T_t \mathcal{R}_n(r-t)\left( DU_n(r,\cdot)\widetilde{G}_n\widetilde{B}_n(\cdot)+F_n(\cdot)\right)(x)dr ,\quad t\in [0,T], \ x\in H_n.
\end{equation}
\begin{proposition}\label{SolBKn}
Assume that Hypotheses \ref{hyp:finito-dimensionale} hold true. For every $n\in\N$, equation \eqref{Back-Kolmon} admits a unique solution $U_n\in C_b^{0,1}([0,T]\times H_n;H_n)$ such that the map $x\rightarrow DU_n(t,x)\widetilde{G}_n$ belongs to $C^1_b(H_n;\mathcal{L}(U;H_n))$ for every $t\in [0,T]$. Moreover, for every $n\in\N$ and $t\in [0,T]$ we have
\begin{align}\label{stima-n}
\sup_{t\in[0,T]}\left(\norm{U_n(t,\cdot)}_{C^1_b(H_n;H_n)}+\|DU_n(t,\cdot)\widetilde{G}_n\|_{C^1_b(H_n;\mathcal{L}(U;H_n))}\right)&\leq M_{T} \|F_n\|_{C^\theta_b(H_n;H_n)},
 \end{align}
where $M_T$ is a positive constant such that
\[
\lim_{T\rightarrow0}M_T=0.
\]
\end{proposition}
\begin{proof}
By Proposition \ref{WPK}, with $\X=H_n$, $E=U$, $\mathcal N=\widetilde B_n$ and $I=\widetilde G_n$, for every $n\in\N$ equation \eqref{Back-Kolmon} has a unique solution $U_n\in C_b^{0,1}([0,T]\times H_n;H_n)$ such that for every $t\in [0,T]$ we have
\begin{align*}
\sup_{t\in[0,T]}\left(\norm{U_n(t,\cdot)}_{C^1_b(H_n;H_n)}+\|DU_n(t,\cdot)\widetilde{G}_n\|_{C^1_b(H_n,\mathcal{L}(U,H_n))}\right)\leq M_{T,n} ||F_n||_{C^\theta_b(H_n;H_n)},
 \end{align*}
 where
\begin{align*}
&M_{T,n}:=C_{T,n}e^{C_{T,n}\|\widetilde{B}_n\|_{C^\theta_b(H_n;H_n)}},\quad  C_{T,n}:=\int^T_0K_{t,n}dt,\\
& K_{t,n}:=1+\left(1+\norm{\widetilde{G}_n}_{\mathcal{L}(U;H_n)}+\|\Gamma_{t,n}\widetilde{G}_n\|_{\mathcal{L}(U;H_n)}\right)\|e^{tA_n}\|^\theta_{\mathcal{L}(H_n;H_n)}\|\Gamma_{t,n}\|^{1-\theta}_{\mathcal{L}(H_n;H_n)}.
\end{align*}
By \eqref{coefficienti-n}, \eqref{StimaUGamma} and Lemma \ref{Lemma:semi-n}, for every $n\in\N$ we get
\[
C_{T,n}\leq \int^T_0 1+\left(1+\norm{\widetilde{G}}_{\mathcal{L}(U;H)}+\|\Gamma_{t}\widetilde{G}\|_{\mathcal{L}(U;H)}\right)\|e^{tA}\|^\theta_{\mathcal{L}(H;H)}\|\Gamma_{t}\|^{1-\theta}_{\mathcal{L}(H;H)}dt
\]
By \eqref{supercontron} we obtain the statement.
\end{proof}

For every $n\in\N$, let $U_n$ be the unique solution to \eqref{Back-Kolmon} given by Proposition \ref{SolBKn}. For every $v\in H_n$ we define the function $U^v_{n}=\scal{U_n}{v}_H$.
We note that by \eqref{deri1vet}, if $v,w\in H_n$ then
\begin{equation}\label{scambio}
\scal{DU_n(x)\widetilde{G}_nv}{w}=\scal{\widetilde{G}_n^*\nabla U^w_n(x)}{v},\quad x\in H_n,\; n\in\N.
\end{equation}
From \eqref{vet-real} and \eqref{scambio}, $U^v_{n}$ belongs to $C^{0,1}_b([0,T]\times H_n)$ and solves the integral equation
\begin{equation}\label{Back-Kolmonk}
  U^v_{n}(t,x)=\int^T_t R_n(r-t)\left( \scal{\widetilde{G}_n^*\nabla U^v_{n}(r,\cdot)}{\widetilde{B}_n(\cdot)}_H+\scal{F_{n}}{v}_H\right)(x)dr ,\quad t\in [0,T], \ x\in H_n.
\end{equation}
Arguing as in the proof of Proposition \ref{SolBKn}, we infer that
\begin{align}\label{stima-n-scalare}
\sup_{t\in[0,T]}\left(\norm{U^v_n(t,\cdot)}_{C^1_b(H_n)}+\|\widetilde{G}_n^*\nabla U_n(t,\cdot)\|_{C^1_b(H_n;U)}\right)\leq M_T\|\langle F_n,v\rangle_H\|_{C^\theta_b(H_n)},\quad
t\in[0,T].
\end{align}
In particular, if $\{g_n : n\in\N\}$ is the orthonormal basis of $H$ introduced in Remark \ref{remarkhyp1}(iv), then the functions $U_{n,k}=\scal{U_n}{g_k}_H$, with $k=1,..., s_n$ verify
\begin{equation}\label{serieUk}
U_n=\sum^{s_n}_{k=1}U_{n,k}g_k,
\end{equation}
where $s_n={\rm Dim}(H_n)$. 
\begin{rmk}\label{2-derivate}
If $U=H$ and $\widetilde{G}=\Id_H$ (see the example in Subsection \ref{Heat-case}), then the solution $U_n$ to \eqref{Back-Kolmon} belongs to $C^{0,2}_b([0,T]\times H_n;H_n)$. Moreover, for every $v\in H_n$ the function $U^v_{n}=\scal{U_n}{v}_H$ is a classical solution to 
{\small
\begin{align*}
\left\{
\begin{array}{ll}
\displaystyle \frac{\partial U^v_{n}(t,x)}{\partial t}+\frac{1}{2}{\rm Trace}[\nabla^2U^v_{n}(t,x)]+\langle A_nx,\nabla U^v_{n}(t,x)\rangle_H+\langle B_n(x),\nabla U
^v_{n}(r,x)\rangle_H +\scal{F_n(x)}{v}_H=0, \\
U^v_{n}(T,x)=0.
\end{array}
\right.
\end{align*}
}
We underline that in the general case this fact is not true.
\end{rmk}

\subsection{Proof of Theorem \ref{Ito-Tanaka} }\label{dim-Ito}
As in the previous subsection, we fix $T>0$ and a weak solution $(X,W)$ to \eqref{eqFO}. We now prove Theorem \ref{Ito-Tanaka}.

\begin{proof}[Proof of Theorem \ref{Ito-Tanaka}]

Fix $n\in\N$. Let $U_n$ be the solution of \eqref{Back-Kolmon} with $F_n=B_n$ given by Proposition \ref{SolBKn}. Let $\{g_1,\ldots,g_{s_n}\}$ be the orthonormal basis of $H_n$ introduced in Remark \ref{remarkhyp1}(iv) with $s_n={\rm Dim}(H_n)$. The function $U_{n,k}=\scal{U_n}{g_k}_H$ is the solution to \eqref{Back-Kolmonk} with $F_n=B_n$, $v=g_k$ and $k=1,...,s_n$. Fix $k\in \{1,\ldots,s_n\}$. We set
\[
f(r,x):=\scal{\widetilde{G}_n^*\nabla U_{n,k}(r,x)}{\widetilde{B}_n(x)}_H+\scal{B_n(x)}{g_k}_H, \quad r\in[0,T],\; x\in H_n. 
\]
By Proposition \ref{SolBKn}, the map $x\rightarrow f(r,x)$ belongs to $C^\theta_b(H_n)$ for every $r\in [0,T]$. Let $(f_h)_{h\in\N}\subseteq C^{0,2}_b([0,T]\times H_n)$ be the sequence introduced in Theorem \ref{OU-finito}. From such a theorem, the function
\[
U_{n,k,h}(t,x)=\int_t^TR_n(r-t)f_h(r,x)dr, \qquad t\in[0,T],
\]
is the strict solution to the parabolic equation 

\begin{align}\label{Parabolica-n}
\left\{
\begin{array}{ll}
\displaystyle \frac{\partial u(t,x)}{\partial t}+\frac{1}{2}{\rm Trace}\left[G_nG_n^*\nabla^2u(t,x)\right]+\scal{A_nx}{\nabla u(t,x)}+f_h(t,x)=0, \qquad t\in(0,T], \ x\in H_n, \vspace{3mm} \\
u(T,x)=0,\qquad  x\in \R^n.
\end{array}
\right.
\end{align}
Let $\{X_n(t,x)\}_{t\in [0,T]}$ be the process defined by \eqref{approsimazione-mild}. By the It\^o formula, for every $t\in[0,T]$ we get
\begin{align*}
dU_{n,k,h}(t,X_n(t,x))&=\frac{\partial U_{n,k,h}}{\partial t}(t,X_n(t,x))dt\\
&+ \scal{\nabla U_{n,k,h}(t,X_n(t,x))}{A_nX_n(t,x)+B_n(X(t,x))}dt\\
&+\frac{1}{2}{\rm Trace}[G_nG_n^*\nabla^2U_{n,k,h}(t,X_n(t,x))]dt + \scal{\nabla U_{n,k,h}(t,X_n(t,x))}{G_ndW(t)}
\end{align*}
which, combined with \eqref{Parabolica-n}, for every $t\in[0,T]$ gives
\begin{align*}
dU_{n,k,h}(t,X_n(t,x))
= & \scal{\nabla U_{n,k,h}(t,X_n(t,x))}{B_n(X(t,x))}dt-f_h(t,X_n(t,x))dt \\
&+ \scal{\nabla U_{n,k,h}(t,X_n(t,x)) }{G_ndW(t)}, \qquad \mathbb P\textup{-a.s.}
\end{align*}
Therefore, for every $t\in[0,T]$ it holds that
\begin{align*}
U_{n,k,h}(t,X_n(t,x))-U_{n,k,h}(0,P_nx)
= & \int_0^t\scal{\nabla U_{n,k,h}(s,X_n(s,x))}{B_n(X(s,x))}ds-\int_0^tf_h(s,X_n(s,x))ds \\
&+ \int_0^t\scal{\nabla U_{n,k,h}(s,X_n(s,x)) }{G_ndW(s)}, \qquad \mathbb P\textup{-a.s.}
\end{align*}
Letting $h$ tend to infinity, from Theorem \ref{OU-finito} with $u=U_{n,k}$, we infer that, for every $t\in[0,T]$,
\begin{align*}
U_{n,k}(t,X_n(t,x))-U_{n,k}(0,P_nx)
= & \int_0^t\scal{\nabla U_{n,k}(s,X_n(s,x))}{B_n(X(s,x))}ds-\int_0^tf(s,X_n(s,x))ds \\
&+ \int_0^t\scal{\nabla U_{n,k}(s,X_n(s,x)) }{G_ndW(s)}, \qquad \mathbb P\textup{-a.s.},
\end{align*}
from which it follows that, for every $t\in[0,T]$,
\begin{align*}
dU_{n,k}(t,X_n(t,x))&= \scal{\nabla U_{n,k}(t,X_n(t,x))}{B_n(X(t,x))-B_n(X_n(t,x))}dt\\
&- \scal{B_n(X_n(t,x))}{g_k}dt+ \scal{\nabla U_{n,k}(t,X_n(t,x)) }{G_ndW(t)}, \qquad \mathbb P\textup{-a.s.}
\end{align*}
Hence, for every $t\in[0,T]$,
\begin{align}
-\scal{B_n(X_n(t,x))}{g_k}dt
= & dU_{n,k}(t,X_n(t,x))-\scal{\nabla U_{n,k}(t,X_n(t,x))}{B_n(X(t,x))-B_n(X_n(t,x))}dt\notag\\
&-\scal{\nabla U_{n,k}(t,X_n(t,x))}{G_ndW(t)}, \qquad \mathbb P\textup{-a.s.}
\label{Bnk}
\end{align}
Summing up $k$ from $1$ to $s_n$ in both the sides of \eqref{Bnk}, by \eqref{serieUk} (see also \eqref{RderiU} and \eqref{RderiU-var}) we obtain, for every $t\in[0,T]$,
\begin{align}\label{Bn}
-B_n(X_n(t,x))dt
= & dU_{n}(t,X_n(t,x))
-DU_{n}(t,X_n(t,x))\left(B_n(X(t,x))-B_n(X_n(t,x))\right)dt\notag\\
&-D U_n(t,X_n(t,x))G_ndW(t), \qquad \mathbb P\textup{-a.s.}
\end{align}
Adding and subtracting $B_n(X_n(t,x))dt$ in \eqref{SDEn-intro}, from \eqref{Bn} we get, for every $t\in[0,T]$, 
\begin{align*}
dX_n(t)&=A_nX_n(t,x)dt+B_n(X_n(t,x))dt+\left(B_n(X(t,x))-B_n(X_n(t,x))\right)dt+G_ndW(t)\\
&=A_nX_n(t,x)dt-dU_{n}(t,X_n(t,x))dt+DU_{n}(t,X_n(t,x))\left(B_n(X(t,x))-B_n(X_n(t,x))\right)dt\\
&+\left(B_n(X(t,x))-B_n(X_n(t,x))\right)dt+D U_n(t,X_n(t,x))G_ndW(t)+G_ndW(t), \qquad \mathbb P\textup{-a.s.}
\end{align*}
By applying the variation of constants formula, it follows that, for every $t\in[0,T]$,
\begin{align*}
X_n(t)&=e^{tA_n} P_nx-\int^t_0e^{(t-s)A_n}dU_{n}(s,X_n(s,x))+\int^t_0e^{(t-s)A_n} \left(B_n(X(s,x))-B_n(X_n(s,x))\right)ds\\
&+\int^t_0e^{(t-s)A_n}DU_{n}(s,X_n(s,x))\left(B_n(X(s,x))-B_n(X_n(s,x))\right)ds\\
&+\int^t_0e^{(t-s)A_n}DU_n(t,X_n(t,x))G_ndW(s)+ \int^t_0e^{(t-s)A_n}G_ndW(s), \qquad \mathbb P\textup{-a.s.}
\end{align*}
Finally, integrating by parts the first integral we conclude that, for every $t\in[0,T]$,
\begin{align*}
X_n(t)&=e^{tA_n}( P_nx+U_n(0,P_nx))-U_n(t,X_n(t,x))-A_n\int^t_0e^{(t-s)A_n}U_{n}(s,X_n(s,x))ds\\
&+\int^t_0e^{(t-s)A_n} \left(B_n(X(s,x))-B_n(X_n(s,x))\right)ds\\
&+\int^t_0e^{(t-s)A_n}DU_{n}(s,X_n(s,x))\left(B_n(X(s,x))-B_n(X_n(s,x))\right)ds\\
&+\int^t_0e^{(t-s)A_n}DU_n(s,X_n(s,x))G_ndW(s)+ \int^t_0e^{(t-s)A_n}G_ndW(s), \qquad \mathbb P\textup{-a.s.}
\end{align*}
\end{proof}

\subsection{Proof of Theorem \ref{pathwiseuniqueness}}\label{dim-unicità}
Fix $T>0$ and let $(X_1,W)$ and $(X_2,W)$ be two weak solutions to \eqref{eqFO} defined on the same probability space $(\Omega,\mathcal{F},\{\mathcal{F}_t\}_{t\in [0,T]},\mathbb{P})$. Theorem \ref{pathwiseuniqueness} is verified if we prove that
\begin{equation}\label{unicitàL2}
\Delta:=\mathbb{E}\left[\int_0^T\norm{X_1(t,x)-X_2(t,x)}_H^2dt\right]=\int_0^T\mathbb{E}\left[\norm{X_1(t,x)-X_2(t,x)}_H^2\right]dt=0,\quad x\in H.
\end{equation}
Indeed, by Hypotheses \ref{hyp:finito-dimensionale}(iv) and \cite[Theorem 5.11]{Dap-Zab14} the processes $X_1$ and $X_2$ admit a continuous modification. Hence, by \eqref{unicitàL2} there exists $\Omega_0\subseteq\Omega$ such that $\mathbb{P}(\Omega_0)=1$ and 
\[
X_1(t,x)(\omega)=X_2(t,x)(\omega),\quad \forall\; (t,\omega)\in [0,T]\times\Omega_0,
\]
namely $X_1=X_2$ on $L^2([0,T]\times \Omega, \mathcal{B}([0,T])\times \mathcal{F},\lambda\times\mathbb{P})$.
\begin{proof}[Proof of \eqref{unicitàL2}]
Let $x\in H$, $n\in\N$ and let $U_n$ be the solution to \eqref{Back-Kolmon} given by Proposition \ref{SolBKn} with $F_n=B_n$. By Proposition \ref{Ito-Tanaka}, for every $t\in [0,T]$ we have
\begin{equation}\label{unicitaL2n}
\Delta_n:=\int_0^T\mathbb{E}\left[\norm{X_{1,n}(t,x)-X_{2,n}(t,x)}_H^2\right]dt\leq 7(I_1+I_2+\sum_{i=1}^2I_{3,i}+\sum_{i=1}^2I_{4,i}+I_5),
\end{equation}
where $\{X_{1,n}(t,x)\}_{t\in[0,T]}$ and $\{X_{2,n}(t,x)\}_{t\in[0,T]}$ are the processes defined in \eqref{approsimazione-mild}, with $X$ replaced by $X_1$ and $X_2$, respectively, and 
\begin{align}
&I_1:=\int_0^T\mathbb{E}\left[\|U_n(t,X_{1,n}(t,x))-U_n(t,X_{2,n}(t,x))\|^2_H\right]dt,\notag\\
&I_2:=\int_0^T\mathbb{E}\left[\norm{A_n\int^t_0e^{(t-s)A_n}\left(U_{n}(s,X_{1,n}(s,x))-U_{n}(s,X_{2,n}(s,x))\right)ds}_H^2\right]dt,\notag\\
&I_{3,i}:=\int_0^T\mathbb{E}\left[\norm{\int^t_0e^{(t-s)A_n} \left(B_n(X_i(s,x))-B_n(X_{i,n}(s,x))\right)ds}^2_H\right]dt, \qquad i=1,2,\notag\\
&I_{4,i}:=\int_0^T\mathbb{E}\left[\norm{\int^t_0e^{(t-s)A_n} DU_{n}(s,X_{i,n}(s,x))\left(B_n(X_i(s,x))-B_n(X_{i,n}(s,x))\right)ds}^2_H\right]dt \quad i=1,2, \notag\\    &I_5:=\int_0^T\mathbb{E}\left[\norm{\int_0^te^{(t-s)A_n}(DU_n(s,X_{1,n}(s,x))-DU_n(s,X_{2,n}(s,x)))G_ndW(s)}^2_H\right]dt\label{I7}.
\end{align}
Before to estimate the above integrals, we recall that
\begin{equation}\label{holderBn}
    \|B_n\|_{C^\theta_b(H;H_n)}\leq \|B\|_{C^\theta_b(H;H)},\quad  n\in\N.
\end{equation}    
Let us estimate $I_1$. From \eqref{stima-n} and \eqref{holderBn} we get
\begin{equation}\label{SI1}
    I_1\leq M_{T}^2 \|B\|_{C^\theta_b(H;H)}^2\Delta_n.
\end{equation}
To estimate $I_2$, we take advantage of  \eqref{stima_conv_fourier_part_n}, \eqref{stima-n}, \eqref{holderBn} and Fubini--Tonelli's Theorem to infer that
\begin{align}\label{SI2}
I_2&=\mathbb{E}\left[ \int_0^T\norm{A_n\int^t_0e^{(t-s)A_n}\left(U_{n}(s,X_{1,n}(s,x))-U_{n}(s,X_{2,n}(s,x))\right)ds}_H^2dt\right]\notag\\
&\leq \mathbb{E}\left[ C^2_{T,A}\int_0^T\norm{U_{n}(s,X_{1,n}(s,x))-U_{n}(s,X_{2,n}(s,x))}_H^2dt\right]\notag\\
&=C^2_{T,A}\int_0^T\mathbb{E}\left[\norm{U_n(t,X_{1,n}(t,x))-U_n(t,X_{2,n}(t,y))}_H^2\right]dt\notag \\
&\leq C^2_{T,A}M_{T}^2 \|B\|_{C^\theta_b(H;H)}^2\Delta_n.
\end{align}
In the same way by \eqref{P5}, \eqref{stima-n} and \eqref{holderBn} we get
\begin{equation}\label{SI36}
    \sum_{i=1}^2(I_{3,i}+I_{4,i})\leq T\sup_{t\in[0,T]}\norm{e^{tA}}_{\mathcal{L}(H)}^2\left(1+M^2_{T} ||B||_{C^\theta_b(H;H)} ^2\right)||B||^2_{C^\theta_b(H;H)}\Pi_{n},
\end{equation}
where 
\begin{equation*}
\Pi_{n}:=\int_0^T\mathbb{E}\left[\norm{X_{1,n}(s,x)-X_1(s,x)}_H^{2\theta}\right]ds+\int_0^T\mathbb{E}\left[\norm{X_{2,n}(s,x)-X_2(s,x)}_H^{2\theta}\right]ds.
\end{equation*}
It remains to deal with \eqref{I7}. By applying the It\^o isometry, we get
\begin{align}\label{itooo}
I_5 & = \int^T_0\int_0^t\mathbb{E}\left[\norm{e^{(t-s)A_n}(DU_n(s,X_{1,n}(s,x))-DU_n(s,X_{2,n}(s,x)))G_n}_{\mathcal L_2(U;H)}^2\right]dsdt.
\end{align}
If Hypotheses \ref{hyp:traccia-finita} hold true, then from \eqref{P5}, \eqref{stima-n} and \eqref{itooo} we obtain
\begin{align}\label{SI8-traccia}
I_5 & \leq T\sup_{t\in[0,T]} \norm{e^{tA}}_{\mathcal{L}(H)}^2M^2_{T} ||B||^2_{C^\theta_b(H;H)} \norm{\mathcal{V}}_{\mathcal{L}_2(U;U)}\|\widetilde G\|_{\mathcal L(U;H)}\Delta_n,
\end{align}
since, given two operators $A\in\mathcal{L}(U;H)$ and $B\in\mathcal{L}_2(U;U)$, it follows that
\[
\norm{AB}_{\mathcal{L}_2(U;H)}\leq \norm{A}_{\mathcal{L}(U;H)}\norm{B}_{\mathcal{L}_2(U;U)}.
\]
If Hypotheses \ref{hyp:goal-addo} hold true, then the computations are more involving. Let us fix an orthonormal basis $\{u_\ell:\ell\in\N\}$ of $U$. We set 
\[
T_n(s,t):=e^{(t-s)A_n}(DU_n(s,X_{1,n}(s,x))-DU_n(s,X_{2,n}(s,x)))G_n
\]
for every $s,t\in[0,T]$ with $s\leq t$. By Hypotheses \ref{hyp:goal-addo} we get 
\begin{align*}
&\int_0^t\mathbb{E}\left[\norm{T_n}_{L_2(U;H)}^2\right]ds
=\int_0^t\mathbb{E}\left[\sum_{\ell=1}^\infty\norm{T_n u_\ell}_H^2\right]ds
\\= & \int_0^t\mathbb{E}\left[\sum_{\ell=1}^\infty\sum_{k=1}^n\left(\sum_{h=1}^{d_k}\langle T_n u_\ell,e^k_h\rangle_H^2+\sum^{d_k}_{i,j=1,\;i\neq j}\langle T_n u_\ell,e^k_i\rangle_H\langle T_n u_\ell,e^k_j\rangle_H\scal{e^k_i}{e^k_j}_H\right)\right]ds \\
\\ \leq & \int_0^t\mathbb{E}\left[\sum_{\ell=1}^\infty\sum_{k=1}^n\left(\sum_{h=1}^{d_k}\langle T_n u_\ell,e^k_h\rangle_H^2+\frac{1}{2}\sum^{d_k}_{i,j=1,\;i\neq j}\left(\langle T_n u_\ell,e^k_i\rangle_H^2+\langle T_n u_\ell,e^k_j\rangle_H^2\right)\right)\right]ds \\
\\ = & \int_0^t\mathbb{E}\left[\sum_{\ell=1}^\infty\sum_{k=1}^n\left(\sum_{h=1}^{d_k}\langle T_n u_\ell,e^k_h\rangle_H^2+\frac{1}{2}(d_k-1)\left(\sum^{d_k}_{i=1}\langle T_n u_\ell,e^k_i\rangle_H^2+\sum^{d_k}_{j=1}\langle T_n u_\ell,e^k_j\rangle_H^2\right)\right)\right]ds \\
\\ = & \int_0^t\mathbb{E}\left[\sum_{\ell=1}^\infty\sum_{k=1}^n d_k\sum_{i=1}^{d_k}\langle T_n u_\ell,e^k_i\rangle_H^2\right]ds\leq d\int_0^T\mathbb{E}\left[\sum_{\ell=1}^\infty\sum_{k=1}^n \sum_{i=1}^{d_k}\langle T_n u_\ell,e^k_i\rangle_H^2\right]ds \\
= & d\int_0^t\sum_{\ell=1}^\infty\sum_{k=1}^n\sum_{i=1}^{d_k}\mathbb{E}\left[\langle e^{(t-s)A_n}(DU_n(s,X_{1,n}(s,x))-DU_n(s,X_{2,n}(s,x)))G_n u_\ell,e^k_i\rangle_H^2\right]ds \\
= & d\int_0^t\sum_{\ell=1}^\infty\sum_{k=1}^n\sum_{i=1}^{d_k}e^{2(t-s){\rm Re}(\rho^k_{i})}\mathbb{E}\left[\langle (DU_n(s,X_{1,n}(s,x))-DU_n(s,X_{2,n}(s,x))) G_nu_\ell,e^k_i\rangle_H^2\right]ds.
\end{align*}
Setting $U_{n,k,i}:=\scal{U_n}{e^k_i}_H$ and $B^k_i:=\scal{B}{e^k_i}_H$ by \eqref{coefficienti-n}, \eqref{scambio} and \eqref{stima-n-scalare}, we obtain 
\begin{align*}
& \int_0^t\mathbb{E}\left[\norm{T_n(s,t)}_{L_2(U;H)}^2\right]ds \\
\leq &   d\int_0^t\sum_{\ell=1}^\infty\sum_{k=1}^n\sum_{i=1}^{d_k}e^{2(t-s){\rm Re}(\rho^k_{i})}\mathbb{E}\left[\langle \nabla U_{n,k,i}(s,X_{1,n}(s,x))-\nabla U_{n,k,i}(s,X_{2,n}(s,x)),\widetilde{G}_n\mathcal{V} u_\ell\rangle_H^2\right]ds \\
= &   d\int_0^t\sum_{\ell=1}^\infty\sum_{k=1}^n\sum_{i=1}^{d_k}e^{2(t-s){\rm Re}(\rho^k_{i})}\mathbb{E}\left[\langle \mathcal{V}^*\widetilde{G}_n^*(\nabla U_{n,k,i}(s,X_{1,n}(s,x))-\nabla U_{n,k,i}(s,X_{2,n}(s,x))),u_\ell\rangle_U^2\right]ds \\
\leq & d\|\mathcal{V}^*\|_{\mathcal L(U;U)}^2 \int_0^t\sum_{k=1}^n\sum_{i=1}^{d_k}e^{2(t-s){\rm Re}(\rho^k_{i})}\mathbb{E}\left[\|\widetilde{G}_n^*\nabla U_{n,k,i}(s,X_{1,n}(s,x))-\widetilde{G}_n^*\nabla U_{n,k,i}(s,X_{2,n}(s,x))\|_U^2\right]ds \\
\leq & dM_{T}^2 \|\mathcal{V}^*\|_{\mathcal L(U;U)}^2 \int_0^t\sum_{k=1}^n\sum_{i=1}^{d_k}e^{2(t-s){\rm Re}(\rho^k_{i})}\|B^k_i\|_{C_b^\theta(H_n)}^2\mathbb{E}\left[\|X_{1,n}(s,x)-X_{2,n}(s,x)\|_H^2\right]ds.
\end{align*}
Integrating with respect to $t$ between $0$ and $T$, by the Fubini-Tonelli theorem we get
\begin{align}
I_5
\leq & dM_{T}^2\|\mathcal{V}^*\|_{\mathcal L(U;U)}^2\int_0^T\int_0^t\sum_{k=1}^n\sum_{i=1}^{d_k}e^{2(t-s){\rm Re}(\rho^k_{i})}\|B^k_i\|_{C_b^\theta(H_n)}^2\mathbb{E}\left[\|X_{1,n}(s,x)-X_{2,n}(s,x)\|_H^2\right]dsdt\notag \\
\leq & dM_{T}^2\|\mathcal{V}^*\|_{\mathcal L(U;U)}^2Z\Delta_n,
\label{SS1}
\end{align}
where
\[
Z=\int_0^T\sum_{k=1}^n\sum_{i=1}^{d_k}e^{2t{\rm Re}(\rho^k_{i})}\|B^k_i\|_{C_b^\theta(H_n)}^2dt=-\sum_{k=1}^n\sum_{i=1}^{d_k}\frac{\|B^k_i\|_{C_b^\theta(H_n)}^2}{2{\rm Re}(\rho_i^k)}(1-e^{2T{\rm Re}(\rho_i^k)}).
\]
By \eqref{conv_serie_holder} we infer that there exists a positive constant $\overline C$, independent of $n\in\N$, such that
\begin{align}\label{SSS2}
\int_0^T\sum_{k=1}^n\sum_{j,i=1}^{d_k}e^{2t{\rm Re}(\rho^k_{i})}\|B^k_j\|_{C_b^\theta(H_n)}^2dt\leq \overline C,
\end{align}
so, combing \eqref{SS1} and \eqref{SSS2} we conclude that
\begin{align}\label{SI8-addo}
I_5
\leq dM_{T}^2 \|\mathcal{V}^*\|_{\mathcal L(U;U)}^2\overline C\Delta_n.   \end{align}
By \eqref{SI1}, \eqref{SI2}, \eqref{SI36} and one between \eqref{SI8-traccia} and \eqref{SI8-addo}, we infer that there exists a positive constant $K$, independent of $n\in\N$, such that \begin{equation}\label{SS3}
\Delta_n\leq K\left[M_{T}^2\left(2\Delta_n+\Pi_{n}\right)+\Pi_n\right].
\end{equation}
From \eqref{Convergenza-mild}, the sequence $(\Delta_n)$ converges to $\Delta$ as $n$ goes to $\infty$. Moreover, since $\theta<1$, from H\"older's inequality and \eqref{Convergenza-mild} we deduce that $\Pi_n\to0$ as $n\to\infty$. Therefore, letting $n\to \infty$ in \eqref{SS3} we get
\begin{equation*}
\int_0^T\mathbb{E}\left[\norm{X_1(t,x)-X_2(t,x)}^2_H\right]dt\leq 2KM_{T}^2 \int_0^T\mathbb{E}\left[\norm{X_1(t,x)-X_2(t,x)}^2_H\right]dt.
\end{equation*}
Finally, noticing that by Proposition \ref{SolBKn} $M_T\rightarrow 0$ as $T\rightarrow 0$, if $T>0$ is small enough, then we conclude that
\begin{equation*}
\int_0^{T}\mathbb{E}\left[\norm{X_1(t,x)-X_2(t,x)}^2_H\right]dt=0.
\end{equation*}
The statement for general $T>0$ follows by standard arguments.
\end{proof}

Finally, we can prove corollary \ref{Strong}.

{\begin{proof}[Proof of Corollary \ref{Strong}]
By \cite[Proposition 3]{Cho-Gol1995}, it follows that there exists a weak mild solution to \eqref{eqFO}. Further, Theorem \eqref{pathwiseuniqueness} gives pathwise uniqueness for \eqref{eqFO}. Therefore, from \cite{On04}, which states that weak existence and pathwise uniqueness for equation \eqref{eqFO} imply strong existence, we obtain the desired result.
\end{proof}}

\subsection{An application of Theorem \ref{Ito-Tanaka} to an approximation result}
Let $T>0$. Assume that the assumptions of Corollary \ref{Strong} hold true. Hence, for every $x\in H$, the SPDE \eqref{eqFO} has a unique strong mild solution $\{X(t,x)\}_{t\in [0,T]}$. For every $n\in\N$ and $x\in H$, we are concerned with the unique mild solution $\{\widehat{X}_n(t,x)\}_{t\in [0,T]}$ to 
\begin{align*}
\left\{
\begin{array}{ll}
\displaystyle  d\widehat{X}_n(t)=A_n\widehat{X}_n(t)dt+B_n(\widehat{X}_n(t))dt+G_ndW(t), \quad t\in[0,T],  \vspace{1mm} \\
\widehat{X_n}(0)= P_nx\,
\end{array}
\right.
\end{align*}
where $A_n$, $G_n$ and $B_n$ are given by \eqref{coefficienti-n}. If $B$ is a Lipschitz continuous function, then it is easy to prove that for every fixed $T>0$ and $x\in H$ we have
\begin{equation}\label{Convergenza-mildDF}
\lim_{n \to +\infty}\int^T_0\mathbb{E}\left[\|\widehat{X}_{n}(t,x)-X(t,x)\|_H^2\right]dt=0.
\end{equation}  
When $B$ is only H\"older continuous, it is not trivial to prove that \eqref{Convergenza-mildDF} is verified. In this subsection, exploiting the computations in the proof of Theorem \ref{Ito-Tanaka}, we show that \eqref{Convergenza-mildDF} holds true even in the case where $B$ is only $\theta$-H\"older continuous.

\begin{proposition}
\label{prop:conv_appr_buone}
Assume that Hypotheses \ref{hyp:finito-dimensionale} and one between \ref{hyp:traccia-finita} and \ref{hyp:goal-addo} hold true. Then, for every $T>0$ and $x\in H$, \eqref{Convergenza-mildDF} holds true.
\end{proposition}
\begin{proof}
Let $n\in\N$, let $\{g_1,\ldots,g_{s_n}\}$ be the orthonormal basis of $H_n$ introduced in Remark \ref{remarkhyp1}(iv) (with $s_n={\rm Dim}(H_n)$) and let $U_{n}$ be the solution to \eqref{Back-Kolmon} with $F_n=B_n$.  Arguing as in the proof of Theorem \ref{Ito-Tanaka}, we deduce that 
\begin{align}
\widehat{X}_n(t,x)&=e^{tA_n}( P_nx+U_n(0, P_nx))-U_n(t,\widehat{X}_n(t,x))-A_n\int^t_0e^{(t-s)A_n}U_{n}(s,\widehat{X}_n(s,x))ds\notag\\
&+\int^t_0e^{(t-s)A_n}DU_n(s,\widehat{X}_n(s,x))G_ndW(s)+ \int^t_0e^{(t-s)A_n}G_ndW(s),\quad t\in [0,T],\; \mathbb{P}-{\rm a.s.}\label{I-Tcappuccio}
\end{align}
Let $\{X_n(t,x)\}_{t\in [0,T]}$ be the process defined in \eqref{approsimazione-mild}. By the triangular inequality, for every $n\in\N$ we get 
\begin{align}
\int^T_0\mathbb{E}\left[\|\widehat{X}_{n}(t,x)-X(t,x)\|_H^2\right]dt
\leq &  \int^T_0\mathbb{E}\left[\|X(t,x)-X_n(t,x)\|_H^2\right]dt\notag \\
& +\int^T_0\mathbb{E}\left[\|\widehat{X}_{n}(t,x)-X_n(t,x)\|_H^2\right]dt\label{triangoliamo}.
\end{align} 
From  \eqref{Convergenza-mild}, the first addend in the right-hand side of \eqref{triangoliamo}  goes to $0$. Further, by \eqref{I-T} and \eqref{I-Tcappuccio} we obtain 
\begin{align*}
\int^T_0\mathbb{E}&\left[\|\widehat{X}_{n}(t,x)-X_n(t,x)\|_H^2\right]dt \leq \int_0^T\mathbb{E}\left[\|U_n(t,X_n(t,x))-U_n(t,\widehat{X}_n(t,x))\|^2_H\right]dt\\
&+\int_0^T\mathbb{E}\left[\norm{A_n\int^t_0e^{(t-s)A_n}\left(U_{n}(s,X_{n}(s,x))-U_{n}(s,\widehat{X}_n(s,x))\right)ds}_H^2\right]dt\\
&+\int_0^T\mathbb{E}\left[\norm{\int^t_0e^{(t-s)A_n} \left(B_n(X(s,x))-B_n(X_{n}(s,x))\right)ds}^2_H\right]dt \\
&+\int_0^T\mathbb{E}\left[\norm{\int^t_0e^{(t-s)A_n} DU_{n}(s,X_{n}(s,x))\left(B_n(X(s,x))-B_n(X_{n}(s,x))\right)ds}^2_H\right]dt \\
&\int_0^T\mathbb{E}\left[\norm{\int_0^te^{(t-s)A_n}(DU_n(s,X_{n}(s,x))-DU_n(s,\widehat{X}_n(s,x)))G_ndW(s)}^2_H\right]dt,
\end{align*} 
Arguing as in the proof of Theorem \ref{pathwiseuniqueness}, we infer that there exists a positive constant $K$, independent of $n\in\N$, such that
\begin{align*}
\widehat\Delta_n\leq K\left[M_{T}^2\left(2\widehat\Delta_n+\Pi_{n}\right)+\Pi_n\right], 
\qquad n\in\N,
\end{align*}
where 
\begin{align*}
\widehat \Delta_n:=\int_0^T\mathbb{E}\left[\norm{\widehat X_{n}(t,x)-X_{n}(t,x)}_H^2\right]dt, \quad    
\Pi_{n}:=\int_0^T\mathbb{E}\left[\norm{X_{n}(t,x)-X(t,x)}_H^{2\theta}\right]dt
\qquad n\in\N. 
\end{align*}
Recalling that $M_T\to0$ as $T\to0$ and $\Pi_n\to0$ as $n\to\infty$,  choosing $T$ small enough we deduce that $\widehat\Delta_n\to 0$ as $n\to\infty$. The statement for general $T>0$ follows from standard arguments.
\end{proof}

\section{Applications to specific models}
\label{sec:appl}

\subsection{Stochastic damped wave and Euler-Bernoulli beam equations}\label{Damped}

We consider the following semilinear stochastic differential equation:
\begin{align}
\label{eq_damped_ex}
\left\{
\begin{array}{ll}
\displaystyle \frac{\partial^2y}{\partial t^2}(t)
= -\Lambda y(t)-\rho\Lambda^{\alpha}\left(\frac{\partial y}{\partial t}(t)\right)+C\left(y(t),\frac{\partial y}{\partial t}(t)\right)+\Lambda^{-\gamma}\dot{W}(t), & t\in[0,T], \vspace{1mm} \\
y(0)=y_0\in U, & \vspace{1mm} \\
\displaystyle \frac{\partial y}{\partial t}(0)=y_1\in U,
\end{array}
\right.
\end{align}
where $y$ is a $U$-valued function with $U$ real separable Hilbert space, $\rho,\gamma$ are positive constants and $\alpha\in[0,1)$. Here, $\Lambda:D(\Lambda)\subseteq U\to U$ is a positive self-adjoint operator and $W$ is $U$-valued cylindrical Wiener process. Equations of the form of \eqref{eq_damped_ex} describe elastic systems with structural damping.

Equation \eqref{eq_damped_ex} can be rewritten as \eqref{eqFO} by considering the real separable Hilbert space $H:=U\times U$ and the operators $A:D(A)\subseteq H\to H$ and $\widetilde G:U\to H$ are defined as
\begin{align}
\label{damped_def_op_A_G}
D(A):= & \left\{\begin{pmatrix}
h_1 \\ h_2    
\end{pmatrix}:h_2\in D(\Lambda^{\frac12}), \ h_1+\rho\Lambda^{\alpha-\frac12}h_2\in D(\Lambda^{\frac12})\right\}, \ A=\mathcal A_{\alpha,\rho}:=\begin{pmatrix}
0 & \Lambda^{\frac12} \\
-\Lambda^{\frac12} & -\rho \Lambda^{\alpha}
\end{pmatrix}, \notag \\
\widetilde G:= & \begin{pmatrix}
0 \\ {\rm Id}    
\end{pmatrix}, \quad \mathcal V:=\Lambda^{-\gamma}, \quad G:= \widetilde G\Lambda^{-\gamma}.
\end{align}
For every $t\in[0,T]$ and $h:=\begin{pmatrix}
    h_1 \\ h_2
\end{pmatrix}\in H$, we set
\begin{align*}
X(t):=\begin{pmatrix}
\Lambda^{\frac12}y(t) \\ \frac{\partial y}{\partial t}(t)    
\end{pmatrix}, \qquad \widetilde B(h):=C(\Lambda^{-\frac12}h_1,h_2), \qquad   B(h):=\widetilde G \widetilde B(h) \in H.
\end{align*} 
It follows that equation \eqref{eq_damped_ex} reads as \eqref{eqFO}. Moreover, if $C\in C_b^\theta(H;U)$, then for every $h:=\begin{pmatrix}
h_1 \\ h_2
\end{pmatrix}, k:=\begin{pmatrix}
k_1 \\ k_2
\end{pmatrix}\in H$ we get
\begin{align}
\label{hold_tilde_B}
\|\widetilde B(h)-\widetilde B(k)\|_U^2
\leq & \|C\|_{C^\theta_b(H;U)}^2|(\|\Lambda^{-\frac12}(h_1-k_1)\|_U^2+\|h_2-k_2\|_U^2)^{\theta} \notag \\
\leq & \|C\|_{C^\theta_b(H;U)}^2|(\|\Lambda^{-\frac12}\|_{\mathcal L(U)}^{2}+1)^{\theta}\|h-k\|_H^{2\theta},
\end{align}
which shows that $\widetilde B\in C^\theta_b(H;U)$.

\begin{rmk}
Another possible approach to equations as \eqref{eq_damped_ex}, both in the deterministic and stochastic setting and also for problems without the damping term, is to consider $H:=V\times U$ ($H:U\times V'$ in the stochastic case) and the operator $\widetilde A:D(\widetilde A)\subseteq H\to H$ defined as
\begin{align*}
\widetilde A:=\begin{pmatrix}
 0 & {\rm Id}  \\
 -\Lambda & -\rho\Lambda^\alpha
\end{pmatrix},    
\end{align*}
where $V$ is a suitable real separable Hilbert (see \cite{MasPri17,MasPri23}) and $V'$ is its dual, not identified with it. However, a good choice of $V$ makes these two different approaches equivalent, as shown in \cite{lasi-trig}.
\end{rmk}

We take advantage from a spectral decomposition introduced in \cite{chen-russ} and exploited in \cite{lasi-trig,trig} to construct the space $H_n$, $n\in\N$, which we have used in the approximating procedure. If necessary, we complexify both the spaces and the operators which we deal with. We stress that the operator $A$ generates a strongly continuous semigroup $\{e^{tA}\}_{t\geq0}$ which is analytic if $\alpha\in\left[\frac12,1\right)$.

Let $(e_n)_{n\in\N}$ be a sequence of ({\it non normalized}) eigenvectors of $\Lambda$ with corresponding simple eigenvalues $(\mu_n)_{n\in\N}$ such that $\{e_n:n\in\N\}$ is a basis of $U$ and $(\mu_n)_{n\in\N}$ increases to $\infty$ an $n$ goes to $\infty$. For every $n\in\N$, the values $\lambda^+_n,\lambda^-_n$ defined by
\begin{align}
\label{AvlA}
\lambda^{\pm}_n:=\frac{-\rho\mu_n^\alpha\pm\sqrt{\rho^2\mu^{2\alpha}_n-4\mu_n}}{2}, \qquad \lambda^+_n+\lambda_n^-=-\rho\mu_n^\alpha, \quad \lambda_n^+\lambda_n^-= \mu_n,  
\end{align}
are the eigenvalues of the operator $A=\mathcal A_{\alpha,\rho}$ with corresponding normalized (in $H$) eigenvectors
\begin{align*}
\Phi_n^+=\begin{pmatrix}
\mu_n^{\frac12} e_n \\ \lambda^+_n e_n    
\end{pmatrix},
\quad \Phi_n^-=\chi_n\begin{pmatrix}
\mu_n^{\frac12} e_n \\ \lambda_n^- e_n    
\end{pmatrix}, \qquad n\in\N.
\end{align*}
Further, the adjoint operator $A^*$ of $A$ admits the representation
\begin{align*}
A^*=\begin{pmatrix}
0 & -\Lambda^{\frac12} \\
\Lambda^{\frac12} & -\rho\Lambda^\alpha
\end{pmatrix}
\end{align*}
and has eigenvalues $\lambda_n^+,\lambda_n^-$, with corresponding normalized eigenvectors
\begin{align*}
\Psi_n^+=\begin{pmatrix}
-\mu^{\frac12}_n e_n \\ \lambda^+_n e_n    \end{pmatrix},
\quad \Psi_n^-=\chi_n\begin{pmatrix}
-\mu^{\frac12}_n e_n \\ \lambda_n^- e_n   
\end{pmatrix}, \qquad n\in\N.
\end{align*}
Assume that
\begin{align}
\label{cond_rho_mu_n_alpha}
\rho^2\neq4\mu_n^{1-2\alpha},\quad \forall\; n\in\N.
\end{align}
Under this condition, the eigenvalues $\lambda_n^+,\lambda_n^-$ of $A$ are simple for every $n\in\N$. Further, condition $\|\Phi_n^+\|_H=\|\Phi_n^-\|_H=1$ for every $n\in\N$ imply
\begin{align*}
\|e_n\|_{U}^2(\mu_n^2+|\lambda_n^{+}|^2)=1, \qquad \chi_n^2\|e_n\|^2_U(\mu_n^2+|\lambda_n^-|^2)=1,   \qquad n\in\N, 
\end{align*}
which give
\begin{align*}
\chi_n^2=\frac{\mu_n+|\lambda_n^+|^2}{\mu_n+|\lambda_n^-|^2}, \qquad n\in\N.    
\end{align*}
Each system $\{\Phi_n^+:n\in\N\}$ and $\{\Phi_n^-:n\in\N\}$ is orthonormal in $H$. Further, \eqref{cond_rho_mu_n_alpha} implies that $\{\Phi_n^+,\Phi_n^-:n\in\N\}$ is a (non-orthogonal) basis of $H$ and $H=H^+ + H^-$ (non-orthogonal, direct sum), where
\begin{align*}
H^+:=\overline{{\rm span}\{\Phi_n^+:n\in\N\}}, \qquad H^-:=\overline{{\rm span}\{\Phi_n^-:n\in\N\}}.   
\end{align*}
Each element $h\in H$ can be uniquely decomposed as $h=h^++h^-$ with $h^+\in H^+$ and $h^-\in H^-$. \\
Let us consider the decomposition of $A$, $e^{tA}$, $R(\lambda,A)$ and $\widetilde G$ on $H^+$ and $H^-$. We get
\begin{align*}
Ah = & \sum_{n=1}^\infty \lambda_n^+\langle h^+,\Phi_n^+\rangle_H\Phi_n^++\sum_{n=1}^\infty \lambda_n^-\langle h^-,\Phi_n^-\rangle_H\Phi_n^-, \\
e^{tA}h
= & \sum_{n=1}^\infty e^{\lambda_n^+t}\langle h^+,\Phi_n^+\rangle_H\Phi_n^++\sum_{n=1}^\infty e^{\lambda_n^-t}\langle h^-,\Phi_n^-\rangle_H\Phi_n^-, \qquad t\geq0, \\
R(\lambda,A)h
= & \sum_{n=1}^\infty \frac{1}{\lambda-\lambda_n^+}\langle h^+,\Phi_n^+\rangle_H\Phi_n^++\sum_{n=1}^\infty \frac{1}{\lambda-\lambda_n^-}\langle h^-,\Phi_n^-\rangle_H\Phi_n^-, \qquad \lambda\in\rho(A)
\end{align*}
for every $h\in H$ with $h=h^++h^-$.
From the definition of $\widetilde G$ (see \eqref{damped_def_op_A_G}), it follows that
\begin{align*}
\widetilde Gu=\begin{pmatrix}
0 \\ u   
\end{pmatrix}    
= \sum_{n=1}^\infty (b_n^+u_n\Phi_n^++b_n^-u_n\Phi_n^-), \qquad u\in U,
\end{align*}
where $u_n=\langle u,e_n/\|e_n\|_U\rangle_U$, 
$b_n^++\chi_nb_n^-=0$ and
$(\lambda_n^+b_n^++\chi_n\lambda_n^-b_n^-)\|e_n\|_U^2=1$ for every $n\in\N$. This implies that
\begin{align*}
b_n^+=\frac{1}{\|e_n\|_U(\lambda_n^+-\lambda_n^-)}, \qquad b_n^-=-\frac{b_n^+}{\chi_n}
= \frac{1}{\chi_n\|e_n\|_U(\lambda_n^--\lambda_n^+)}
\end{align*}
for every $n\in\N$. In particular, the explicit expression of $A$ and of $\{e^{tA}\}  _{t\geq0}$ implies that the semigroup $\{e^{tA}\}_{t\geq0}$ is immediately differentiable also for $\alpha\in\left(0,\frac12\right)$, even if not analytic for $\alpha$ belonging to this interval, that for every $\theta>0$ we can define the positive powers $(-A)^\theta$ of $-A$ as
\begin{align*}
(-A)^\theta h:= & \sum_{n=1}^\infty (-\lambda_n^+)^\theta\langle h^+,\Phi_n^+\rangle_H\Phi_n^++\sum_{n=1}^\infty (-\lambda_n^-)^\theta\langle h^-,\Phi_n^-\rangle_H\Phi_n^-
\end{align*}
for every $h=h^++h^-\in H$ such that the above series converge in $H$, that $e^{tA}h\in D((-A)^\theta)$ for every $h\in H$, every $t>0$ and every $\theta>0$, and that, if $\alpha\in\left[\frac12,1\right)$, then for every $T>0$ there exists a positive constant $C$, depending on $\theta$ and $T$, such that
\begin{align}
\label{damped_stima_sing_A}
\|(-A)^\theta e^{tA}\|_{\mathcal L(H)}\leq \frac{C}{t^{\theta}} \qquad \forall t\in(0,T].
\end{align}

From the above construction, it follows that, if $\alpha\in\left(0,\frac12\right]$, then we have
\begin{align}
\notag 
& |\lambda_n^{+}|,|\lambda_n^-|\sim \mu_n^{\frac12}, \quad  \,  \|e_n\|_{U}\sim \mu_n^{-\frac12}, \quad {\rm Re}(\lambda_n^+),{\rm Re}(\lambda_n^-)\sim -\mu_n^\alpha, \quad |\lambda_n^+-\lambda_n^-|\sim \mu_n^{\frac12}, \\
& b_n^{+},b_n^-\sim {\rm const}(b), \quad \chi_n\sim {\rm const}(\chi),
\label{1stime_coefficienti_avl}
\end{align}
definitely with respect to $n\in\N$ (see also formulae \cite[(2.3.14)-(2.3.18)]{trig}). The case $\alpha\in\left[\frac12,1\right)$ is analogously treated. We simply remark that in this case the asymptotic behaviour in \eqref{1stime_coefficienti_avl} is replaced by
\begin{align}
\notag 
& |\lambda_n^+|\sim \mu_n^{1-\alpha}, \quad |\lambda_n^-|\sim\mu_n^\alpha, \quad  \|e_n\|_{U}\sim \mu_n^{-\frac12}, \quad {\rm Re}(\lambda_n^+)\sim -\mu_n^{1-\alpha}, \quad {\rm Re}(\lambda_n^-)\sim -\mu_n^{\alpha}\\
& |\lambda^-_n-\lambda_n^+|\sim \mu_n^{\alpha}, \quad b_n^-\sim {\rm const}(b), \quad \chi_n,b_n^+\sim \mu_n^{\frac12-\alpha},
\label{1stime_coeff_avl_2}
\end{align}
definitely with respect to $n\in\N$. 

Since $\{e^{tA}\}_{t\geq0}$ is an analytic semigroup when $\alpha\in\left[\frac12,1\right)$ (but it does not for $\alpha\in\left[0,\frac12\right)$), Hypothesis \ref{hyp:finito-dimensionale}(i) is fulfilled.
Let us show that, if $\alpha\in\left(0,\frac12\right)$, then estimate \eqref{stima_ris_per_fourier_1} fails. Indeed, for every $k\in\N$, we consider $AR(z_k,A)\Phi_k^+$, where $z_k=a+i{\rm Im}\lambda_k^+$ for some $a\in\R\setminus\{{\rm Re}\lambda_k^+\}$. This implies that
\begin{align*}
\|AR(z_k,A)\Phi_k^+\|_H=\left|\frac{\lambda_k^+}{z_k-\lambda_k^+}\right|
= \left|\frac{\lambda_k^+}{a-{\rm Re}{\lambda_k^+}}\right|\sim\frac{\mu_k^{\frac12}}{\mu_k^{\alpha}}=\mu_k^{\frac12-\alpha},
\end{align*}
and the sequence $(\mu_k^{\frac12-\alpha})_{k\in\N}$ blows up as $k$ goes to $\infty$ since $(\mu_k)_{k\in\N}$ tends to $\infty$ as $k$ goes to $\infty$.

\begin{rmk}
\label{rmk:damped_hyp_eigenvector-A*}
We stress that Hypotheses \ref{hyp:finito-dimensionale}(iii) and \ref{hyp:goal-addo}(a)-(b) are fulfilled. The choice 
\[
H_n:={\rm span}\{\Phi_k^+,\Phi_k^-\; :\; k\in\{1,\ldots n\}\},\qquad n\in\N,
\]
gives Hypotheses \ref{hyp:finito-dimensionale}(iii). Recalling that $A^*=\begin{pmatrix}
0 & -\Lambda^{\frac12} \\  \Lambda^{\frac12} & -\rho\Lambda^\alpha   
\end{pmatrix}$ and noticing that
\[
{\rm span}\{\Psi_k^+,\Psi_k^-\; :\; k\in\{1,\ldots n\}\}={\rm span}\{\Phi_k^+,\Phi_k^-\; :\; k\in\{1,\ldots n\}\}, \qquad n\in\N,
\]
Hypotheses \ref{hyp:goal-addo}(a)-(b) are fulfilled with $d_n=d=2$, $e_1^n=\Psi_n^+$ and $e_2^n=\Psi_n^-$ for every $n\in\N$.
\end{rmk}

\subsubsection{The stochastic convolution}
Now we prove that the stochastic convolution
\begin{align*}
W_A(t):=\int_0^t e^{(t-s)A}GdW(s)    
\end{align*}
is well-defined for every $t\in[0,T]$, where $A$ and $G$ have been introduced in \eqref{damped_def_op_A_G}. In particular, we show that Hypothesis \ref{hyp:finito-dimensionale}(i) is verified.
\begin{proposition}
\label{prop:damped_conv_stoc_1}
Let $A$ and $G$ be as in \eqref{damped_def_op_A_G}. Assume that one of the following conditions holds true:
\begin{enumerate}[\rm(i)]
    \item $\Lambda^{-2\gamma}$ is a trace-class operator on $U$;
    \item $\alpha>0$ and there exist $\delta>0$ and a positive constant $c$ such that for every $n\in\N$ we have $\mu_n\leq c n^\delta$ and $\delta>\frac{1}{2\gamma+\alpha}$.
\end{enumerate}
Therefore, there exists $\eta\in(0,1)$ such that
\begin{align*}
\int_0^T t^{-\eta}{\rm Trace}_H\left[e^{tA}GG^*e^{tA^*}\right]dt<\infty.   
\end{align*}
\end{proposition}
\begin{proof}
Assume that condition $(i)$ holds true. Then,
\begin{align*}
\|e^{tA}G\|_{\mathcal L_2(U;H)}^2
\leq \|e^{tA}\|^2_{\mathcal L(H)}\|\widetilde G\|^2_{\mathcal L(U;H)}\|\Lambda^{-\gamma}\|_{\mathcal L_2(U;U)}<\infty,
\end{align*}
uniformly with respect to $t\in[0,T]$. Since
${\rm Trace}_H[e^{tA}GG^*e^{tA^*}]=\|e^{tA}G\|_{\mathcal L_2(U;H)}^2$ for every $t\in[0,T]$, it follows that $\sup_{t\in[0,T]}{\rm Trace}_H[e^{tA}GG^*e^{tA^*}]<\infty$ and the thesis follows for any choice of $\eta\in(0,1)$.

\vspace{2mm}
We prove the thesis when (ii) holds true. Suppose that there exist $\delta>0$ and a positive constant $c$ such that for every $n\in\N$ we have $\mu_n\leq cn^\delta$, and $\delta>\frac{1}{2\gamma+\alpha}$. The decomposition of $H$ gives
\begin{align}
\|e^{tA}G\|_{\mathcal L_2(U;H)}^2
= & \sum_{n=1}^\infty \|e^{tA}G(e_n/\|e_n\|_U)\|_H^2 \notag \\
= & \sum_{n=1}^\infty \|e^{\lambda_n^+t}\mu_n^{-\gamma}b_n^+ \Phi_n^++e^{\lambda_n^-t}\mu_n^{-\gamma}b_n^-\Phi_n^-\|_H^2 \notag  \\
\leq & 2 \sum_{n=1}^\infty \mu_n^{-2\gamma}(\|e^{\lambda_n^+t}b_n^+ \Phi_n^+\|_H^2+\|e^{\lambda_n^-t}b_n^-\Phi_n^-\|_H^2).
\label{stima_HS_norm_1}
\end{align}
for every $t\in(0,\infty)$. We separately consider the cases $\alpha\in\left(0,\frac12\right)$ and $\alpha\in\left[\frac12,1\right)$. 

$\bullet \ \alpha\in\left(0,\frac12\right)$. From \eqref{1stime_coefficienti_avl} and \eqref{stima_HS_norm_1}, we get
\begin{align*}
\|e^{tA}G\|_{\mathcal L_2(U;H)}^2
\leq 
& 2|{\rm const}(b)|\sum_{n=1}^\infty \mu_n^{-2\gamma}(e^{2{\rm Re}(\lambda_n^+)t}+e^{2{\rm Re}(\lambda_n^-)t}) \\
\sim & 2|{\rm const}(b)|\sum_{n=1}^\infty \mu_n^{-2\gamma}e^{-2\rho\mu_n^\alpha t}
\sim 2|{\rm const}(b)|\sum_{n=1}^\infty n^{-2\delta \gamma}e^{-2\rho n^{\delta\alpha} t}
\end{align*}
for every $t\in(0,\infty)$. Hence, there exists a positive constant $\hat c$, which may vary line to line, such that
\begin{align*}
\int_0^T t^{-\eta}\|e^{tA}G\|_{\mathcal L_2(U;H)}^2dt
\leq &  \hat c\int_0^T t^{-\eta}\sum_{n=1}^\infty n^{-2\delta\gamma}e^{-2\rho n^{\delta\alpha} t}dt \\
\sim & \hat c\int_0^Tt^{-\eta}\int_1^\infty x^{-2\delta\gamma} e^{-2\rho x^{\delta\alpha}t}dxdt.
\end{align*}
Let us consider the change of variables $y=x^{\delta\alpha}t$. We get
\begin{align}
\int_0^T t^{-\eta}\|e^{tA}G\|_{\mathcal L_2(U;H)}^2dt
\leq & \hat c \int_0^Tt^{\frac{2\gamma}\alpha-\frac1{\delta\alpha}-\eta}\int_t^\infty y^{-\frac{2\gamma}\alpha+\frac1{\delta\alpha}-1} e^{-2\rho y}dy dt\notag    \\
\leq & \hat c\int_0^Tt^{\frac{2\gamma}\alpha-\frac1{\delta\alpha}-\eta}\int_t^T y^{-\frac{2\gamma}\alpha+\frac1{\delta\alpha}-1} dy dt  \notag  \\
& +\hat c\int_0^Tt^{\frac{2\gamma}\alpha-\frac1{\delta\alpha}-\eta}\int_T^\infty y^{-\frac{2\gamma}\alpha+\frac1{\delta\alpha}-1} e^{-2\rho y}dy dt   \notag \\
\leq & \hat c \left(\int_0^Tt^{-\eta}dt+\int_0^Tt^{\frac{2\gamma}\alpha-\frac{1}{\delta\alpha}-\eta}dt\right).\label{conto-damped}
\end{align}
Let us notice that 
\begin{align*}
\frac{2\gamma}{\alpha}-\frac{1}{\delta\alpha}>-1 \Longleftrightarrow
\frac{1}{\alpha}\left(2\gamma-\frac{1}{\delta}\right)>-1
\Longleftrightarrow \delta>\frac{1}{2\gamma+\alpha},
\end{align*}
which means that, under our assumptions, choosing $\eta\in\left(0,\left(1+\frac{2\gamma}\alpha-\frac{1}{\delta\alpha}\right)\wedge 1\right)$, we get 
\begin{align*}
\int_0^T t^{-\eta}\|e^{tA}G\|_{\mathcal L_2(U;H)}^2dt<\infty.
\end{align*}

$\bullet$ $\alpha\in\left[\frac12,1\right)$. From \eqref{1stime_coeff_avl_2} and  \eqref{stima_HS_norm_1}, we infer that
\begin{align*}
\|e^{tA}G\|_{\mathcal L_2(U;H)}^2
\leq 
& 2\sum_{n=1}^\infty \mu_n^{-2\gamma}\left[(\mu_n^{\frac12-\alpha}e^{{\rm Re}(\lambda_n^+)t})^2+|{\rm const}(b)^2|e^{2{\rm Re}(\lambda_n^-)t}\right] \\
& \sim 2\sum_{n=1}^\infty \mu_n^{-2\gamma}\left[\mu_n^{1-2\alpha}e^{-2\rho\mu_n^{1-\alpha}t}+|{\rm const}(b)|^2e^{-2\rho \mu_n^\alpha t}\right] \\
& \sim\sum_{n=1}^\infty n^{\delta(1-2\gamma-2\alpha)}e^{-2\rho n^{\delta(1-\alpha)}t}
\end{align*}
for every $t\in(0,\infty)$. Hence, there exists a positive constant $\hat c$, which may vary line to line, such that
\begin{align*}
\int_0^T t^{-\eta}\|e^{tA}G\|_{\mathcal L_2(U;H)}^2dt
\leq &  \hat c\int_0^T t^{-\eta}\sum_{n=1}^\infty n^{\delta(1-2\gamma-2\alpha)}e^{-2\rho n^{\delta(1-\alpha)} t}dt \\
\sim & \hat c\int_0^Tt^{-\eta}\int_1^\infty x^{\delta(1-2\gamma-2\alpha)} e^{-2\rho x^{\delta(1-\alpha)}t}dxdt.
\end{align*}
Let us consider the change of variables $y=x^{\delta(1-\alpha)}t$. We get
\begin{align}
\int_0^T t^{-\eta}\|e^{tA}G\|_{\mathcal L_2(U;H)}^2dt
\leq & \hat c \int_0^Tt^{-\frac{1-2\gamma-2\alpha}{1-\alpha}-\frac{1}{\delta(1-\alpha)}-\eta}\int_t^\infty y^{\frac{1-2\gamma-2\alpha}{1-\alpha}+\frac{1}{\delta(1-\alpha)}-1} e^{-2\rho y}dy dt\notag    \\
\leq & \hat c\int_0^Tt^{-\frac{1-2\gamma-2\alpha}{1-\alpha}-\frac{1}{\delta(1-\alpha)}-\eta}\int_t^T y^{\frac{1-2\gamma-2\alpha}{1-\alpha}+\frac{1}{\delta(1-\alpha)}-1}  dy dt  \notag  \\
& +\hat c\int_0^Tt^{-\frac{1-2\gamma-2\alpha}{1-\alpha}-\frac{1}{\delta(1-\alpha)}-\eta}\int_T^\infty y^{\frac{1-2\gamma-2\alpha}{1-\alpha}+\frac{1}{\delta(1-\alpha)}-1}  e^{-2\rho y}dy dt   \notag \\
\leq & \hat c \left(\int_0^Tt^{-\eta}dt+\int_0^Tt^{-\frac{1-2\gamma-2\alpha}{1-\alpha}-\frac{1}{\delta(1-\alpha)}-\eta}dt\right).
\label{conto-damped2}
\end{align}
Let us notice that 
\begin{align*}
-\frac{1-2\gamma-2\alpha}{1-\alpha}-\frac{1}{\delta(1-\alpha)}>-1
\Longleftrightarrow
1-2\gamma-2\alpha+\frac{1}{\delta}<1-\alpha
\Longleftrightarrow \delta>\frac{1}{2\gamma+\alpha},
\end{align*}
which means that, under our assumptions, choosing $\eta\in\left(0,\left(1-\frac{1-2\gamma-2\alpha}{1-\alpha}-\frac{1}{\delta(1-\alpha)}\right)\wedge 1\right)$, we get 
\begin{align*}
\int_0^T t^{-\eta}\|e^{tA}G\|_{\mathcal L_2(U;H)}^2dt<\infty.
\end{align*}
\end{proof}

\subsubsection{The control problem}

We recall that we have set $H_n:={\rm span}\{\Phi_k^+,\Phi_k^-:k=1,\ldots,n\}$ for every $n\in\N$. From the above discussion, Hypotheses \ref{hyp:finito-dimensionale}(v) is satisfied. Further, since $H_n\cap H_{n-1}^\perp={\rm span}\{\Psi_n^+,\Psi_n^-\}$ for every $n\in\N$, also Hypotheses \ref{hyp:goal-addo}(a)-(b) are fulfilled with $d_n=2$ for every $n\in\N$. It remains to prove that Hypotheses \ref{hyp:finito-dimensionale}(vi) are verified.

To this aim, we take advantage from the technique applied in \cite[Proposition 1.3]{Zab08} in finite dimension, and generalized in infinite dimension in \cite{MasPri17,MasPri23}, for the case of wave equation, and in \cite{AddMasPri23} for the case of damped equation. Here, we apply such a method in finite dimension looking for estimates which are independent of the dimension.

For every $t>0$ we consider the control problem
\begin{align}
\label{damped_contr_prob_n}
\left\{
\begin{array}{ll}
Y'(\tau )=AY(\tau )+Gu(\tau ), & \tau \in(0,t], \vspace{1mm} \\
Y(0)=h\in H, & 
\end{array}
\right.
\end{align}
where $u:[0,t]\to U$. Let us notice that for every $v\in U$ we get (recall that $e_k$'s are not normalized) 
\begin{align*}
Gv
= \begin{pmatrix}
 0 \\ \Lambda^{-\gamma}u   
\end{pmatrix}
= \sum_{k=1}^{+\infty}\begin{pmatrix}
0 \\ \widetilde v_ke_k/\|e_k\|_U    
\end{pmatrix}, \qquad n\in\N,
\end{align*}
where $\widetilde v_k:=\langle \Lambda^{-\gamma}v,e_k/\|e_k\|_U\rangle_U$ for every $k\in\N$, and
\begin{align*}
\widetilde G\Lambda^{-\gamma}v
= \begin{pmatrix}
0 \\ \Lambda^{-\gamma}v  
\end{pmatrix}
= \sum_{k=1}^{+\infty}\begin{pmatrix}
 0 \\ \widetilde v_ke_k/\|e_k\|_U   
\end{pmatrix}, \qquad n\in\N.
\end{align*}

We say that problem \eqref{damped_contr_prob_n} is {\it null-controllable} if for every $t>0$ and every $h\in H$ there exists a control $u\in L^2(0,t;U)$ such that \eqref{damped_contr_prob_n} admits a unique mild solution $Y$ and $Y(t)=0$. A mild solution $Y$ to \eqref{damped_contr_prob_n} is a function $Y:[0,t]\to H$ which fulfills
\begin{align}
\label{damped_mild_sol_contr_prob}
y(\tau )=e^{\tau A}h+\int_0^\tau e^{(\tau -s)A}Gu(s)ds, \qquad \tau \in[0,t].    
\end{align}
The null-controllability of \eqref{damped_contr_prob_n} is equivalent to the fact that $e^{tA}(H)\subset Q_{t}^{\frac12}(H)$ (see for instance \cite[Theorem 2.3]{Zab08}). If we denote by $\mathcal E_{C}(t,h)$ the minimal energy to steer $h$ to $0$ at time $t$, i.e.,
\begin{align*}
\mathcal E_{C}(t,h):=\inf\{\|u\|_{L^2(0,t;U)}:\textrm{ \eqref{damped_contr_prob_n} admits a unique solution $Y$ with $Y(t)=0$} \},    
\end{align*}
then  $\|Q_{t}^{-\frac12}e^{tA}h\|_{H}=\mathcal E_{C}(t,h)$, see again \cite[Theorem 2.3]{Zab08}. 

In order to apply the abstract results, we need an estimate of $\|Q_{t}^{-\frac12}e^{tA}h\|_{H_n}$ when $h\in H$ and $h=\widetilde Gu:=P_n\widetilde Gu$ for some $u\in U$. At first, we provide an estimate for $\mathcal E_{C}(t,h)$ with $h\in H$.

\begin{thm}
\label{thm:damped_contr_>1/2}
Let $\alpha\in\left[\frac12,1\right)$ and $\gamma\geq0$. Then, system \eqref{damped_contr_prob_n} is null-controllable. Further, for every $t>0$ there exists a positive constant $\overline c$, which depends on $\alpha$ and $\gamma$ but is independent of $t$, if $t$ varies in a bounded interval, such that
\begin{align}
\label{stima_energia_damped_<1/2}
\mathcal E_{C}(t,h)
\leq 
\begin{cases}
\displaystyle \frac{\overline c\|h\|_{H}}{t^{\frac12+(\gamma+\alpha-\frac12)/(1-\alpha)}}, & \gamma+2\alpha\geq \frac32, \\[4mm] 
\displaystyle \frac{\overline c\|h\|_{H}}{t^{\frac32}}, & \gamma+2\alpha< \frac32, 
\end{cases}
\qquad \forall h\in H.
\end{align}
\end{thm}
\begin{proof}
As already noticed, we adapt the method exploited in \cite[Proposition 1.3]{Zab08}. For reader's convenience, we split the proof into three steps. In the first step we show that system \eqref{damped_contr_prob_n} is null-controllable, providing, for every $h\in H$ and $t>0$, an explicit control $u$, of the form $u=u_1+u_2$, which steers $h$ to $0$ at time $t$. In the second step we compute the $L^2$-estimate of $u_1$, while in the last step we estimate the $L^2$-norm of $u_2$. Combining these estimates, we conclude the proof.

{\bf Step $\bm1$}. 
Let us fix $t>0$. If we consider the matrix representation of the operators $A$ and $G$ on $H$, then the $2\times 2$ matrix 
\begin{align*}
[G|AG]=\begin{pmatrix}
0 & \Lambda^{\frac12-\gamma} \\
\Lambda^{-\gamma} & -\rho\Lambda^{\alpha-\gamma}
\end{pmatrix}
= \Lambda^{-\gamma}
\begin{pmatrix}
0 & \Lambda^\frac12 \\
{\rm Id} & -\rho\Lambda^{\alpha}
\end{pmatrix}
\end{align*}
has the following (formal) inverse,
\begin{align*}
K:=[G|AG]^{-1}=\begin{pmatrix}
\rho\Lambda^{\alpha-\frac12+\gamma} & \Lambda^{\gamma} \\
\Lambda^{-\frac12+\gamma} & 0
\end{pmatrix}=\Lambda^\gamma\begin{pmatrix}
\rho\Lambda^{\alpha-\frac12} & {\rm Id} \\
\Lambda^{\frac12} & 0
\end{pmatrix}.    
\end{align*}
We denote by $K_i$, $i=1,2$, the $i$-th row of $K$, and we consider the control $u:[0,t]\to U$ defined as
\begin{align}
\label{control}
u(\tau)=
\begin{cases}
K_1 \psi_t(\tau)+K_2\psi_t'(\tau), & \forall \tau\in(0,t), \vspace{1mm} \\
0, & \tau=0, \ \tau=t,
\end{cases}
\end{align}
where $\psi_t(\tau)=-\Phi_t(\tau)e^{\tau A_n}h$ for every $\tau \in(0,t)$ and $\Phi_t:[0,t]\to \mathbb R$ is defined as $\Phi_t(\tau )=\overline c_m \tau ^m(t-\tau )$ for every $\tau\in[0,t]$, $\overline c_m$ is a normalizing constant which gives $\|\Phi_t\|_{L^1(0,t)}=1$ and $m\in\N$ satisfies $-2(\gamma+\alpha-\frac12)/(1-\alpha)+2m>-1$. Let us notice that $\psi_t$ is differentiable in $(0,t)$ and $\psi_t'(\tau )=-\Phi_t'(\tau )e^{\tau A}h-\Phi_t(\tau )A_ne^{\tau A}h$. Further, the operators $K_1,K_2:H\to U$ are meant as
\begin{align*}
K_1k
= & K_1\begin{pmatrix}
    k_1 \\ k_2
\end{pmatrix}= \rho \Lambda^{\alpha-\frac12+\gamma} k_1+\Lambda^{\gamma}k_2\in U, \\   
K_2k
= & K_2\begin{pmatrix}
    k_1 \\ k_2
\end{pmatrix}=  \Lambda^{-\frac12+\gamma} k_1\in U
\end{align*}
for every $k=\begin{pmatrix}
 k_1 \\ k_2   
\end{pmatrix}\in H$. It follows that
\begin{align*}
GK_1k+AGK_2k=
\begin{pmatrix}
0 \\ \rho\Lambda^{\alpha-\frac12}k_1+k_2    
\end{pmatrix}
+\begin{pmatrix}
k_1 \\ -\rho\Lambda^{\alpha-\frac12}k_1  
\end{pmatrix}
=k,
\qquad k=\begin{pmatrix}
 k_1 \\ k_2   
\end{pmatrix}\in H.    
\end{align*}
We show that $u$ steers $h$ to $0$ at $t$. Integrating by parts $\psi_t'$, it follows that 
\begin{align*}
\int_0^te^{(t-s)A}Gu(s)ds
= & \int_0^te^{(t-s)A}GK_1\psi_t(s)ds
+ \int_0^te^{(t-s)A}GK_2\psi_t'(s)ds \\
= & \int_0^te^{(t-s)A}GK_1\psi_t(s)ds
+ e^{(t-s)A}GK_2\psi_t(s)\Big|^{t}_0 \\
& + \int_0^te^{(t-s)A}AGK_2\psi_t(s)ds \\
= &  \int_0^te^{(t-s)A}\psi_t(s)ds \\
= & -e^{tA}h\int_0^t\Phi_t(s)ds =-e^{tA}h.
\end{align*}
Replacing this equality in \eqref{damped_mild_sol_contr_prob} we get $Y(t)=0$. 

{\bf Step $\bm 2$}. Here, we estimate the $L^2$-norm of $u_1=K_1\psi_t$. We stress that
\begin{align*}
\|K_1\psi_t(\tau )\|_U^2
=  &  |\Phi_t(\tau )|^2\|\rho\Lambda^{\alpha-\frac12+\gamma}(e^{\tau A}h)_1+\Lambda^{\gamma}(e^{\tau A}h)_2\|^2_U, \qquad \tau\in(0,t).
\end{align*}
From \eqref{AvlA} we get
\begin{align}
& \rho\Lambda^{\alpha-\frac12+\gamma}(e^{\tau A}h)_1+\Lambda^{\gamma}(e^{\tau A}h)_2 \notag \\
= & \sum_{k=1}^{+\infty}[e^{\lambda_k^+\tau }\langle h^+,\Phi_k^+\rangle_H(\rho \mu_k^{\alpha-\frac12+\gamma}\mu_k^{\frac12}+\mu_k^\gamma\lambda_k^+)e_k
+e^{\lambda_k^-\tau }\langle h^-,\Phi_k^- \rangle_H \chi_k(\rho \mu_k^{\alpha-\frac12+\gamma}\mu_k^{\frac12}+\mu_k^\gamma\lambda_k^-)e_k] \notag \\
= & 
\sum_{k=1}^{+\infty}\mu_k^\gamma[e^{\lambda_k^+\tau }\langle h^+,\Phi_k^+\rangle_H(\rho \mu_k^{\alpha}+\lambda_k^+)e_k
+e^{\lambda_k^-\tau }\langle h^-,\Phi_k^- \rangle_H \chi_k(\rho \mu_k^{\alpha}+\lambda_k^-)e_k] \notag \\
= & -\sum_{k=1}^{+\infty}\mu_k^\gamma[e^{\lambda_k^+\tau }\langle h^+,\Phi_k^+\rangle_H\lambda_k^-e_k
+e^{\lambda_k^-\tau }\langle h^-,\Phi_k^- \rangle_H \chi_k\lambda_k^+e_k].
\label{forma_K_1psi_t}
\end{align}
As far as the first addend is concerned, from \eqref{damped_stima_sing_A} and \eqref{1stime_coeff_avl_2} we infer that
\begin{align*}
\bigg\|\sum_{k=1}^{+\infty}\mu_k^\gamma e^{\lambda_k^+\tau }\langle h^+,\Phi_k^+\rangle_H\lambda_k^-e_k\bigg\|_U^2
\sim & \sum_{k=1}^{+\infty} \langle h^+,\Phi_k^+\rangle_H^2|\mu_k^{\gamma+\alpha-\frac12}e^{\lambda_k^+\tau }|^2 \\
\sim & \sum_{k=1}^{+\infty} \langle h^+,\Phi_k^+\rangle_H^2|(\lambda_k^+)^{(\gamma+\alpha-\frac12)/(1-\alpha)}e^{\lambda_k^+\tau }|^2 \\
\leq & \sum_{k=1}^{+\infty} \langle h^+,\Phi_k^+\rangle_H^2\|(-A)^{(\gamma+\alpha-\frac12)/(1-\alpha)}e^{A\tau }\|_{\mathcal L(H)}^2 \\
\leq & C^2\tau ^{-2(\gamma+\alpha-\frac12)/(1-\alpha)}\|h\|_{H}^2
\end{align*}
for every $\tau \in(0,t]$, and $C$ is a positive constant which does not depend on $\tau $ and $n$. As far as the second sum is considered, we get
\begin{align*}
\bigg\|\sum_{k=1}^{+\infty} \mu_k^\gamma e^{\lambda_k^-\tau }\langle h^-,\Phi_k^-\rangle_H\lambda_k^+\chi_k e_k\bigg\|_U^2
\sim & \sum_{k=1}^{+\infty} \langle h^-,\Phi_k^-\rangle_H^2|\mu_k^{\gamma+1-2\alpha}e^{\lambda_k^-\tau }|^2 \\
\sim & \sum_{k=1}^{+\infty} \langle h^-,\Phi_k^-\rangle_H^2|(-\lambda_k^-)^{(\gamma+1-2\alpha)/\alpha}e^{\lambda_k^-\tau }|^2 \\
\leq & \sum_{k=1}^{+\infty} \langle h^-,\Phi_k^-\rangle_H^2\|(-A)^{(\gamma+1-2\alpha)/\alpha}e^{A\tau }\|_{\mathcal L(H)}^2 \\
\leq & C^2\tau ^{(-2(\gamma+1-2\alpha)/\alpha)\wedge 0}\|h\|_{H}^2
\end{align*}
for every $\tau\in(0,t]$, and $C$ is a positive constant which does not depend on $\tau $. We have
\begin{align*}
\frac{\gamma+\alpha-\frac12}{1-\alpha}
\geq \frac{\gamma+1-2\alpha}{\alpha}
\Longleftrightarrow
-\alpha^2+\alpha\left(2\gamma+\frac52\right)-1-\gamma\geq0.
\end{align*}
The solution of the associated homogeneous equation are $\alpha_1=\frac12$ and $\alpha_2=2\gamma+2>1$. Since $\alpha\in\left[\frac12,1\right)$, it follows that
\begin{align}
\int_0^t\|K_1\psi_t(\tau )\|_U^2d\tau 
\leq & C^2\int_0^t\tau^{-2(\gamma+\alpha-\frac12)/(1-\alpha)}|\Phi_t(\tau)|^2d\tau 
\|h\|_{H_n}^2 
\leq \widetilde c_1^2t^{-2(\gamma+\alpha-\frac12)/(1-\alpha)-1}\|h\|_{H}^2,
\label{damped_stima_controllo_<1/2_1}
\end{align}
where we have used the fact that $-2\left(\gamma+\alpha-\frac12\right)/(1-\alpha)+2m>-1$.

{\bf Step $\bm3$}. Let us consider the second addend which defines $u$. We have
\begin{align}
K_2\psi_t'(\tau )
= & -\Lambda^{-\frac12+\gamma}(\Phi'_t(\tau )(e^{\tau A_n}h)_1+\Phi_t(\tau )(Ae^{\tau A}h)_1) \notag \\
= & -\sum_{k=1}^{+\infty}\mu_k^{-\frac12+\gamma}[\Phi'_t(\tau )(e^{\lambda_k^+\tau }\langle h^+,\Phi_k^+\rangle_H+e^{\lambda_k^-\tau }\langle h^-,\Phi_k^-\rangle_H\chi_k)\mu_k^{\frac12}e_k \notag \\
& +\Phi_t(\tau )(\lambda_k^+e^{\lambda_k^+\tau }\langle h^+,\Phi_k^+\rangle_H+\lambda_k^-e^{\lambda_k^-\tau }\langle h^-,\Phi_k^-\rangle_H\chi_k)\mu_k^{\frac12}e_k] \notag \\
= & -\sum_{k=1}^{+\infty}\mu_k^{\gamma}[\Phi'_t(\tau )(e^{\lambda_k^+\tau }\langle h^+,\Phi_k^+\rangle_H+e^{\lambda_k^-\tau }\langle h^-,\Phi_k^-\rangle_H\chi_k)e_k \notag \\
& +\Phi_t(\tau )(\lambda_k^+e^{\lambda_k^+\tau }\langle h^+,\Phi_k^+\rangle_H+\lambda_k^-e^{\lambda_k^-\tau }\langle h^-,\Phi_k^-\rangle_H\chi_k)e_k]
\label{forma_K_2psi_t'}
\end{align}
for every $\tau \in(0,t]$. 
If we separately estimate the two series, from \eqref{1stime_coeff_avl_2} we obtain
\begin{align*}
& \bigg\|\sum_{k=1}^{+\infty}\mu_k^{\gamma}\Phi'_t(\tau )(e^{\lambda_k^+\tau }\langle h^+,\Phi_k^+\rangle_H+e^{\lambda_k^-\tau }\langle h^-,\Phi_k^-\rangle_H\chi_k)e_k \bigg\|_U^2 \\
\sim & |\Phi_t'(\tau )|^2\sum_{k=1}^{+\infty}(\langle h^+,\Phi_k^+\rangle_H^2|(\lambda_k^+)^{\gamma/(1-\alpha)}e^{\lambda_k^+\tau }|^2+\langle h^-,\Phi_k^-\rangle_H^2|(\lambda_k^-)^{(\gamma+\frac12-\alpha)/\alpha} e^{\lambda_k^-\tau }|^2)\|e_k\|_U^2 \\
\sim & |\Phi_t'(\tau )|^2\sum_{k=1}^{+\infty}(\langle h^+,\Phi_k^+\rangle_H^2|(\lambda_k^+)^{(\gamma-\frac12)/(1-\alpha)}e^{\lambda_k^+\tau }|^2+\langle h^-,\Phi_k^-\rangle_H^2|(\lambda_k^-)^{(\gamma-\alpha)/\alpha}e^{\lambda_k^-\tau }|^2) \\
\leq & |\Phi_t'(\tau )|^2\sum_{k=1}^{+\infty}(\langle h^+,\Phi_k^+\rangle_H^2 \|(-A)^{(\gamma-\frac12)/(1-\alpha)}e^{\tau A}\|_{\mathcal L(H)}^2
+\langle h^-,\Phi_k^-\rangle_H^2\|(-A)^{(\gamma-\alpha)/\alpha}e^{\tau A}\|_{\mathcal L(H)}^2) \\
\leq & C^2|\Phi_t'(\tau )|^2\sum_{k=1}^{+\infty}(\tau ^{(-2(\gamma-\frac12)/(1-\alpha))\wedge 0}\langle h^+,\Phi_k^+\rangle_H^2
+\tau ^{(-2(\gamma-\alpha)/\alpha)\wedge 0}\langle h^-,\Phi_k^-\rangle_H^2) \\
\leq & C^2|\Phi_t'(\tau )|^2\tau ^{(-2(\gamma-\frac12)/(1-\alpha))\wedge 0}\|h\|_{H}^2
\end{align*}
for every $\tau \in(0,t]$, since if $0\leq \gamma\leq \frac12\leq \alpha$ then both the exponents of $\tau$ are $0$, if $\gamma\in\left(\frac12,\alpha\right)$ then the first exponent is negative and the second one is $0$, and if $\gamma\geq \alpha$ then
\begin{align*}
\frac{\gamma-\frac12}{1-\alpha}\geq \frac{\gamma-\alpha}{1-\alpha}\geq \frac{\gamma-\alpha}{\alpha}\geq 0.    
\end{align*}
Here, $C$ is a positive constant which does not depend on $\tau $ and $n$. Similar arguments give
\begin{align*}
& \bigg\|\sum_{k=1}^{+\infty}\mu_k^{\gamma}\Phi_t(\tau )(\lambda_k^+e^{\lambda_k^+\tau }\langle h^+,\Phi_k^+\rangle_H+\lambda_k^-e^{\lambda_k^-\tau }\langle h^-,\Phi_k^-\rangle_H\chi_k)e_k 
\bigg\|_U^2 \\
\sim & |\Phi_t(\tau )|^2\sum_{k=1}^{+\infty} (\langle h^+,\Phi_k^+\rangle_H^2|(\lambda_k^+)^{1+\gamma/(1-\alpha)}e^{\lambda_k^+\tau }|^2+\langle h^-,\Phi_k^-\rangle_H^2|(\lambda_k^-)^{1+(\gamma+\frac12-\alpha)/\alpha}e^{\lambda_k^-\tau }|^2)\|e_k\|_U^2 \\
\sim & |\Phi_t(\tau )|^2\sum_{k=1}^{+\infty} (\langle h^+,\Phi_k^+\rangle_H^2|(\lambda_k^+)^{1+(\gamma-\frac12)/(1-\alpha)}e^{\lambda_k^+\tau }|^2+\langle h^-,\Phi_k^-\rangle_H^2|(\lambda_k^-)^{\gamma/\alpha}e^{\lambda_k^-\tau }|^2) \\
\leq & |\Phi_t(\tau )|^2\sum_{k=1}^{+\infty} (\langle h^+,\Phi_k^+\rangle_H^2\|(-A)^{1+(\gamma-\frac12)/(1-\alpha)}e^{\tau A}\|_{\mathcal L(H)}^2+\langle h^-,\Phi_k^-\rangle_H^2\|(-A)^{\gamma/\alpha}e^{\tau A}\|_{\mathcal L(H)}^2) \\
\leq &  C^2|\Phi_t(\tau )|^2\sum_{k=1}^{+\infty} (\langle h^+,\Phi_k^+\rangle_H^2\tau ^{(-2-2(\gamma-\frac12)/(1-\alpha))\wedge 0}+\langle h^-,\Phi_k^-\rangle_H^2\tau ^{-2\gamma/\alpha})
\end{align*}
for every $\tau \in(0,t]$, and $C$ is a positive constant which does not depend on $\tau $ and $n$. We claim that the second addend in $K_2\psi_t'$ is a little $o$ of the first addend of $K_2\psi_t'$ as $\tau$ tends to $0$. Indeed, $|\Phi_t(\tau)|\sim\frac{\tau^{m}}{t^{m+1}}$ and $|\Phi_t'(\tau)|\sim\frac{\tau^{m-1}}{t^{m+1}}$ for $\tau$ near $0$. Hence,
\begin{align*}
|\Phi_t'(\tau)|^2\tau^{(-2(\gamma-\frac12)(1-\alpha))\wedge 0}
\sim \frac{1}{t^{2m+2}}\tau^{2m-2+(-2(\gamma-\frac12)/(1-\alpha))\wedge 0}, \qquad \tau\in(0,t),
\end{align*}
and
\begin{align*}
& |\Phi_t(\tau)|^2
(\tau^{(-2-2(\gamma-\frac12)/(1-\alpha))\wedge 0}+\tau^{-2\gamma/\alpha})
\leq  
\frac{C}{t^{2m+2}}
\begin{cases}
\tau^{2m-2}, & \gamma\in\left[0,\frac12\right), \\[1mm]
\tau^{2m-2-2(\gamma-\frac12)/(1-\alpha)}, & \gamma \geq\frac12,
\end{cases}
\end{align*}
for some positive constant $C$. The claim is so proved. \\
Since $2m-2-2(\gamma-\frac12)/(1-\alpha)
= 2m-2(\gamma+\frac12-\alpha)/(1-\alpha)\geq 2m-2(\gamma+\alpha-\frac12)/(1-\alpha)>-1$, it follows that
\begin{align}
\int_0^t\|K_2\psi_t'(\tau )\|_U^2d\tau 
\leq & \frac{C^2}{t^{2m+2}}\int_0^t \tau ^{2m-2+(-2(\gamma-\frac12)/(1-\alpha))\wedge 0}d\tau \|h\|_{H}^2
\leq C^2t^{-3+(-2(\gamma-\frac12)/(1-\alpha))\wedge0}\|h\|_{H}^2.
\label{damped_stima_controllo_<1/2_2_1}
\end{align}
It remains to compare the exponents of $t$ in \eqref{damped_stima_controllo_<1/2_1} and \eqref{damped_stima_controllo_<1/2_2_1}. We split the cases $\gamma\in\left[0,\frac12\right)$ and $\gamma \geq\frac12$. In the first situation, the exponent in \eqref{damped_stima_controllo_<1/2_1} is $-1-2(\gamma+\alpha-\frac12)/(1-\alpha)$ and that in \eqref{damped_stima_controllo_<1/2_2_1} is $3$. Hence,
\begin{align*}
-1-\frac{2\gamma+2\alpha-1}{1-\alpha}\leq -3
\Longleftrightarrow 2\gamma+2\alpha-1\geq2-2\alpha \Longleftrightarrow \gamma+2\alpha\geq \frac32. 
\end{align*}
If $\gamma\geq \frac12$, then $\alpha\in\left[\frac12,1\right)$ implies that
\begin{align*}
-1-\frac{2\gamma+2\alpha-1}{1-\alpha}
\leq -1-\frac{2\gamma-2\alpha+1}{1-\alpha}=-3-\frac{2\gamma-1}{1-\alpha},
\end{align*}
which means that the exponent of $t$ in \eqref{damped_stima_controllo_<1/2_1} is smaller than that in \eqref{damped_stima_controllo_<1/2_2_1}. \\ Finally, we notice that $\gamma\geq \frac12$ implies $\gamma+2\alpha\geq \frac32$, since $\alpha\in\left[\frac12,1\right)$. Hence, we obtain
\begin{align*}
\mathcal E_{C}(t,h)\leq \frac{\overline c\|h\|_{H}}{t^{\frac12+(\gamma+\alpha-\frac12)/(1-\alpha)}}, \qquad \gamma+2\alpha\geq \frac32,     
\end{align*}
and
\begin{align*}
\mathcal E_{C}(t,h)\leq \frac{\overline c\|h\|_{H}}{t^{\frac32}}, \qquad \gamma+2\alpha< \frac32.   
\end{align*}
\end{proof}
\begin{remark}
In \cite{trig} it has been proved that, near $t=0$, $\mathcal E_{C,n}(t,h)$ behaves like $t^{-\frac32}$ if $\alpha\in\left[\frac12,\frac34\right)$, and like $t^{-\frac{\alpha}{2(1-\alpha)}}$ if $\alpha\in \left[\frac34,1\right)$. If we compare this result with \eqref{stima_energia_damped_<1/2} when $\gamma=0$, we notice that we recover the same estimates as \cite{trig}.
\end{remark}

The second part of this subsection is devoted to estimate $\|Q_{t}^{-\frac12}e^{tA}\widetilde{G}a\|_{H}$ with $a\in U$. To this aim, we need an explicit expression for $h^+\in H^+$ and $h^-\in H^-$ which appear in the decomposition of $h\in H$. Notice that
\begin{align}
h
= & \left(\begin{matrix}
h_1 \\ h_2
\end{matrix}\right) 
= \left(\begin{matrix}
\displaystyle \sum_{n\in\N}(h_1)_n\frac{e_n}{\|e_n\|_U} \\[2mm]
\displaystyle \sum_{n\in\N}(h_2)_n\frac{e_n}{\|e_n\|_U}
\end{matrix}
\right)=\sum_{n\in\N}(\langle h^+,\Phi_n^+\rangle_H\Phi_n^++\langle h^-,\Phi_n^-\rangle_H\Phi_n^-) \notag \\
= & \sum_{n\in\N}\left(
\begin{matrix}
\mu_n^{\frac12}(\langle h^+,\Phi_n^+\rangle_H+\chi_n\langle h^-,\Phi_n^-\rangle_H)e_n \\[2mm]
(\lambda_n^+\langle h^+,\Phi_n^+\rangle_H+\chi_n\lambda_n^-\langle h^-,\Phi_n^-\rangle_H)e_n
\end{matrix}
\right).
\label{dec_h_+h_-}
\end{align}
By comparing the corresponding components in \eqref{dec_h_+h_-}, we infer that, for every $n\in\N$,
\begin{align*}
(h_1)_n
= & \mu_n^{\frac12}(\langle h^+,\Phi_n^+\rangle_H+\chi_n\langle h^-,\Phi_n^-\rangle_H)\|e_n\|_U,  \\
(h_2)_n
= & (\lambda_n^+\langle h^+,\Phi_n^+\rangle_H+\chi_n\lambda_n^-\langle h^-,\Phi_n^-\rangle_H)\|e_n\|_U.
\end{align*}
It follows that, for every $n\in\N$,
\begin{align}
\label{expl_h+h-}
\langle h^+,\Phi_n^+\rangle_H
= & \frac{\lambda_n^-(h_1)_n-\mu_n^{\frac12}(h_2)_n}{\mu_n^{\frac12}(\lambda_n^--\lambda_n^+)\|e_n\|_U},\qquad
\langle h^-,\Phi_n^-\rangle_H
=  \frac{\lambda_n^+(h_1)_n-\mu_n^{\frac12}(h_2)_n}{\chi_n\mu_n^{\frac12}(\lambda_n^+-\lambda_n^-)\|e_n\|_U}.
\end{align}

Now we are able to estimate the norm of $Q_{t}^{-\frac12}e^{tA}$ along the directions of $\widetilde{G}$.
\begin{thm}
\label{thm:damped_contr_>1/2_direct}
Let $\alpha\in\left[\frac12,1\right)$ and $\gamma\geq0$. Then, for every $t>0$ there exists a positive constant $\overline c$, which depends on $\alpha$ and $\gamma$ but is independent of $t$, if $t$ varies in a bounded interval, such that 
\begin{align}
\label{stima_energia_damped_>1/2_direct}
\mathcal E_{C}(t,\widetilde{G}a)
\leq \frac{\overline{c}\|\widetilde{G}a\|_{H}}{t^{\frac12+\frac\gamma{1-\alpha}}}, \qquad a\in U.
\end{align}
\end{thm}
\begin{proof}
Fix $a\in U$. To prove the estimate, we consider the control $u$ defined in \eqref{control}, with $h$ replaced by $\widetilde{G}a$ and $m$ in the definition of $\Phi_t$ which fulfills $2m-2\gamma/(1-\alpha)>-1$. We stress that, from Step $1$ in the proof of Theorem \ref{thm:damped_contr_>1/2}, we already know that $u$ steers $\widetilde{G}a$ at $0$ at time $t$. Hence, we only need to prove the estimate. As in the proof of Theorem \ref{thm:damped_contr_>1/2}, we split the proof into two steps. In the former we estimate the $L^2$-norm of $K_1\psi_t$, while in the latter we deal with $K_2\psi_t'$. At first, we recall that
\begin{align*}
h=\widetilde{G}a=\begin{pmatrix}
0 \\ \displaystyle \sum_{k=1}^{+\infty} a_k \frac{e_k}{\|e_k\|_U}    
\end{pmatrix}, \qquad a_k=\langle a,\frac{e_k}{\|e_k\|_U}\rangle_U, \quad k\in\N,    
\end{align*}
which replaced in \eqref{expl_h+h-} gives, for every $k\in\N$,
\begin{align}
\label{exp_Ga^+Ga^-}
\langle (G_na)^+,\Phi_k^+\rangle_H
= & \frac{-a_k}{(\lambda_k^--\lambda_k^+)\|e_k\|_U},\qquad
\langle (G_na)^-,\Phi_k^-\rangle_H
=  \frac{-a_k}{\chi_k(\lambda_k^+-\lambda_k^-)\|e_k\|_U}.    
\end{align}

{\bf Step $\bm1$}. From \eqref{forma_K_1psi_t} and taking \eqref{exp_Ga^+Ga^-} into account, we infer that
\begin{align*}
K_1\psi_t(\tau)
= & \Phi_t(\tau)\sum_{k=1}^{+\infty}\frac{\mu_k^\gamma a_k}{(\lambda^-_k-\lambda^+_k)}[e^{\lambda_k^+\tau}\lambda_k^--e^{\lambda_k^-\tau}\lambda_k^+]\frac{e_k}{\|e_k\|_U}, \qquad \tau\in(0,t),
\end{align*}
which gives
\begin{align*}
\|K_1\psi_t(\tau)\|_U^2
\leq & C|\Phi_t(\tau)|^2\sum_{k=1}^{+\infty} [\mu_k^{2\gamma}e^{-2\rho \mu_k^{1-\alpha}\tau}+\mu_k^{2\gamma+2-4\alpha}e^{-2\rho\mu_k^\alpha\tau}]a_k^2 \\
\leq & C|\Phi_t(\tau)|^2[\tau^{-\frac{2\gamma}{1-\alpha}}+\tau^{(-2\frac{\gamma+1-2\alpha}{\alpha})\wedge 0}]\sum_{k=1}^{+\infty}a_k^2 \\
\leq & C|\Phi_t(\tau)|^2\tau^{-\frac{2\gamma}{1-\alpha}}\|\widetilde{G}a\|_{H}^2, \qquad \tau\in(0,t),
\end{align*}
since if $\gamma\leq 2\alpha-1$ then the second exponent of $\tau$ is $0$, and 
\begin{align*}
\frac{2\gamma}{1-\alpha}\geq \frac{2\gamma+2-4\alpha}{1-\alpha} \geq \frac{2\gamma+2-4\alpha}{\alpha}\geq0, \qquad \gamma>2\alpha-1.    
\end{align*}
It follows that
\begin{align}
\label{stima_dirct_1}
\int_0^t\|K_1\psi_t(\tau)\|_U^2d\tau
\leq C t^{-2m-2}\int_0^t\tau^{2m}\tau^{-\frac{2\gamma}{1-\alpha}}d\tau \sum_{k=1}^{+\infty}a_k^2
\leq C t^{-1-\frac{2\gamma}{1-\alpha}}\|\widetilde{G}a\|_{H}^2.
\end{align}

{\bf Step $\bm2$}. Now we estimate $K_2\psi_t'$. Taking advantage from \eqref{forma_K_2psi_t'} and \eqref{exp_Ga^+Ga^-} it follows that
\begin{align*}
K_2\psi_t'(\tau)
= & \sum_{k=1}^{+\infty}\frac{\mu_k^{\gamma}}{(\lambda_k^--\lambda_k^+)\|e_k\|_U}[\Phi_t'(\tau)(e^{\lambda_k^+\tau}-e^{\lambda_k^-\tau})e_k+\Phi_t(\tau)(\lambda_k^+ e^{\lambda_k^+\tau}-\lambda_k^-e^{\lambda_k^-\tau})e_k]a_k \\
= & \sum_{k=1}^{+\infty}\frac{\mu_k^{\gamma}}{(\lambda_k^--\lambda_k^+)\|e_k\|_U}[\Phi_t'(\tau)e^{\lambda_k^+\tau}(1-e^{(\lambda_k^--\lambda_k^+)\tau})e_k+\Phi_t(\tau)(\lambda_k^+ e^{\lambda_k^+\tau}-\lambda_k^-e^{\lambda_k^-\tau})e_k]a_k
\end{align*}
for every $\tau\in(0,t)$, which gives
\begin{align*}
\|K_2\psi_t'(\tau)\|_U^2
\leq & C\sum_{k=1}^{+\infty} \mu_k^{2\gamma}\bigg(\left|\Phi_t'(\tau)\frac{1-e^{(\lambda_k^--\lambda_k^+)\tau}}{\lambda_k^--\lambda_k^+}\right|^2e^{-2\rho\mu_k^{1-\alpha}\tau} \\
& +|\Phi_t(\tau)|^2(\mu_k^{2(1-2\alpha)}e^{-2\rho\mu_k^{1-\alpha}\tau}+e^{-2\rho\mu_k^\alpha\tau})\bigg)a_k^2 \\
\leq & C\sum_{k=1}^{+\infty} \big(|\tau\Phi_t'(\tau)|^2\mu_k^{2\gamma}e^{-2\rho\mu_k^{1-\alpha}\tau}+|\Phi_t(\tau)|^2\mu_k^{2\gamma+2-4\alpha}e^{-2\rho\mu_k^{1-\alpha}\tau}+\mu_k^{2\gamma}e^{-2\rho\mu_k^\alpha\tau}\big)a_k^2 \\
\leq & C \big(|\tau^{-\frac{\gamma}{1-\alpha}+1}\Phi_t'(\tau)|^2+(\tau^{(-\frac{2\gamma+2-4\alpha}{1-\alpha})\wedge 0}+\tau^{-\frac{2\gamma}{\alpha}})|\Phi_t(\tau)|^2\bigg)\sum_{k=1}^{+\infty}a_k^2 \\
\leq & Ct^{-2m-2}  \big(\tau^{-\frac{2\gamma}{1-\alpha}+2m}+\tau^{2m}(\tau^{(-\frac{2\gamma+2-4\alpha}{1-\alpha})\wedge 0}+\tau^{-\frac{2\gamma}{\alpha}})\bigg)\|\widetilde{G}a\|_{H}^2 \\
\leq & C t^{-2m-2}\tau^{2m-\frac{2\gamma}{1-\alpha}}\|\widetilde{G}a\|_{H}^2, \qquad \tau\in(0,t),
\end{align*}
since $\gamma/(1-\alpha)\geq \gamma/\alpha$ for $\alpha\in\left[\frac12,1\right)$ and $2\gamma/(1-\alpha)\geq (2\gamma+2-4\alpha)/\alpha$ (see Step $1$). It thus follows that
\begin{align}
\label{stima_dirct_2}
\|K_2\psi_t'\|_U^2
\leq Ct^{-2m-2}\int_0^t\tau^{2m-\frac{2\gamma}{1-\alpha}}d\tau\|G_na\|_{H_n}^2 \leq Ct^{-1-\frac{2\gamma}{1-\alpha}}\|G_na\|_{H_n}^2.
\end{align}
Combining \eqref{stima_dirct_1} and \eqref{stima_dirct_2} we get the thesis.
\end{proof}

\begin{rmk}
Let us notice that the singularity which appears in the estimate of $\Gamma_{t}$ along the direction of $\widetilde{G}$ is integrable at $0$ if $\frac\gamma{1-\alpha}<\frac12$, i.e.,  if $\gamma<\frac12-\frac\alpha2$. \end{rmk}

\subsubsection{The main result}

\begin{thm}
\label{thm:damped_main_result}
Assume that: 
\begin{enumerate}[\rm(i)]
\item $\alpha\in\left[\frac12,1\right)$ and $\gamma\in\left[0,\frac12-\frac\alpha2\right)$;
\item $\theta\in\left(\frac23\cdot \frac{\gamma+1-\alpha}{1-\alpha},1\right)$ if $\gamma+2\alpha<\frac32$, $\theta\in\left(\frac{4\gamma+2\alpha-1}{2\gamma+\alpha},1\right)$ if $\gamma+2\alpha\geq \frac32$ and $C\in C_b^\theta(H;U)$.
\end{enumerate}
Then, if {\rm (a)} $\Lambda^{-2\gamma}:U\to U$ is a trace-class operator or {\rm (b)} $\delta>\frac{1}{2\gamma+\alpha}$ and
\begin{align}
\label{damped_cond_path_uniq_1}
\sum_{n\in\N}\mu_n^{-\alpha}\left\|\langle C(\cdot), \frac{e_n}{\|e_n\|_U}\rangle_U\right\|_{C_b^\theta(H)}^2<\infty,    
\end{align}
pathwise uniqueness holds true for \eqref{eq_damped_ex}.
\end{thm}
\begin{proof}
Let us notice that, under these assumptions on $\alpha,\gamma$, $\Lambda$ and $\delta$, from Proposition \ref{prop:damped_conv_stoc_1} Hypotheses \ref{hyp:finito-dimensionale}(iv)  is satisfied. Further, from Proposition \ref{prop:eq_ipotesi} and Theorems \ref{thm:damped_contr_>1/2} and \ref{thm:damped_contr_>1/2_direct}, Hypotheses \ref{hyp:finito-dimensionale}(vi) are verified. We also recall that, for $\alpha\in[\frac12,1)$, the operator $A$ generates a strongly continuous and analytic semigroup on $H$, and Remark \ref{rmk:damped_hyp_eigenvector-A*} shows that Hypotheses \ref{hyp:goal-addo}(a)-(b) hold true. If (a) is satisfied, i.e., $\Lambda^{-2\gamma}$ is a trace-class operator, then it immediately follows that $G=\widetilde G\Lambda^{-\gamma}$ is a Hilbert-Schmidt operator, and so Hypotheses \ref{hyp:traccia-finita} holds true.

It remains to check Hypotheses \ref{hyp:goal-addo}(c) under condition (b). From the definition of $\widetilde G$, $\Psi_n^+$ and $\Psi_n^-$ and recalling that $B=\widetilde G \widetilde B$, we get $\langle B(x),\Psi_n^+\rangle_H=\langle \widetilde G\widetilde B(x),\Psi_n^+\rangle_H=\lambda_n^+\langle \widetilde B(x),e_n\rangle_U$ and $\langle B(x),\Psi_n^-\rangle_H=\chi_n\lambda_n^-\langle C(x),e_n\rangle_U$ for every $n\in\N$. Hence, condition \eqref{conv_serie_holder} reads as
\begin{align}
\label{damped_serie_holder}
\sum_{n\in\N}\left(|\lambda_n^+|^2\|e_n\|_U^2\frac{\|\langle \widetilde B(\cdot),\frac{e_n}{\|e_n\|_U}\rangle_U\|_{C_b^\theta(H)}^2}{-{\rm Re}(\lambda_n^+)}+|\chi_n|^2|\lambda_n^-|^2\|e_n\|_U^2\frac{\|\langle \widetilde B(\cdot),\frac{e_n}{\|e_n\|_U}\rangle_U\|_{C_b^\theta(H)}^2}{-{\rm Re}(\lambda_n^-)}\right)<\infty.    
\end{align}
From \eqref{1stime_coeff_avl_2} we get $|\lambda_n^+|
\sim \mu_n^{1-\alpha}$, $|\lambda_n^-|
\sim \mu_n^\alpha$, $\|e_n\|_U\sim \mu_n^{-\frac12}$, $\chi_n\sim \mu_n^{\frac12-\alpha}$, $-{\rm Re}(\lambda_n^+)\sim \mu_n^{1-\alpha}$ and $-{\rm Re}(\lambda_n^+)\sim \mu_n^\alpha$ as $n$ goes to $\infty$. Therefore, the series in \eqref{damped_serie_holder}
behaves like
\begin{align*}
\sum_{n\in\N}\mu_n^{-\alpha}\left\|\langle \widetilde B(\cdot),\frac{e_n}{\|e_n\|_U}\rangle_U\right\|_{C_b^\theta(H)}^2.
\end{align*}
Arguing as in \eqref{hold_tilde_B}, we infer that 
\begin{align*}
\left\|\langle \widetilde B(\cdot),\frac{e_n}{\|e_n\|_U}\rangle_U\right\|_{C_b^\theta(H)}^2\leq \left\|\langle C(\cdot),\frac{e_n}{\|e_n\|_U}\rangle_U\right\|_{C_b^\theta(H)}^2(\|\Lambda^{-\frac12}\|_{\mathcal L(U)}^{2}+1)^{\theta}, \qquad n\in\N,
\end{align*}
which implies that if \eqref{damped_cond_path_uniq_1} holds true then \eqref{damped_serie_holder} is verified.
\end{proof}

We split the applications of Theorem \ref{thm:damped_main_result} into different statements. 

\begin{corollary}[Stochastic damped wave equation in dimension $1$]
\label{coro:damped_wave_dim_1}
If $U=L^2(0,\pi)$ and $\Lambda$ is minus the realization of the Laplace operator with homogeneous Dirichlet boundary conditions in $L^2(0,\pi)$, $\alpha\in\left[\frac12, 1\right)$, $\gamma\in\left(\frac14-\frac\alpha2,\frac12-\frac\alpha2\right)\cap[0,\infty)$, $\theta\in\left(\frac23\cdot\frac{\gamma+1-\alpha}{1-\alpha},1\right)$ if $\gamma+2\alpha<\frac32$ and $\theta \in\left(\frac{4\gamma+2\alpha-1}{2\gamma+\alpha},1\right)$ if $\gamma+2\alpha\geq \frac32$, $C\in C_b^\theta(H;U)$ and \eqref{damped_cond_path_uniq_1} is verified, then pathwise uniqueness for \eqref{eq_damped_ex} holds true.  
\end{corollary}
\begin{proof}
Under these assumptions we get $\mu_n\sim n^2$. Hence, $2>\frac{1}{2\gamma+\alpha}$ is satisfied if and only if $\gamma>\frac14-\frac\alpha2$. This implies that conditions (i), (ii) and (b) in Theorem \ref{thm:damped_main_result} are verified and so we get pathwise uniqueness for \eqref{eq_damped_ex}. 
\end{proof}

\begin{rmk}
Even if we set $\alpha=\frac12$, no choice of $\gamma<\frac12-\frac14=\frac14$ makes $\Lambda^{-2\gamma}$ a trace-class operator. Further, we do not cover dimension $2$ because in this case condition $\delta>\frac1{2\gamma+\alpha}$ reads as $2\gamma+\alpha>1$, which has empty intersection with $\gamma<\frac12-\frac\alpha2$.
\end{rmk}

\begin{rmk}
Corollary \ref{coro:damped_wave_dim_1} has important consequences.
\begin{enumerate}
    \item if $\alpha>\frac12$, then we can choose $\gamma=0$, which means that we can consider the white noise.
    \item If $\alpha>\frac12$ then \eqref{damped_cond_path_uniq_1} is satisfied for every $C\in C_b^\theta(H;U)$, since $\mu_n^{-\alpha}\sim n^{-2\alpha}$ whose series converges.
    \item If $\alpha=\frac12$ then the series of $\mu_n^{-\alpha}$ does not converge. Hence, a contribution from $C$ is needed in order to get \eqref{damped_cond_path_uniq_1}. If for every $(h,k)\in H=U\times U$ we have
\begin{align*}
C((h,k))(\xi)=g_1(\xi)\int_0^1f_1(\xi')\min\{h(\xi'),r\}^{\theta}d\xi'+g_2(\xi)\int_0^1f_2(\xi')\min\{k(\xi'),r\}^{\theta}d\xi'   
\end{align*}
for every $\xi\in[0,1]$, where $g_1,g_2,f_1,f_2\in L^\infty(0,1)$ and $r>0$ is fixed, then
\begin{align*}
\sum_{n\in\N} \left\|\langle C(\cdot),\frac{e_n}{\|e_n\|_U}\rangle_U\right\|_{C_b^\theta(H)}^2<\infty,   
\end{align*}
see for instance \cite{Dap-Fla2010} and the forthcoming Proposition \ref{prop:drift_holder_series}. 
\end{enumerate}
\end{rmk}

\begin{corollary}
\label{coro:beam_eq}
Let $m=1,2,3$, $U=L^2((0,\pi)^m)$ and $\Lambda=(-\Delta)^2$, where $\Delta$ is the realization of the Laplace operator with Dirichlet homogeneous boundary conditions in $L^2((0,\pi)^m)$, and assume that \eqref{damped_cond_path_uniq_1} is verified. Therefore, 
pathwise uniqueness for \eqref{eq_damped_ex} holds true when

\vspace{2mm}
\begin{center}
$\alpha\in \left[\frac12,1\right),\quad  \gamma\in\left(\frac m8-\frac\alpha2,\frac12-\frac\alpha2\right)\cap[0,\infty)$, 

\vspace{2mm}
$\theta\in \left(\frac23\cdot\frac{\gamma+1-\alpha}{1-\alpha},1\right) \textrm{ if $\gamma+2\alpha<\frac32$, }\theta\in \left(\frac{4\gamma+2\alpha-1}{2\gamma+\alpha},1\right) \textrm{ if $\gamma+2\alpha\geq\frac32$}$.
\end{center}
\end{corollary}

\begin{proof}
\begin{enumerate}[m=1.] 
\item  Let us notice that, in this situation, $\mu_n\sim n^4$. This means that $\delta=4$ and condition $\delta>\frac{1}{2\gamma+\alpha}$ reads as $\gamma>\frac18-\frac\alpha2$. Hence, the assumptions of Theorem \ref{thm:damped_main_result} are fulfilled and pathwise uniqueness follows.
\vspace{2mm}
\item  Under this condition, we get $\mu_n\sim n^2$, as for the stochastic damped wave equation in dimension $1$. Analogous computations to those in the proof of Corollary \ref{coro:damped_wave_dim_1} give the thesis.
\vspace{2mm}
\item The Laplace operator in $L^2((0,\pi)^3)$ has eigenvalues which behave like $n^{\frac23}$, whence $\mu_n\sim n^{\frac43}$. Let us notice that $\frac43>\frac{1}{2\gamma+\alpha}$ if and only if $\gamma>\frac{3}{8}-\frac{\alpha}{2}$. Therefore, the assumptions ensure that conditions (i), (ii) and (b) in Theorem \ref{thm:damped_main_result} are fulfilled. This implies that for \eqref{eq_damped_ex} pathwise uniqueness holds true.
\end{enumerate}    
\end{proof}

\begin{rmk}
If $m=1$, $\alpha\in\left[\frac12,\frac{3}{4}\right)$ and $\gamma\in \left(\frac{1}{8},\frac12 -\frac{\alpha}{2}\right)$ we can avoid assuming that \eqref{damped_cond_path_uniq_1} holds true, since $\Lambda^{-2\gamma}$ has finite trace.
\end{rmk}

\subsubsection{Counterexample to uniqueness  in the deterministic case for damped wave equation}
\label{subsec:counter}
In this section we show that the deterministic damped wave equation with H\"older continuous nonlinear term   could  be ill-posed. Therefore, Corollary \ref{coro:damped_wave_dim_1} 
gives  in fact results on the regularizing effect of the noise. 

We consider a semilinear deterministic damped wave equation in $L^2(0,\pi)$ with $\alpha=\frac7{12}$ and $\rho=1$, which fulfill Hypothesis \ref{cond_rho_mu_n_alpha}. 
 
Let us consider the semilinear deterministic equation
\begin{align}
\label{Count_det_damped_wave_eq_1}
\left\{
\begin{array}{lll}
\displaystyle \frac{\partial^2y}{\partial\tau^2}(\tau,\xi)
= \frac{\partial^2y}{\partial\xi^2}(\tau,\xi)-\left(-\frac{\partial^2}{\partial_\xi^2}\right)^{\frac7{12}} \!\frac{\partial y}{\partial\tau} (\tau,\xi)+c(\xi,y(\tau,\xi)), & \xi\in [0,\pi], & \tau\in[0,1],  \\
y(\tau,0)=y(\tau,\pi)=0, & & \tau\in[0,1], \\
\displaystyle y(0,\xi)=\frac{\partial y}{\partial t}(0,\xi)=0, & \xi \in [0,\pi],
\end{array}
\right.
\end{align}
where, for every $\xi\in[0,\pi]$ and $y\in\R$,
\begin{align*}
c(\xi,y)
=  & \varphi(y)\left(56({\rm sgn}(\sin(2\xi)))|\sin(2\xi)|^{\frac14}|y|^{\frac34} 
+8\cdot 4^{\frac7{12}}({\rm sgn}(\sin(2\xi)))|\sin(2\xi)|^{\frac18}|y|^{\frac7{8}} +4y\right).
\end{align*}
Here, $\varphi\in C_c^\infty(\R)$ satisfies $0\leq \varphi\leq 1$, $\varphi\equiv1$ in $(-2,2)$ and $\varphi\equiv0$ in $(-3,3)^c$.\\
We claim that equation \eqref{Count_det_damped_wave_eq_1} is not well-posed. Indeed, $e_n(x):=\sin (2n\xi)$ is an eigenfunction of $-\frac{\partial^2}{\partial \xi^2}$ and $-\frac{\partial^2}{\partial \xi^2} e_n=4n^2e_n$ for any $n\in\N$. Therefore, 
\begin{align*}
\left(-\frac{\partial^2}{\partial \xi^2}\right)^{\frac7{12}}\!e_n= 4^{\frac7{12}}n^{\frac7{6}}e_n, \qquad  n\in\N.
\end{align*}
$c$ is $\frac34$-H\"older continuous with respect to $y$, uniformly with respect to $\xi$, and both $y(\tau,\xi)=0$ and $y(\tau,\xi)=\tau^8\sin(2\xi)$ are solutions to \eqref{Count_det_damped_wave_eq_1}. However, if we perturb \eqref{Count_det_damped_wave_eq_1} by means of white noise, then the assumptions of Corollary \ref{coro:damped_wave_dim_1} are fulfilled with $\alpha=\frac7{12}$, $\gamma=0$ and $\theta=\frac{3}{4}$. Hence, pathwise uniqueness holds true for the stochastic version of \eqref{Count_det_damped_wave_eq_1}.

\subsection{Stochastic heat equation}\label{Heat-case}
In this subsection, we compare some results in \cite{Dap-Fla2010} with those in this paper. First of all, we note that Theorem \ref{pathwiseuniqueness} applies to all the examples contained in \cite{Dap-Fla2010}. 

Consider the SPDE introduced in \cite[Subsection 6.1]{Dap-Fla2010} given by
\begin{align}\label{eqFObeta}
\left\{
\begin{array}{ll}
\displaystyle  dX(t)=-(-\Delta)^{\beta}X(t)dt+B(X(t))dt+(-\Delta)^{-\gamma/2}dW(t), \quad t\in[0,T],  \vspace{1mm} \\
X(0)=x\in H,
\end{array}
\right.
\end{align}
where $B\in C^\theta_b(H;H)$, for some $\theta\in (0,1)$, $\beta,\gamma\geq 0$ and $\Delta$ is the realization of the Laplace operator with periodic boundary conditions in $H=L^2([0,2\pi]^m)$ with $m=1,2,3$. We are going to show that pathwise uniqueness holds true for SPDE \eqref{eqFObeta} with less restrictive assumptions on $\beta$ and $\gamma$ than those assumed in \cite{Dap-Fla2010}. In particular, the hypotheses of this paper cover the case of the stochastic heat equation in $L^2([0,2\pi]^3)$ which is not contained in the assumptions of \cite{Dap-Fla2010}.

We recall that there exists an orthonormal basis $\{e_k:k\in\N\}$ of $H$ consisting in eigenvectors of $\Delta$. So the spaces $\{H_n\}_{n\in\N}$ given by
\[
H_0=\{0\},\quad H_n:={\rm span}\{e_1,...,e_n\},\quad n\in\N,
\]
so Hypotheses \ref{hyp:finito-dimensionale}(v) and Hypotheses \ref{hyp:goal-addo}(a)-(b) (with $d_n=1$ for every $n\in\N$) hold true. Moreover, for every $k\in\N$, we have
\begin{equation}\label{autovalori}
\Delta e_k=-\lambda_ke_k, \quad \lambda_k\sim k^{\frac2m}.
\end{equation}
By easy computations, for every $t>0$ and $n\in\N$ we obtain that
\begin{align*}
& Q_t:=\int_0^te^{-2s(-\Delta)^{\beta}}(-\Delta)^{-\gamma}ds=\frac12(-\Delta)^{-(\beta+\gamma)}(\Id_H-e^{-2t(-\Delta)^{\beta}}).
\end{align*}
\begin{proposition}\label{c-hyp1} $ $
  \begin{enumerate}
      \item If $(m-2\gamma)/2\beta<1$, then there exists $\eta>0$ such that for every $t>0$ we have
      \begin{equation}\label{stima-heat}
          \int_0^t{\rm Trace}_H\frac{1}{s^\eta}\left[e^{-2s(-\Delta)^{\beta}}(-\Delta)^{-\gamma}\right]ds<\infty.
      \end{equation}
      \item There exists a constant $c>0$ such that for every $n\in\N$ and $t>0$ we have
\begin{align*}
 \norm{\Gamma_t}_{\mathcal{L}(H)}=\norm{Q^{-\frac12}_te^{-t(-\Delta)^{\beta}}}_{\mathcal{L}(H)}\leq \frac{c}{t^{\frac12+\frac{\gamma}{2\beta}}}.
\end{align*}
      \item If $(m-2\beta)/2<\gamma< \beta\theta/(2-\theta)$, then Hypotheses \ref{hyp:finito-dimensionale}, with $\widetilde G={\rm Id}_H$ and $\widetilde B=B$, hold true.
  \end{enumerate}
  \end{proposition}
\begin{proof}
  \begin{enumerate}
      \item Let $t>0$ and $\eta\in(0,1)$. By \eqref{autovalori} we have
\begin{align*}
    \int_0^t\frac{1}{s^\eta}{\rm Trace}_H\left[e^{-2s(-\Delta)^{\beta}}(-\Delta)^{-\gamma}\right]ds&\leq C_0\int^t_0\frac{1}{s^\eta}\sum_{k=1}^{\infty}e^{-2sk^{2\beta/m}}k^{-2\gamma/m}ds\\
    &\leq C_1\int^t_0\frac{1}{s^\eta}\int_{1}^{\infty}e^{-2sx^{2\beta/m}}x^{-2\gamma/m}dxds\\
    &\leq C_2\int^t_0\frac{1}{s^{(m+\eta-2\gamma)/2\beta}}\int_{s}^{\infty}e^{-2y}y^{(m-2\gamma-2\beta)/2\beta}dyds,
\end{align*}
where $C_0,C_1,C_2$ are positive constants. So by the same arguments used in \eqref{conto-damped} and \eqref{conto-damped2}, \eqref{stima-heat} holds true for every $(m-2\gamma)/2\beta<1$ and $\eta<\eta_0$ for some $\eta_0<1$.

\item Follows by \cite[Proposition 2.1.1]{Lun95}.
      
\item Follows combining points $(1)$ and $(2)$
  \end{enumerate}
    
\end{proof}

Let $g,h\in L^{\infty}([0,2\pi]^m)$, $r>0$ and $\theta\in (0,1)$. Consider the function $B_r:H\rightarrow H$ given by
\[
B_r(f)(\xi):=g(\xi)\int_{[0,2\pi]^m}h(\xi')\min\{f(\xi') ,r\}^{\theta}d\xi',\quad f\in H.
\]
By a slight modification of \cite[Lemma 8]{Dap-Fla2010} we obtain the following result.
\begin{proposition}
\label{prop:drift_holder_series}
    $B_r$ belongs to $C_b^\theta(H;H)$ and verifies \eqref{conv_serie_holder}.
\end{proposition}
Using Theorem \ref{pathwiseuniqueness} in the case of the Laplace operator (namely $\beta=1$ in the SPDE \eqref{eqFObeta}) we deduce the following result.
\begin{thm}\label{Heatequation}
Assume that $\beta=1$ and $B=B_r$. Pathwise uniqueness holds true for SPDE \eqref{eqFObeta} in the following cases:
\begin{align*}
& 0\leq \gamma < \frac{\theta}{2-\theta}, && 0<\theta<1,&& m=1; \\
& 0< \gamma < \frac{\theta}{2-\theta}, && 0<\theta<1, && m=2; \\
& \frac12< \gamma < \frac{\theta}{2-\theta},&&\frac{2}{3}<\theta<1, && m=3. \\
\end{align*}
\end{thm}

Theorem \ref{Heatequation} covers the case $m=3$ which, instead, is not contained in \cite[Proposition 10]{Dap-Fla2010}.

\begin{proposition}
Assume that $(m-2\beta)/2<\gamma< \beta\theta/(2-\theta)$ and that $B=(-\Delta)^{-\gamma/2}F$ with $F\in C^\theta_b(H;H)$. Then Hypotheses \ref{hyp:finito-dimensionale} and \ref{hyp:goal-addo} hold true and so pathwise uniqueness holds true for SPDE \eqref{eqFObeta}.
\end{proposition}
\begin{proof}\label{structure-heat}
Since  $(m-2\beta)/2<\gamma< \beta\theta/(2-\theta)$, by Proposition \ref{c-hyp1} it follows that Hypotheses \ref{hyp:finito-dimensionale} hold true. In this framework, the series in \eqref{conv_serie_holder} reads as
\[
\sum_{k=1}^{\infty}\dfrac{\|\scal{(-\Delta)^{-\gamma/2}F}{e_k}\|^2_{C^\theta_b(H)}}{\lambda^\beta_k}=\sum_{k=1}^{\infty}\dfrac{\|\scal{F}{e_k}\|^2_{C^\theta_b(H)}}{\lambda_k^{\beta+\gamma}}.
\]
From \eqref{autovalori} and $(m-2\beta)/2<\gamma$, we obtain 
\[
\sum_{k=1}^{\infty}\dfrac{\|\scal{(-\Delta)^{-\gamma/2}F}{e_k}\|^2_{C^\theta_b(H)}}{\lambda^\beta_k}\sim \|F\|_{C^\theta_b(H;H)}\sum_{k=1}^{\infty}\dfrac{1}{k^{(2\beta+2\gamma)/m}}<\infty.
\]
Hence, Hypotheses \ref{hyp:goal-addo}(c) holds true.
\end{proof}

\begin{rmk}
Let $\mathcal{O}$ be a bounded subset of $R^m$. It is possible to extend this example to an operator $A$ which is the realization in $L^2(\mathcal{O})$ of a general second-order differential operator, see for instance the class of operators defined in \cite[Section 6.1]{Cer2001}. This generalization is due to the fact that we remove the assumption that $A$ is self-adjoint, which is instead considered in \cite{Dap-Fla2010}.
\end{rmk}

\begin{appendices}
\section{}\label{sec:Equazioni-Kolmogorov}

Let $\X$ be a real separable Hilbert space. In this section, we recall some preliminary results about $L^2$-{\it maximal regularity} for analytic semigroups and  Ornstein-Uhlenbeck semigroups. The results in this section are known, see \cite{desimon64,weis01} for $L^2$-{\it maximal regularity}  and see \cite{Ang-Big-Fer2023, Big-Fer2021, Big-Fer2023,Cer-Dap2012, Cer-Lun2019, Cer-Lun2021, Pri2009, Pri-Zam2000} for properties of  Ornstein-Uhlenbeck type semigroup in infinite dimension.  Some of the following results are known in the literature; however, in order to apply the method presented in this paper, it is essential to have explicit constants in the estimates provided in this appendix. Therefore, we have chosen to include the proof of all the estimates we will use.


\subsection{Maximal $L^2$-regularity}\label{MaxL2}
In this subsection, we consider the complexification of $\X$, and we still denote it by $\X$. Let $\C:{\rm Dom}(\C)\subseteq{\X}\rightarrow\X$ be the infinitesimal generator of a strongly continuous and analytic semigroup $\{e^{t\C}\}_{t\geq 0}$ on $\X$.
The strong continuity of $\{e^{t\C}\}_{t\geq 0}$ implies that there exist $\omega\in\R$ and $M\geq 1$ such that $\|e^{t\C}\|_{\mathcal L(\X)}\leq Me^{\omega t}$ for every $t\geq0$. From the analyticity of $\{e^{t\C}\}_{t\geq 0}$ we deduce that there exist $\theta\in(\frac{\pi}2,\pi)$ and $c>0$ such that
\begin{align}
\notag
& \rho(C)\supset S_{\theta,\omega}:=\{\lambda\in\mathbb C: \lambda\neq \omega, \ |{\rm arg}(z-\omega)|<\theta\}, \\    
& \|R(\lambda,\C)\|_{\mathcal L(\X)}
\leq \frac{c}{|\lambda-\omega|} \qquad \forall \lambda\in S_{\theta,\omega}.
\label{stima_risolvente}
\end{align}
In particular, this implies that for every $\zeta>\omega$ there exist $c_{\zeta,1},c_{\zeta,2}>0$ such that
\begin{align}
\label{stima_ris_per_fourier_1}
\|R(\lambda,\C)\|_{\mathcal L(\X)}
\leq c_{\zeta,1}, \qquad \|\C R(\lambda,\C)\|_{\mathcal L(\X)}
\leq c_{\zeta,2} \qquad \forall \lambda\in\{z\in \mathbb{C}:{\rm Re}z\geq \zeta\}.
\end{align}
Indeed, from the definition of $R(\lambda,\C)$ it follows that for every $\zeta>\omega$ we get
\begin{align*}
& \|R(\lambda,\C)\|_{\mathcal L(\X)}
\leq \frac{c}{|\lambda-\omega|}\leq \frac{c}{\zeta-\omega}=:c_{\zeta,1}, \\
& \|\C R(\lambda,\C)\|_{\mathcal L(\X)}
= \|{\rm Id}-\lambda R(\lambda,\C)\|_{\mathcal L(\X)}\leq 1+\frac{c|\lambda|}{|\lambda-\omega|}
\leq c_{\zeta,2} \qquad \forall \lambda\in \mathbb C, \ {\rm Re}\lambda \geq \zeta, 
\end{align*}
for a suitable positive constant $c_{\zeta,2}$. 

The following optimal estimate has been already obtained in \cite{desimon64, weis01}. We provide the proof since we need to show that the constant which appears in the final estimate only depends on the constants which appear in \eqref{stima_ris_per_fourier_1}.

We introduce the Fourier transform for vector-valued functions. Given $f\in L^1(\R;\X)$, the Fourier transform of $f$ is given by
\begin{align*}
\mathcal F(f)(z)=\int_{\R}f(s)e^{-isz}ds \qquad \forall z\in \mathbb C.    
\end{align*}
It is well-known that $\mathcal F$ maps $L^1(\R;\X)\cap L^2(\R;\X)$ into $L^2(\R;\X)$ and $\|\mathcal F(f)\|_{L^2(\R;\X)}=\sqrt{2\pi}\|f\|_{L^2(\R;\X)}$. Further, for every $f,g\in L^2(\R;\X)$ we get $\mathcal F(f*g)=\mathcal F(f)\mathcal F(g)$, where
\begin{align*}
(f*g)(t)=\int_{\R}\langle f(s),g(t-s)\rangle_{\X} ds  \qquad \textrm{a.e. }t\in \R. 
\end{align*}

\begin{proposition}
\label{lemm:fourier}
For every $f\in L^2(0,T;\X)$ and every $\zeta>\omega$ the function $g:[0,T]\to \X$, defined as
\begin{align*}
g(t):=\int_0^te^{(t-s)\C}f(s)ds \qquad \forall t\in[0,T],
\end{align*}
belongs to $L^2([0,T];D(\C))$ and

\begin{align}
\label{stima_conv_fourier_gen}
\|e^{-\zeta\cdot}g\|_{L^2([0,T];D(\C))}\leq 2\pi(c_{\zeta,1}+c_{\zeta,2})\|e^{-\zeta\cdot}f\|_{L^2([0,T];\X)},
\end{align}
where $c_{\zeta,1}$ and $c_{\zeta,2}$ are the constants introduced in 
\eqref{stima_ris_per_fourier_1}.
In particular, \eqref{stima_conv_fourier_gen} gives
\begin{align}
\label{stima_conv_fourier_part}
\|g\|_{L^2([0,T];D(\C))}\leq 2\pi(c_{\zeta,1}+c_{\zeta,2})e^{2|\zeta| T}\|f\|_{L^2([0,T];\X)}. \end{align}
\end{proposition}
\begin{proof}
Let us prove \eqref{stima_conv_fourier_gen}, since \eqref{stima_conv_fourier_part} immediately follows from it. If we consider the trivial extension $\widetilde f$ of $f$ on $\mathbb R$ and we define $\widetilde g$ as $g$ with $f$ replaced by $\widetilde f$, then $\widetilde S$ is defined on the whole $\mathbb R$. Let us consider the Fourier transform of $t\mapsto e^{-\zeta t}\widetilde g(t)$; since $e^{-\zeta \cdot}\widetilde g= (e^{-\zeta \cdot}e^{\cdot \C})* (e^{-\zeta \cdot}\widetilde f)$, the properties of the Fourier transform give
\begin{align*}
\mathcal F(e^{-\zeta\cdot}\widetilde g)(z)
=
\mathcal F(e^{-\zeta \cdot} e^{\cdot \C})(z)\mathcal F(e^{-\zeta\cdot}\widetilde f)(z)
\end{align*}
for every $z\in \mathbb C$, where we have set $e^{t\C}=0$ for $t<0$. Let us notice that for every $\eta\in\mathbb R$ we have
\begin{align*}
\mathcal F(e^{-\zeta \cdot}e^{\cdot \C})(\eta)
= 
\int_{\mathbb R}
e^{t\C}e^{-(\zeta+i\eta) t}dt,
\end{align*}
and we recall that
\begin{align*}
R(\lambda,\C)=\int_0^\infty e^{-\lambda t}e^{t\C}dt   
\end{align*}
for every $\lambda\in \mathbb C$ with ${\rm Re}\lambda\geq \zeta$. We stress that these integrals are well-defined from the assumptions on $e^{t\C}$ and $\omega$, and so $\mathcal F(e^{-\zeta\cdot}e^{\cdot \C})(\eta)=R(\lambda,\C)$ with $\lambda=\zeta+i\eta$. Further, for every $\lambda\in \mathbb C$ with ${\rm Re}\lambda\geq \zeta$, from 
\eqref{stima_ris_per_fourier_1} we infer that
\begin{align*}
\|R(\lambda,\C)\|_{\mathcal L(\X;D(\C))}
= \|R(\lambda,\C)\|_{\mathcal L(\X)}+\|\C R(\lambda,\C)\|_{\mathcal L(\X)}
\leq c_{\zeta,1}+c_{\zeta,2}.
\end{align*}

It follows that
\begin{align*}
\|\mathcal F(e^{-\zeta\cdot}\widetilde g)(\eta)\|_{D(\C)}^2
\leq 
(c_{\zeta,1}+c_{\zeta,2})^2\|\mathcal F(e^{-\zeta\cdot}\widetilde f)(\eta)\|_{L^2(\X)}^2
\end{align*}
for every $\eta\in\mathbb R$. Integrating on $\mathbb R$ with respect to $\eta$ and recalling that $\|\mathcal F(h)\|_{L^2(\R;\X)}=\sqrt{2\pi}\|h\|_{L^2(\R;\X)}$ for every $h\in L^1(\R;\X)\cap L^2(\R;\X)$, we conclude that
\begin{align*}
\|e^{-\zeta\cdot}\widetilde g\|_{L^2(\mathbb R;D(\C))}^2
\leq 2\pi(c_{\zeta,1}+c_{\zeta,2})^2\|e^{-\zeta\cdot }\widetilde f\|_{L^2(\mathbb R;\X)}^2,
\end{align*}
which gives the thesis since both $\widetilde g$ and $\widetilde f$ vanish for $t<0$ and $t>T$.
\end{proof}

\subsection{The Ornstein-Uhlenbeck semigroup on $B_b(\X)$}
Let $\C:{\rm Dom}(\C)\subseteq \X\rightarrow \X$ be the infinitesimal generator of a strongly continuous semigroup $\{e^{t\C}\}_{t\geq0}$ on $\X$. Let $Q: \X\rightarrow \X$ be a linear bounded self-adjoint non-negative operator.

\begin{hyp1}\label{controllabilità}
For every $t>0$ the following conditions hold true:
\begin{align*}
&{\rm Trace}[Q_t]<\infty, \qquad e^{t\C}(\X)\subseteq Q_t^{\frac12}(\X), \qquad Q_t:=\int^t_0e^{s\C}Qe^{s\C^*}ds.
\end{align*}
\end{hyp1}

We introduce the Ornstein-Uhlenbeck semigroup $\{R(t)\}_{t\geq 0}$ given by
\begin{equation}\label{OUS}
(R(t)\varphi)(x)=\int_{\X}\varphi(e^{t\C}x+y)\mu_t(dy),\quad t>0,\;\varphi\in B_b(\X),\; x\in \X,
\end{equation}
where $\mu_t$ is the Gaussian measure on $\mathcal{B}(\X)$ with mean $0$ and covariance operator $Q_t$.
Under Hypotheses \ref{controllabilità} the semigroup $\{R(t)\}_{t\geq 0}$ verifies some regularity properties that we state in the subsequent propositions.
For any $t>0$ we set
\begin{equation}\label{gammat}
\Gamma_t=Q^{-\frac12}_t e^{t\C}.
\end{equation}
\begin{proposition}[Theorem 6.2.2 and Proposition 6.2.9 of \cite{Dap-Zab2002}]\label{stimeclassiche}
Assume that Hypotheses \ref{controllabilità} hold true. Then
\[
R(t)(B_b(\X))\subseteq C_b^{\infty}(\X),\quad t>0.
\]
For every $\varphi\in B_b(\X)$, $t>0$ and $x,h,k\in \X$ we have
\begin{align*}
D(R(t)\varphi)(x)h&=\int_{\X}\scal{\Gamma_{t}h}{Q_{t}^{-\frac12}y}_H\varphi(e^{t\mathcal{C}}x+y)\mu_{t}(dy),\\
D^2(R(t)\varphi)(x)(h,k)&=\int_\X\Big(\scal{\Gamma_{t}h}{Q_{t}^{-\frac12}y}_\X\scal{\Gamma_{t}k}{Q_{t}^{-\frac12}y}_\X-\scal{\Gamma_{t}h}{\Gamma_{t}k}_\X\Big)\varphi(e^{t\mathcal{C}}x+y)\mu_{t}(dy),
\end{align*}
and if, in addition, $\varphi\in C^1_b(H_n)$, then
\begin{align*}
DR(t)\varphi(x)h&=\int_\X\scal{\nabla \varphi(e^{t\mathcal{C}}x+y)}{e^{t\mathcal{C}}h}\mu_{t}(dy),\\
D^2R(t)\varphi(x)(h,k)&=\int_\X\scal{\Gamma_{t}h}{Q_{t}^{-\frac12}y}_H\scal{\nabla \varphi(e^{t\mathcal{C}}x+y)}{e^{t\mathcal{C}}k}\mu_{t}(dy).
\end{align*}
\end{proposition}

By the previous proposition we easily deduce the following estimates.
\begin{proposition}
Assume that Hypotheses \ref{controllabilità} hold true. For every $\varphi\in B_b(\X)$, $t>0$ and $x,h,k\in \X$ we have
\begin{align*}
   &|D(R(t)\varphi)(x)h|\leq \|\Gamma_th\|_\X\|f\|_{\infty},\\
   &|D^2(R(t)\varphi)(x)(h,k)|\leq \sqrt{2}\|\Gamma_th\|_\X\|\Gamma_tk\|_\X\|\varphi\|_{\infty}.
\end{align*}
For every $\varphi\in UC_b^1(\X)$, $t>0$ and $x,h,k\in \X$ we have
\begin{align*}
   &|D(R(t)\varphi)(x)h|\leq \|e^{t\C}h\|_{\X}\|f\|_{C^1_b(\X)},\\
   &|D^2(R(t)\varphi)(x)(h,k)|\leq \|e^{t\C}k\|_{\X}\|\Gamma_th\|_\X\|\varphi\|_{C^1_b(\X)}.
\end{align*}
\end{proposition}
By interpolation, see for instance \cite[Proposition 2.3.3]{Dap-Zab2002}, we deduce the following result.
\begin{proposition}\label{SchauderR}
Assume that Hypotheses \ref{controllabilità} hold true, and let $\theta\in (0,1)$. For every $\varphi\in C_b^\theta(\X)$ and $x,h,k\in\X$ we have
\begin{align*}
   &|D(R(t)\varphi)(x)h|\leq \|e^{t\C}h\|_{\X}^\theta \|\Gamma_th\|^{1-\theta}_\X\|\varphi\|_{C^\theta_b(\X)},\\
   &|D^2(R(t)\varphi)(x)(h,k)|\leq 2^{(1-\theta)/2}\|e^{t\C}k\|_{\X}^\theta\|\Gamma_tk\|_\X^{1-\theta}\|\Gamma_th\|_\X\|\varphi\|_{C^\theta_b(\X)}.
\end{align*}
\end{proposition}

\subsection{The Ornstein-Uhlenbeck semigroup on $B_b(\X;\X)$}
Let $E$ be a separable Hilbert space, let $\{W(t)\}_{t\geq 0}$ be a $E$-cylindrical Wiener process on a normal filtered probability space $(\Omega,\mathcal{F},\{\mathcal{F}_t\}_{t\geq 0},\mathbb{P})$ and let $\mathcal{D}\in\mathcal{L}(E;\X)$. We consider the SPDE 
\begin{align}\label{eqFOL}
\left\{
\begin{array}{ll}
\displaystyle  Z(t)=\C Z(t)dt+\mathcal{D}dW(t), \quad t> 0, \vspace{1mm} \\
Z(0)=x\in \X,
\end{array}
\right.
\end{align}
and we set $Q=\mathcal{D}\mathcal{D}^*:\X\rightarrow\X$.
By Hypotheses \ref{controllabilità}, for every $x\in \X$ the SPDE \eqref{eqFOL} has unique mild solution $\{Z(t)\}_{t\geq 0}$ given by
\[
Z(t)=e^{t\C}x+W_\C(t),\quad {\rm \mathbb{P}-a.s.},\quad \forall\; t>0,
\]
where $\{W_\C(t)\}_{t\geq 0}$ is the stochastic convolution process defined by 
\begin{equation*}
    W_\C (t):=\int^t_0e^{(t-s)\C}\mathcal{D}dW(s), \qquad \mathbb P{\rm -a.s.}, \ \forall t>0.
\end{equation*}
$\{W_\C(t)\}_{t\geq 0}$ is a $\X$-valued Gaussian process and, for every $t\geq0$, the random variable $W_\C(t)$ is a Gaussian random variable with mean $0$ and covariance operator $Q_t$. For an in-depth study of \eqref{eqFOL} we refer to \cite[Chapter 5]{Dap-Zab14}.

We define the vector valued Ornstein--Uhlenbeck semigroup $\{\mathcal{R}(t)\}_{t\geq 0}$ on the space $B_b(\X;\X)$ as
\begin{equation*}
(\mathcal{R}(t)\Phi)(x)=\mathbb{E}\left[\Phi(Z(t,x))\right],\quad \Phi\in B_b(\X;\X),\; t>0,\; x\in \X.
\end{equation*}
Let $\Phi\in B_b(\X;\X)$. For every $v\in \X$ we set
\[
\phi_v(x):=\scal{\Phi(x)}{v}_\X,\qquad x\in \X.
\]
In \cite[Section 3]{Dap-Fla2010} it is proved that 
\[
\mathcal{R}(t)(B_b(\X;\X))\subseteq C^{\infty}_b(\X;\X),\quad t>0,
\]
and for every $\Phi\in B_b(\X;\X)$, $t>0$ and $v,h,k,x\in \X$
\begin{align}
   &\scal{(\mathcal{R}(t)\Phi)(x)}{v}_\X=(R(t)\phi_v)(x),\label{vetOr}\\
   &\scal{D(\mathcal{R}(t)\Phi)(x)h}{v}_\X=D (R(t)\phi_v)(x)h,\label{deri1vet}\\
   &\scal{D^2(\mathcal{R}(t)\Phi)(x)(h,k)}{v}_\X=D^2(R(t)\phi_v)(x)(h,k)\label{derivata2vet}.
\end{align}

\begin{proposition}\label{SchauderRVV}
Assume that Hypotheses \ref{controllabilità} hold true, and let $\theta\in (0,1)$. For every $\Phi\in C_b^\theta(\X;\X)$, $t>0$ and $x,h,k\in\X$ we have
\begin{align*}
   &\|D(\mathcal{R}(t)\Phi)(x)h\|_\X\leq \|e^{t\C}h\|_{\X}^\theta \|\Gamma_th\|^{1-\theta}_\X\|\Phi\|_{C^\theta_b(\X;\X)},\\
   &\|D^2(\mathcal{R}(t)\Phi)(x)(h,k)\|_X\leq 2^{(1-\theta)/2}\|e^{t\C}k\|_{\X}^\theta\|\Gamma_tk\|_\X^{1-\theta}\|\Gamma_th\|_\X\|\Phi\|_{C^\theta_b(\X;\X)}.
\end{align*}
\end{proposition}
\begin{proof}
We prove the statement for the second derivative. We fix $\theta\in (0,1)$, $\Phi\in C_b^\theta(\X;\X)$, $t>0$ and $x,h,k,v\in \X$. From Proposition \ref{SchauderR} and \eqref{derivata2vet} we have
\begin{align*}
   |\scal{D^2(\mathcal{R}(t)\Phi)(x)(h,k)}{v}_\X|&=|D^2(R(t)\phi_v)(x)(h,k)|\leq 2^{(1-\theta)/2}\|e^{t\C}k\|_{\X}^\theta\|\Gamma_tk\|_\X^{1-\theta}\|\Gamma_th\|_\X \|\phi_v\|_{C^\theta_b(\X)}.
\end{align*}
Since $||\phi_v||_{C^\theta_b(\X)}\leq ||\Phi||_{C^\theta_b(\X;\X)}||v||_\X$, we obtain
\begin{align*}
   |\scal{D^2(\mathcal{R}(t)\Phi)(x)(h,k)}{v}_\X|\leq 2^{(1-\theta)/2}\|e^{t\C}k\|_{\X}^\theta\|\Gamma_tk\|_\X^{1-\theta}\|\Gamma_th\|_\X\|v\|_\X\|\Phi\|_{C^\theta_b(\X;\X)},
\end{align*}
and so we conclude 
\begin{align*}
  \|D^2(\mathcal{R}(t)\Phi)(x)(h,k)\|_\X\leq 2^{(1-\theta)/2}\|e^{t\C}k\|_{\X}^\theta\|\Gamma_tk\|_\X^{1-\theta}\|\Gamma_th\|_\X\|\Phi\|_{C^\theta_b(\X;\X)}.
\end{align*}\end{proof}

\subsection{Backward Kolmogorov equation}
Let $\mathcal{I}\in\mathcal{L}(E;\X)$. Let $T>0$, $\mathcal{N}\in C^{\theta}_b(\X;E)$ and $\mathcal{M}\in C^{\theta}_b(\X;\X)$ with fixed $\theta\in(0,1)$. We look for a solution $U:[0,T]\times \X\to \X$ in $C^{0,1}([0,T]\times\X;\X)$ to the integral equation 
\begin{equation}\label{Back-Kolmo}
    U(t,x)=\int^T_t \mathcal{R}(r-t)\left( DU(r,\cdot)\mathcal{I}\mathcal{N}(\cdot)+\mathcal{M}(\cdot)\right)(x)dr ,\quad t>0, \ x\in \X.
\end{equation}
We note that $DU(t,\cdot):\X\rightarrow \mathcal{L}(\X)$ and $DU(t,\cdot)\mathcal{I}\mathcal{N}(\cdot):\X\rightarrow \X$ for every $t\in[0,T]$. Further
\begin{equation}\label{RderiU}
U(t,x)\mathcal{I}\mathcal{N}(x)=\sum^{\infty}_{k=1}\scal{\mathcal{I}\mathcal{N}(x)}{g_k}_\X DU(t,x)g_k, \qquad t\in[0,T], \ x\in \X,
\end{equation}
for any orthonormal basis $\{g_k : k\in\N\}$ of $\X$.
Moreover setting $U_k=\scal{U}{g_k}_\X$ by \eqref{vetOr}, \eqref{deri1vet} and \eqref{RderiU} we have
\begin{equation}\label{RderiU-var}
\scal{DU(t,x)\mathcal{I}\mathcal{N}(x)}{g_k}=\scal{\mathcal{I}\mathcal{N}(x)}{\nabla U_k(t,x)}_\X, \qquad k\in\N,\; t\in[0,T], \; x\in \X.
\end{equation}
Before studying the integral equation \eqref{Back-Kolmo}, we focus on functions $\Phi\in C^1_b(\X;\X)$ such that the map $x\in\X\rightarrow D\Phi(x)\mathcal{I}\in\mathcal{L}(E;X)$ belongs to $C^1_b(\X;\mathcal{L}(E;\X))$.
First of all we note that if $\Phi\in C^2_b(\X;\X)$ then $D\Phi(\cdot)\mathcal{I}\in C^1_b(\X;\mathcal{L}(E;\X))$ and
\[
D\left[D\Phi(x)\mathcal{I}v\right]k=D^2\Phi(\mathcal{I}v,k),\qquad x,k\in\X,\; v\in E.
\]
We refer to \cite{Big-Fer-For-Zan2024} for a review about differentiability along suitable directions.
Finally by Proposition \ref{SchauderRVV} (with $h=\mathcal{I}v$), we deduce the following estimate.

\begin{corollary}\label{SchauderRV}
Assume that Hypotheses \ref{controllabilità} hold true, and let $\theta\in (0,1)$. For every $t>0$, $x,h\in\X$, $v\in E$ and $\Phi\in C^\theta_b(\X;\X)$ we have
\begin{equation*}
\norm{\mathcal{R}(t)\Phi}_{C^1_b(\X;\X)}+\|D\mathcal{R}(t)\Phi(\cdot)\mathcal{I}\|_{C^1_b(\X;\mathcal{L}(E,\X))}\leq K_t\norm{\Phi}_{C^\theta_b(\X;\X)}
\end{equation*}
where
\[
K_t:=1+(1+\norm{\mathcal{I}}_{\mathcal{L}(E;\X)})\|e^{t\mathcal{C}}\|^\theta_{\mathcal{L}(\X;\X)}\|\Gamma_t\|^{1-\theta}_{\mathcal{L}(\X;\X)}+\|e^{t\mathcal{C}}\|^\theta_{\mathcal{L}(\X;\X)}\|\Gamma_t\|^{1-\theta}_{\mathcal{L}(\X;\X)}\|\Gamma_t\mathcal{I}\|_{\mathcal{L}(E;\X)}.
\]
\end{corollary}

Now we can prove that \eqref{Back-Kolmo} is well-posed, but to do so we need an additional assumption (see also Hypotheses \ref{hyp:finito-dimensionale}(vi)).
\begin{hyp1}\label{supercontrollabilità}
Let $\Gamma_t$ be defined by \eqref{gammat}. We assume that Hypotheses \ref{controllabilità} hold true and that 
\[
\int^T_0\|\Gamma_t\|^{1-\theta}_{\mathcal{L}(\X;\X)}\|\Gamma_t\mathcal{I}\|_{\mathcal{L}(E;\X)}dt<\infty,
\]
where both $t\mapsto \|\Gamma_t\|_{\mathcal{L}(\X;\X)}$ and $t\mapsto \|\Gamma_t\mathcal{I}\|_{\mathcal{L}(E;\X)}$ are bounded from below functions in $(0,T)$. Further, we assume that there exists $\theta'<\theta$ such that 
\[
\int^T_0\|\Gamma_t\|^{1-\theta'}_{\mathcal{L}(\X;\X)}dt<\infty.
\]
\end{hyp1}
We note that by Hypotheses \ref{supercontrollabilità} we have
\begin{align}
C_T:=\int^T_0 K_tdt<\infty,\qquad \lim_{T\rightarrow 0}C_T=0\label{K-integrale},
\end{align}
where $K_t$ is the constant given by Corollary \ref{SchauderRV}.
\begin{proposition}\label{WPK}
    Assume that Hypotheses \ref{supercontrollabilità} holds true. Then, equation \eqref{Back-Kolmo} admits a unique solution $U\in C_b^{0,1}([0,T]\times \X;\X)$. Moreover the map $x\rightarrow DU(t,x)\mathcal{I}$ belongs to $C^1_b(\X;\mathcal{L}(E;\X))$ for every $t\in [0,T]$ and it holds
\begin{align}
\label{stima2}
\sup_{t\in[0,T]}\left(\norm{U(t,\cdot)}_{C^1_b(\X;\X)}+\|DU(t,\cdot)\mathcal{I}\|_{C^1_b(\X;\mathcal{L}(E,\X))}\right)&\leq C_{T}e^{C_T\|\mathcal{N}\|_{C^\theta_b(\X;\X)}}\|\mathcal{M}\|_{C^\theta_b(\X;\X)},
 \end{align}
where $C_T$ is the constant defined in \eqref{K-integrale}.
 
\end{proposition}
\begin{proof}
For every $\gamma\geq 0$, we denote by $\mathscr{E}_{T,\gamma}$ the subspace of $C_b^{0,1}([0,T]\times \X;\X)$ endowed with the norm

\begin{align*}
 &\|\Phi\|_{T,\gamma}:=\sup_{t\in [0,T]}e^{\gamma t}\norm{\Phi(t,\cdot)}_{2,\mathcal{I}},\\
&\norm{\varphi}_{2,\mathcal{I}}:=\norm{\varphi}_{C^1_b(\X;\X)}+\|D\varphi\mathcal{I}\|_{C^1_b(\X;\mathcal{L}(E,\X))}
\end{align*}

$(\mathscr{E}_{T,\gamma},\norm{\cdot}_{T,\gamma})$ is a Banach space.
We consider the operator $V$, defined for every $U\in\mathscr E_{T,\gamma}$ by
\begin{equation}\label{operatoreVolterra}
V(U)(t,x)=\int_t^T\mathcal{R}(r-t)\left( DU(r,\cdot)\mathcal{I}\mathcal{N}(\cdot)+\mathcal{M}(\cdot)\right)(x)dr,\quad t\in [0,T],\; x\in \X.
\end{equation}
We prove that a suitable choice of $\gamma$ implies that the operator $V$ is a contraction in $\mathscr E_{T,\gamma}$. 
For every $U\in \mathscr E_{T,\gamma}$, $t\in [0,T]$ and $x,h\in \X$, we have
\begin{align*}
\|DU(t,x+h)\mathcal{I}\mathcal{N}(x+h)-DU(t,x)\mathcal{I}\mathcal{N}(x)\|_\X &\leq \|DU(t,x+h)\mathcal{I}\mathcal{N}(x+h)-DU(t,x+h)\mathcal{I}\mathcal{N}(x)\|_\X\\
&+\|DU(t,x+h)\mathcal{I}\mathcal{N}(x)-DU(t,x)\mathcal{I}\mathcal{N}(x)\|_\X\\
&\leq \|DU(t,x+h)\mathcal{I}\|_{\mathcal{L}(E;\X)}\|\mathcal{N}(x+h)-\mathcal{N}(x)\|_E\\
&+\|DU(t,x+h)\mathcal{I}-DU(t,x)\mathcal{I}\|_{\mathcal{L}(E;X)}\|\mathcal{N}(x)\|_E.
\end{align*}
Hence, for every $t\in [0,T]$ we get
 \begin{align}
\|DU(t,\cdot)\mathcal{I}\mathcal{N}(\cdot)\|_{C^\theta_b(\X;\X)}&\leq \|DU(t,\cdot)\mathcal{I}\|_{C^\theta_b(\X;\mathcal{L}(E;\X))}\|\mathcal{N}\|_{C^\theta_b(\X;E)}\notag\\
&\leq \|DU(t,\cdot)\mathcal{I}\|_{C^1_{b}(\X;\mathcal{L}(E;\X))}\|\mathcal{N}\|_{C^\theta_b(\X;E)}\notag\\
&\leq\|U(t,\cdot)\|_{2,\mathcal{I}}\|\mathcal{N}\|_{C^\theta_b(\X;E)}.\label{stimaDU}
\end{align}
By Corollary \ref{SchauderRV} and \eqref{stimaDU}, for every $\gamma\geq 0$ and $t\in [0,T]$ we get
 \begin{align}
e^{\gamma t}\|V(U)(t,\cdot)\|_{2,\mathcal{I}}
\leq & \int_t^T e^{\gamma t}K_{r-t}\|DU(r,\cdot)\mathcal{I}\mathcal{N}(\cdot)\|_{C^{\theta}_b(\X;\X)}dr+ \|\mathcal{M}\|_{C^\theta_b(\X;\X)}\int_t^T e^{\gamma t}K_{r-t}dr\notag\\
\leq &\|\mathcal{N}\|_{C^{\theta}_b(\X;E)}\|U\|_{T,\gamma}\int_t^T e^{(t-r)\gamma}K_{r-t}dr+ \|\mathcal{M}\|_{C^\theta_b(\X;\X)}e^{\gamma T}\int_0^T K_{s}ds \notag \\
\leq & \|\mathcal{N}\|_{C^{\theta}_b(\X;E)}\|U\|_{T,\gamma}\int_0^T e^{-s\gamma}K_{s}dr+ \|\mathcal{M}\|_{C^\theta_b(\X;\X)}e^{\gamma T}\int_0^T K_{s}ds 
.\label{beta1}
 \end{align}
From \eqref{K-integrale} and \eqref{beta1}, it follows that $V( \mathscr{E}_{T,\gamma})\subseteq  \mathscr{E}_{T,\gamma}$ for every $\gamma\geq 0$.

Let $U_1,U_2\in \mathscr{E}_{T,\gamma}$. By Corollary \ref{SchauderRV}, the definition of $V$ and \eqref{stimaDU},  arguing as in \eqref{beta1}, for every $t\in [0,T]$ and $\gamma\geq 0$ we get
\begin{align*}
e^{\gamma t}\|V(U_1(t,\cdot))-V(U_2(t,\cdot))&\|_{2,\mathcal{I}}
\leq 
\|\mathcal{N}\|_{C^\theta_b(\X;E)}\|U_1-U_2\|_{T,\gamma}\int^T_0 e^{-\gamma s}K_{s}ds.
\end{align*}
By the dominated convergence theorem there exists $\gamma_0\geq 0$ such that $V$ is a contraction on  $\mathscr{E}_{T,\gamma_0}$ and so \eqref{Back-Kolmo} has a unique solution $U\in\mathscr E_{T,\gamma_0}$. Finally, from \eqref{K-integrale}, \eqref{beta1} and Gronwall's lemma we obtain \eqref{stima2}. 
\end{proof}

\subsection{Finite-dimensional Ornstein-Uhlenbeck operator}
In this section we assume that $\X=\R^n$ and we recall some classical results about the Ornstein-Uhlenbeck semigroup $\{R(t)\}_{t\geq 0}$ defined by \eqref{OUS}. 

We consider the (possible degenerate) Ornstein-Uhlenbeck operator defined by
\begin{equation*}
    L\varphi(x):=\frac{1}{2}{\rm Trace}\left[Q\nabla^2\varphi(x)\right]+\scal{\mathcal{C}x}{\varphi(x)},\qquad \varphi\in C^2_b(\R^n),\; x\in\R^n.
\end{equation*}
Moreover, for every $T>0$ and $f\in C_b([0,T]\times\R^n)$ we consider the backward parabolic equation 
\begin{align}\label{Parabolica}
\left\{
\begin{array}{ll}
\displaystyle \frac{\partial u(t,x)}{\partial t}+Lu(t,x)+f(t,x)=0, \qquad t\in(0,T], \ x\in \R^n, \vspace{3mm} \\
u(T,x)=0,\qquad  x\in \R^n.
\end{array}
\right.
\end{align}
\begin{thm}\label{OU-finito}
Let $T>0$, let $\theta\in (0,1)$ and let $f:[0,T]\times\R^n\rightarrow \R$ be a continuous function such that $f(t,\cdot)\in  C^{\theta}_b(\R^n)$ for every $t\in [0,T]$. Assume that Hypotheses \ref{supercontrollabilità} holds true. The
parabolic equation \eqref{Parabolica} has unique strong solution $u\in C^{0,1}_b([0,T]\times\R^n)$ given by
\[
u(t,x)=\int_{t}^TR(r-t)f(r,x)dr. \qquad (t,x)\in[0,T]\times \R^n.
\]
Moreover there exists a sequence $(f_h)_{h\in\N}\subseteq C^{0,2}_b([0,T]\times\R^n;\R)$ and $\theta'<\theta$ such that
\begin{enumerate}
\item for every $t\in [0,T]$ we have
\begin{equation}\label{conv-g}
\sup_{n\in\N}\sup_{t\in[0,T]}\|f_n(t,\cdot)\|_{C^\theta(\R^n)}<\infty, \qquad \lim_{n\rightarrow\infty}\norm{f_h(t,\cdot)-f(t,\cdot)}_{C^{\theta'}_b(\R^n)}=0;
\end{equation}
\item  for every $n\in\N$, the parabolic equation \eqref{Parabolica}, with $f$ replaced by $f_n$, has a unique strict solution $u_n\in C^{1,2}_b([0,T]\times\R^n)$ given by
\[
u_n(t,x)=\int_{t}^TR(r-t)f_n(r,x)dr;
\]
\item for every $t\in [0,T]$ we have
\begin{equation}\label{app-sol-par}
\lim_{n\rightarrow \infty}\norm{u_n(t,\cdot)-u(t,\cdot)}_{C^1_b(\R^n)}=0.
\end{equation}
\end{enumerate}
\end{thm}
\begin{proof}
All the statements are quite classical, we refer to \cite{Lun1997} for a proof and a detailed discussion. In \cite{Lun1997}, in a more general setting, the convergence result \eqref{app-sol-par} is proved with respect to the $C_b$-norm, so here we prove a finer result for this particular case. 

By standard approximation arguments ($f_n$ should be defined by means of convolution with mollifiers for every $n\in\N$), \eqref{conv-g} is verified for every $\sigma<\theta$. Moreover, by Hypotheses \ref{supercontrollabilità} there exists $\theta'<\theta$ such that 
\[
\int_0^T\norm{\Gamma_t}^{1-\theta'}dt<\infty,
\]
so we fix such a $\theta'$.
For every $t\in [0,T]$ and $n\in\N$, by Proposition \ref{SchauderR} there exists a constant $C_{\theta'}>0$ such that, for every $t\in[0,T]$, we get
\begin{align*}
    \norm{u_n(t,\cdot)-u(t,\cdot)}_{C^1_b(\R^n)}\leq C_{\theta'}\int_t^T\norm{\Gamma_{r-t}}^{1-\theta'}\norm{f_n(r,\cdot)-f(r,\cdot)}_{{C}^{\theta'}_b(\R^n)}dr.
\end{align*}
From \eqref{conv-g} and the dominated convergence theorem, we obtain \eqref{app-sol-par}.
\end{proof}

\end{appendices}

\textbf{Acknowledgments} 
The authors are members of GNAMPA (Gruppo Nazionale per l'Analisi Matematica,
la Probabilit\`a e le loro Applicazioni) of the Italian Istituto Nazionale di Alta Matematica (INdAM). D. Addona has been partially supported by the INdaMGNAMPA project 2023 ``Operatori ellittici vettoriali a coefficienti illimitati in spazi $L^p$'' CUP E53C22001930001. 
D. A. Bignamini has been partially supported by the INdaMGNAMPA
Project 2023 ``Equazioni differenziali stocastiche e operatori di Kolmogorov in dimensione infinita'' CUP E53C22001930001. 
The authors have no relevant financial or non-financial interests to disclose.
The authors are really grateful to Enrico Priola for the inspiring discussions about the stochastic damped wave equations.

\section*{Declarations} 

Data sharing not applicable to this article as no datasets were generated or analysed during the current study.

\end{document}